\documentclass[11pt]{article}
\usepackage[T1]{fontenc}
\usepackage{lmodern,amsmath,amssymb,amsthm,mathtools,microtype,mathrsfs}
\usepackage[margin=1in]{geometry}
\usepackage{enumitem,booktabs,needspace,graphicx,tikz,placeins,float}
\usetikzlibrary{arrows.meta,positioning,calc}
\usepackage{xcolor}
\usepackage[hidelinks]{hyperref}
\hypersetup{
  pdftitle={Elementwise Error, Risk, and Level-Set Geometry for Two-Dimensional Total Variation Denoising},
  pdfauthor={Teng Zhang},
  pdfsubject={Statistical theory of two-dimensional total variation denoising},
  pdfkeywords={total variation denoising, fused lasso, elementwise estimation, risk bounds, level-set geometry}
}
\allowdisplaybreaks[2]
\numberwithin{equation}{section}
\newtheorem{lemma}{Lemma}[section]
\newtheorem{theorem}[lemma]{Theorem}
\newtheorem{proposition}[lemma]{Proposition}
\newtheorem{corollary}[lemma]{Corollary}
\theoremstyle{remark}
\newtheorem{remark}[lemma]{Remark}
\newcommand{\R}{\mathbb R}
\newcommand{\Z}{\mathbb Z}
\newcommand{\E}{\mathbb E}
\newcommand{\PP}{\mathbb P}
\newcommand{\one}{\mathbf1}
\newcommand{\TV}{\operatorname{TV}}
\newcommand{\Var}{\operatorname{Var}}

\newcommand{\ip}[2]{\langle #1,#2\rangle}
\newcommand{\norm}[1]{\lVert #1\rVert}
\newcommand{\cutnorm}[1]{\mathcal{D}(#1)}
\newcommand{\bb}{\mathbf{b}}

\newcommand{\be}{\mathbf{e}}

\newcommand{\bu}{\mathbf{u}}
\newcommand{\bv}{\mathbf{v}}

\newcommand{\bx}{\mathbf{x}}
\newcommand{\by}{\mathbf{y}}
\newcommand{\bz}{\mathbf{z}}

\newcommand{\btheta}{\boldsymbol{\theta}}
\newcommand{\beps}{\boldsymbol{\varepsilon}}

\newcommand{\bI}{\mathbf{I}}

\newcommand{\bSigma}{\boldsymbol{\Sigma}}

\newcommand{\bd}{\mathbf{d}}
\newcommand{\bh}{\mathbf{h}}
\title{\Large\bfseries Elementwise Error, Risk, and Level-Set Geometry\\for Two-Dimensional Total Variation Denoising}
\author{Teng Zhang\\[3pt]
\small Dr.\ Yvette Kanouff School of Data,\\
\small Mathematical, and Statistical Sciences\\
\small University of Central Florida\\
\small Orlando, FL 32816, USA\\
\small\href{mailto:teng.zhang@ucf.edu}{\texttt{teng.zhang@ucf.edu}}}
\date{}

\begin{document}
\maketitle
\begin{abstract}
We study elementwise error, global risk, and the geometry of large errors for two-dimensional total variation (TV) denoising of piecewise-constant grid signals. Our elementwise error bounds have three components: regularization bias, side effects, and corner effects; the latter two depend on the horizontal and vertical distances to the boundary. Under independent \(N(0,\sigma^2)\) noise, the elementwise risk bounds are sharp up to logarithmic factors in the worst case over signal values.
Above an explicit threshold, the large-error vertices in each square where the signal is constant form at most four corner-monotone staircase regions, with a constant error sign on each region.

For a grid of \(N\) vertices partitioned into \(k\) rectangles of bounded aspect ratio, an explicit deterministic regularization parameter gives a mean squared error bound of order \(\sigma^2(k/N)^{3/4}\), up to logarithmic factors. A matching lower bound on regular square partitions establishes sharpness among all TV-penalized estimators. The analysis uses local TV problems with bounds uniform over exterior values. Its deterministic analysis tools are a nested-rectangle propagation argument and comparisons with corner-monotone sets, and its probabilistic analysis tool is a directional fluctuation parameter that controls the cumulative fluctuation of Gaussian noise over a grid and has expectation of order \(\sigma(\log N)^{3/4}\).
\end{abstract}


\section{Introduction}\label{sec:introduction}
Total variation denoising estimates a signal by penalizing differences
between neighboring grid vertices. For a piecewise-constant signal,
a global squared-error bound controls the average squared error
over the grid, but does not describe how the error varies within a
region on which the true signal is constant. In two dimensions, a
vertex can be close to one side of such a region, close to a corner,
or far from its boundary. We study how proximity to sides and
corners affects the estimation error and investigate the level-set geometry, i.e., the sets of
vertices whose error magnitudes exceed a specified threshold.

Consider the observation model
\begin{equation}\label{eq:intro-observation-model}
 y_p=\theta_p+\varepsilon_p,\qquad p\in[n]^2,
\end{equation}
where $n\ge4$, the noise 
$\varepsilon_p\stackrel{\mathrm{i.i.d.}}{\sim}N(0,\sigma^2)$,
$\sigma>0$, and $\btheta$ is piecewise constant on an unknown
partition of $[n]^2$ into axis-aligned rectangles. Let $N=n^2$ be the number
of observations. We study the total variation (TV)-penalized estimator
\begin{equation}\label{eq:whole-image-model}
 \widehat{\btheta}^{\lambda}
 =\arg\min_{\bx\in\R^{[n]^2}}
 \left\{\sum_{p\in[n]^2}(y_p-x_p)^2
 +\lambda\sum_{\substack{\{p,q\}\subseteq[n]^2\\p\sim q}}
 |x_p-x_q|\right\},
\end{equation}
where $p\sim q$ denotes adjacency in the grid graph, each unordered
edge is counted once, and $\lambda$ is the regularization parameter. Throughout the paper, $\lambda>0$.

\subsection{Main results}
\label{subsec:informal-results}
We first present two informal conclusions, ignoring universal
constants and logarithmic factors. The first bounds the error at each vertex. The notation
$\widetilde{\lesssim}$ denotes an inequality up to such factors.
Throughout this overview, ``with high probability'' means probability
at least $1-n^{-10}$. Write
$e_p^\lambda:=\widehat\theta_p^\lambda-\theta_p$ for the estimation
error at vertex $p$.

\begin{theorem}[Informal: elementwise error on piecewise-constant signals]
\label{thm:informal-piecewise}
Consider the model~\eqref{eq:intro-observation-model} and the total
variation estimator~\eqref{eq:whole-image-model}.
For each vertex $p$, let $m_p$ be the shorter side length of the
rectangle of the partition containing $p$, and let $r_p,s_p$
be the distances from $p$ to the nearer horizontal and vertical
sides of that rectangle, counting boundary vertices at distance one.
With high probability, simultaneously for all vertices $p$ and all
$\lambda\ge C\sigma (\log N)^{3/4}$,
\begin{equation}\label{eq:theorem_11}
 |e_p^\lambda|
 \;\widetilde{\lesssim}\;
 \frac{\lambda}{m_p}
 +\frac{\sigma}{\sqrt{r_p s_p}}
 +\frac{\sigma^2}{\lambda \min(r_p,s_p)}.
\end{equation}
\end{theorem}

For illustration, consider a regular partition into $k=q^2$
squares of side length $m$, where $q\ge3$, so that $N=km^2$.
For sufficiently large $m$, Table~\ref{tab:elementwise-rmse}
compares the worst-case elementwise root mean squared error
(RMSE) at representative vertices of an interior partition square
and the average MSE over that square, up to logarithmic factors.
Here the worst case is over all signals that are  constant on this partition.
\begin{table}[H]
\centering
\caption{Worst-case elementwise RMSE at representative vertices
and average MSE over an interior partition square of side length $m$,
up to logarithmic factors, for different choices of $\lambda$.}
\label{tab:elementwise-rmse}
\begin{tabular}{@{}lccc@{}}
\toprule
Location
& $\lambda=\sigma m^{1/8}$
& $\lambda=\sigma m^{1/4}$
& $\lambda=\sigma m^{1/2}$ \\
\midrule
Interior: $r_p,s_p\asymp m$
& $\sigma m^{-7/8}$
& $\sigma m^{-3/4}$
& $\sigma m^{-1/2}$ \\

Intermediate: $r_p,s_p\asymp m^{1/2}$
& $\sigma m^{-1/2}$
& $\sigma m^{-1/2}$
& $\sigma m^{-1/2}$ \\

Side: $\min(r_p,s_p)=1$, $\max\{r_p,s_p\}\asymp m$
& $\sigma m^{-1/8}$
& $\sigma m^{-1/4}$
& $\sigma m^{-1/2}$ \\

Corner: $r_p=s_p=1$
& $\sigma$
& $\sigma$
& $\sigma$ \\
\midrule
Average MSE over the square
& $\sigma^2 m^{-5/4}$
& $\sigma^2 m^{-3/2}$
& $\sigma^2 m^{-1}$ \\
\bottomrule\end{tabular}
\end{table}
The high-probability bounds in Theorem~\ref{thm:informal-piecewise} and Table~\ref{tab:elementwise-rmse} 
follow from Corollary~\ref{cor:gaussian-signal-elementwise}.
We note that the upper bounds in expectation follow from
Corollary~\ref{cor:elementwise-risk} and the restriction in
Lemma~\ref{lem:patch-restriction}. 
The matching lower bounds follow
from the analysis of local estimator in Theorem~\ref{thm:elementwise-risk-lower} and
Corollary~\ref{cor:power-match}, and adapted to the global estimator
through the checkerboard example in
Section~\ref{subsec:gaussian-rectangle-results}.

The three terms in the upper bound of \eqref{eq:theorem_11} describe different effects.
The bias term $\lambda/m_p$ bounds error caused by TV-regularization across
the rectangle's boundary; this error can remain even when the
observations inside the rectangle are noiseless. The corner term
$\sigma/\sqrt{r_p s_p}$ depends on both vertical and horizontal distances of $p$ to the boundary,
whereas the side term $\sigma^2/(\lambda \min(r_p,s_p))$ depends
on the smaller distance. Increasing $\lambda$ reduces the side
term and increases the bias term.
Corollary~\ref{cor:gaussian-signal-elementwise} gives the formal
statement, including logarithmic factors and the probability bound.
Theorem~\ref{thm:elementwise-risk-lower} establishes sharpness of the
combined three-term elementwise risk bound, up to logarithmic
factors, at every vertex in its stated range of $\lambda$.

The argument uses the local TV problem~\eqref{eq:model} and its
analysis in Theorem~\ref{thm:elementwise}.
Lemma~\ref{lem:patch-restriction} identifies the restriction of the
global estimator~\eqref{eq:whole-image-model} to any square with a
local TV problem. Since the bound in
Theorem~\ref{thm:elementwise} holds uniformly over exterior values,
it applies to the global estimator.

The corresponding expected-risk bounds, combined with square
coverings and summation, yield a bound on the average squared error
over the full grid. Suppose the partition has $k$ rectangles with
aspect ratios at most $\rho\ge1$.
For a sufficiently large fixed universal constant $C_{\rm reg}$,
the deterministic choice
\begin{equation}\label{eq:intro-regularization}
 \lambda=C_{\rm reg}\sigma (\log N)^{5/4}(N/k)^{1/8}
\end{equation}
gives
\[
 \frac{1}{N}\E\norm{\widehat{\btheta}^{\lambda}-\btheta}_2^2
 \;\widetilde{\lesssim}\;
 \rho\sigma^2(k/N)^{3/4}.
\]
Corollary~\ref{cor:rectangle-gaussian} gives the logarithmic factors
and also treats general aspect ratios and other admissible
deterministic values of $\lambda$. The choice of $\lambda$
in~\eqref{eq:intro-regularization} balances the bias and side
contributions. It uses $k$ and $\sigma$, but not the signal's total
variation, rectangle locations, or signal values.

Corollary~\ref{thm:global-risk-sharpness} gives a matching lower
bound on regular square partitions, showing that no deterministic
choice of $\lambda$ improves the power $(k/N)^{3/4}$ in the worst
case. For fixed $k$, the $N^{-3/4}$ exponent also appears in the
small-noise lower bound for constrained TV of Chatterjee and
Goswami~\cite[Theorem~II.3]{CGTV};
Section~\ref{subsec:two-dimensional-risk-comparison} discusses its
connection to penalized TV.

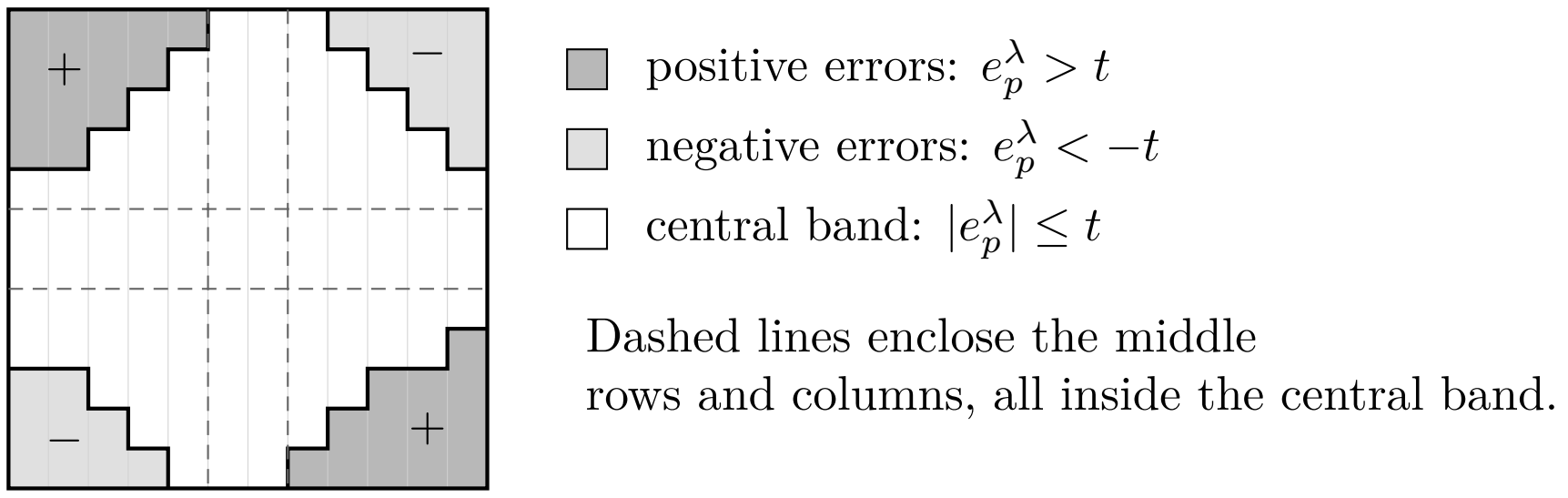
\begin{figure}[!htbp]
\centering
\begin{tikzpicture}[x=0.29cm,y=-0.29cm,font=\small]
 \fill[white] (0,0) rectangle (12,12);
 \foreach \r/\a in {0/5,1/4,2/3,3/2}
   \fill[black!28] (0,\r) rectangle (\a,\r+1);
 \foreach \r/\a in {0/4,1/3,2/2,3/1}
   \fill[black!12] (12-\a,\r) rectangle (12,\r+1);
 \foreach \r/\a in {11/4,10/3,9/2}
   \fill[black!12] (0,\r) rectangle (\a,\r+1);
 \foreach \r/\a in {11/5,10/4,9/3,8/1}
   \fill[black!28] (12-\a,\r) rectangle (12,\r+1);
 \draw[step=1,black!18,line width=0.15pt] (0,0) grid (12,12);
 \draw[thick] (0,0) rectangle (12,12);
 \draw[thick] (5,0)--(5,1)--(4,1)--(4,2)--(3,2)--(3,3)--(2,3)--(2,4)--(0,4);
 \draw[thick] (8,0)--(8,1)--(9,1)--(9,2)--(10,2)--(10,3)--(11,3)--(11,4)--(12,4);
 \draw[thick] (0,9)--(2,9)--(2,10)--(3,10)--(3,11)--(4,11)--(4,12);
 \draw[thick] (12,8)--(11,8)--(11,9)--(9,9)--(9,10)--(8,10)--(8,11)--(7,11)--(7,12);
 \draw[densely dashed,black!60] (5,0)--(5,12) (7,0)--(7,12)
    (0,5)--(12,5) (0,7)--(12,7);
 \node at (1.4,1.5) {$+$};
 \node at (10.5,1.1) {$-$};
 \node at (1.4,10.8) {$-$};
 \node at (10.5,10.5) {$+$};
 \fill[black!28] (14,1) rectangle (15,2);
 \draw (14,1) rectangle (15,2);
 \node[anchor=west] at (15.5,1.5) {positive errors: $e_p^\lambda>t$};
 \fill[black!12] (14,3) rectangle (15,4);
 \draw (14,3) rectangle (15,4);
 \node[anchor=west] at (15.5,3.5) {negative errors: $e_p^\lambda<-t$};
 \draw (14,5) rectangle (15,6);
 \node[anchor=west] at (15.5,5.5) {central band: $|e_p^\lambda|\le t$};
 \node[anchor=west,align=left] at (14,9)
   {Dashed lines enclose the middle\\rows and columns, all inside the central band.};
\end{tikzpicture}
\caption{Illustration of Theorem~\ref{thm:informal-piecewise-geometry}
on the square $Q$. The shaded vertices have $e_p^\lambda>t$ or $e_p^\lambda<-t$,
as indicated by their signs. Each shaded region is corner-monotone,
and a corner region may be empty. The unshaded vertices form the
central band $Q\setminus E_t(Q)$, which contains the middle rows
and columns indicated by the dashed lines.}
\label{fig:corner-structure}
\end{figure}

Our second conclusion describes level-set geometry, i.e., the shape of large errors on a
square where the true signal is constant.

\begin{theorem}[Informal: level-set geometry of the estimation error]
\label{thm:informal-piecewise-geometry}
Under the assumptions of Theorem~\ref{thm:informal-piecewise},
let $Q$ be an axis-aligned $m\times m$ square contained in one
rectangle of the partition.
With high probability, simultaneously for all 
$\lambda\ge C\sigma (\log N)^{3/4}$, and all levels satisfying
\begin{equation}\label{eq:intro-geometry-threshold}
 t\;\widetilde{\gtrsim}\;
 \frac{\lambda}{m}+\frac{\sigma^2}{\lambda},
\end{equation}
the large-error level set
\[
 E_t(Q)=\{p\in Q:|e_p^\lambda|>t\}
\]
is a union of at most four corner-monotone regions, at most one
attached to each corner of $Q$, and $e_p^\lambda$ has a constant sign on
each nonempty region.
\end{theorem}

A set is \emph{corner-monotone} from a corner if, whenever it contains
a vertex, it also contains every grid vertex in the rectangle between
that vertex and the corner. Thus the vertices whose errors exceed
the stated threshold in magnitude lie in staircase regions attached
to the corners of $Q$. The remaining vertices, whose errors have magnitude at most $t$,
form the central band $Q\setminus E_t(Q)$. This band is nonempty
and connected, contains the middle row(s) and column(s) of $Q$,
and intersects each row or column in an interval.
Figure~\ref{fig:corner-structure} illustrates these conclusions. Corollary~\ref{cor:gaussian-signal-geometry} gives the formal
geometry statement derived from Theorem~\ref{thm:geometry}.

Theorem~\ref{thm:informal-piecewise} applies to every rectangle of
the partition, while Theorem~\ref{thm:informal-piecewise-geometry}
applies to any square $Q$ contained in one such rectangle; $Q$ may
coincide with that rectangle when it is square.
Taking $\delta=n^{-10}$ in
Corollaries~\ref{cor:gaussian-signal-elementwise}
and~\ref{cor:gaussian-signal-geometry} gives the probability stated in
Theorems~\ref{thm:informal-piecewise}
and~\ref{thm:informal-piecewise-geometry}.

\subsection{Related work}\label{subsec:intro-related-work}
\paragraph{One-dimensional and graph TV.}
The variational denoising model of Rudin, Osher and
Fatemi~\cite{ROF} is a foundational example of total variation
regularization. Statistical developments include the locally adaptive regression
splines of Mammen and van de Geer~\cite{MvG} and the fused lasso
of Tibshirani et al.~\cite{FusedLasso}.
For one-dimensional TV denoising, Dalalyan, Hebiri and
Lederer~\cite{DHL} derive oracle inequalities that account for the
spacing of the signal's changes. Guntuboyina, Lieu, Chatterjee and
Sen~\cite{GLCS} establish adaptive risk bounds for one-dimensional
TV denoising and higher-order trend filtering.

Zhang and Chatterjee~\cite{ZC} develop elementwise error bounds
for one-dimensional generalized fused lasso under a convex loss.
For squared loss, their bound depends on the distance to a change
point, the length of the containing interval on which the signal is
constant, and the regularization parameter.
Chatterjee~\cite{ChatterjeeMinmax} derives elementwise minmax/maxmin
representations for one-dimensional TV denoising using corrected
local averages over nested intervals. Chatterjee and
Ghoshal~\cite{CGelementwise} establish such representations on
arbitrary finite graphs, using nested connected vertex sets and
signed boundary corrections.

On general graphs, Wang et al.~\cite{WSST} develop trend filtering
through penalties on graph differences, with graph TV as the
order-zero case. Madrid Padilla et al.~\cite{DFS} use a depth-first
search ordering to compare graph variation with variation on a
chain, obtaining error bounds for both the graph fused lasso and a
computationally simpler chain estimator. Fan and Guan~\cite{FG}
study an $\ell_0$ penalty that counts edges across signal changes.
Statistical guarantees under broader noise assumptions are
developed by Madrid Padilla~\cite{MPvariance} for graph fused lasso
and by Madrid Padilla and Chatterjee~\cite{MPquantile} for quantile
trend filtering and TV.

\paragraph{Two-dimensional risk and rectangular adaptation.}
For penalized two-dimensional TV denoising, H\"utter and
Rigollet~\cite[Corollary~5]{HR} establish sharp oracle inequalities
that account for edges across signal changes.
With the constraint radius set to the true total variation,
Chatterjee and Goswami~\cite[Theorem~II.2]{CGTV} establish an
adaptive upper risk bound for constrained TV on rectangular
partitions, giving the $N^{-3/4}$ rate up to logarithmic factors
when the number of rectangles and their aspect ratios are bounded. Their Theorem~II.3 gives a matching $N^{-3/4}$ lower
bound in the small-noise limit at a fixed signal with two equal
rectangular regions, using tangent-cone geometry and Gaussian-width
bounds. Their upper bound concerns the constrained estimator, and
\cite[Section~III-F]{CGTV} poses the corresponding question for the
penalized estimator: ``A possible next step
is to investigate whether a similar $N^{-3/4}$ rate is attained by
the penalized TVD estimator and if so, for what range of the
tuning parameter $\lambda$''. 
The written version of van de Geer's 2019 Oberwolfach
Lecture~\cite[Section~5]{vdGTV} announces the same exponent for
penalized grid TV in a two-rectangle example and sketches a
path-based argument. 
Our result gives expected-risk bounds for general rectangular
partitions of bounded aspect ratio, with an explicit deterministic
regularization rule and matching dependence on $k/N$ on regular
square partitions. 

Other procedures also exploit rectangular structure.
Chatterjee and Goswami~\cite{CGCART} show that optimal regression
trees attain expected squared error per observation bounded by
$C\sigma^2 k\log N/N$ for two-dimensional signals constant on $k$
axis-aligned rectangles. Ortelli and van de Geer~\cite{OvGmixed}
obtain adaptive oracle inequalities for a two-dimensional penalty
built from mixed second differences.

Sadhanala, Wang and Tibshirani~\cite{SWT} study minimax rates over
discrete bounded-variation classes on higher-dimensional grids and
show that TV denoising is rate-optimal up to logarithmic factors.
Sadhanala, Wang, Hu and Tibshirani~\cite{SWHT} place these results
in the framework of Kronecker trend filtering for all derivative
orders and dimensions. At order zero in two dimensions, their
estimator is the TV estimator in~\eqref{eq:whole-image-model},
up to loss normalization. For a total-variation radius of order
$\sqrt N$ and a fixed noise level, the minimax rate is $N^{-1/2}$
up to logarithmic factors.
In the continuum white-noise model,
del \'Alamo, Li and Munk~\cite{ALM} combine TV regularization with
multiscale data-fidelity constraints and prove minimax $L^q$-risk
bounds, up to logarithmic factors, over bounded-variation classes
with uniform supremum-norm control.

\paragraph{Geometry of TV estimates.}
Geometric properties of TV estimators have also been studied in
continuum models. Caselles, Chambolle and Novaga~\cite{CCN} study
discontinuity sets of TV minimizers; Chambolle, Duval, Peyr\'e and
Poon~\cite{CDPP} study extended-support localization; and
Esedo\u{g}lu and Osher~\cite{EO} analyze anisotropic ROF geometry.
For rectilinear anisotropy, \L asica, Moll and Mucha~\cite{LMM}
show that denoising and the associated TV flow preserve piecewise
constancy when the observed datum is constant on a rectangular grid.
The geometry studied here concerns estimation errors inside a
region on which the true signal is constant, even though the
observations fluctuate within that region.

The closest quantitative comparisons are with Zhang and
Chatterjee~\cite{ZC} for elementwise error and with H\"utter and
Rigollet~\cite{HR} and Chatterjee and Goswami~\cite{CGTV} for
two-dimensional risk. Sections~\ref{subsec:one-dimensional-comparison}
and~\ref{subsec:two-dimensional-risk-comparison} develop these
comparisons, including the spatial dependence of the error bounds,
the regularization parameters, and the signal assumptions.

\subsection{Contributions and proof ideas}
\label{subsec:intro-contributions}

Existing work on two-dimensional TV-penalized estimators has
focused primarily on global risk bounds. Our analysis bounds errors
at individual vertices and characterizes the level-set geometry of
large errors. Under Gaussian noise, the local elementwise risk bound
is sharp up to logarithmic factors at every vertex in the worst case
over exterior values. Summing the local bounds gives,
to our knowledge, the first complete proof of an expected-risk upper
bound of order $\sigma^2(k/N)^{3/4}$, up to logarithmic factors,
for the standard penalized two-dimensional TV estimator on general
rectangular partitions with bounded aspect ratios. For fixed $k$,
this gives the $N^{-3/4}$ rate. The explicit deterministic choice
of $\lambda$ depends only on $N$, $k$, and $\sigma$, without
requiring the signal's total variation, rectangle locations, or
signal values. A matching lower bound
on regular square partitions shows that this dependence on $k/N$
cannot be improved by a deterministic choice of the TV
regularization parameter.
The analysis is based on the local TV problem~\eqref{eq:model}:
localization transfers the elementwise and geometric results to the
global estimator, while square coverings and summation yield the
global risk bounds.

The technical components are as follows. For the elementwise
bound, optimality inequalities on nested rectangles limit how far a
high-threshold set can extend into the square. Separate control in
the two coordinate directions gives the corner and side terms
in~\eqref{eq:theorem_11}. For the level-set geometry result,
comparison arguments show that each threshold region is
corner-monotone and has a constant sign. For the Gaussian bound,
lattice stretching balances the two directional perimeters, after
which the Ding--Wirth estimate~\cite{DW} and Gaussian concentration
apply. These arguments may also be useful in studying other graph
structures, loss functions, or higher-dimensional settings.

\subsection{Outline}\label{subsec:intro-outline}
Section~\ref{sec:fluctuation} introduces the local TV problem
used to analyze restrictions of the global estimator, defines the
directional fluctuation parameter, and establishes expectation and
high-probability bounds under Gaussian noise.
Section~\ref{sec:main-results} presents the main deterministic
results, including the local elementwise and geometric bounds and
the global deterministic error bound.
Section~\ref{sec:gaussian-consequences} gives high-probability
guarantees for the global estimator, risk upper bounds for the local
and global estimators, and sharpness results under Gaussian noise.
Section~\ref{sec:examples-comparisons} compares our conclusions with one-dimensional elementwise
bounds and existing two-dimensional risk guarantees.
Section~\ref{sec:discussion} discusses the methodological
contributions, the scope of the risk results, and open questions.
Appendix~\ref{sec:preliminaries} collects shared deterministic tools,
while Appendices~\ref{sec:gaussian-flows}--\ref{sec:completion}
contain the proofs of the main results.

\section{Local TV problem and the directional fluctuation parameter}\label{sec:fluctuation}
Section~\ref{sec:estimator} defines the local TV problem and
relates it to the global estimator~\eqref{eq:whole-image-model}.
Section~\ref{subsec:directional-parameter} introduces the directional
fluctuation parameter, which quantifies noise fluctuations on the grid.
Theorem~\ref{thm:gaussian-D} bounds its expectation and gives a
high-probability bound under i.i.d.\ $N(0,\sigma^2)$ noise.

\subsection{Local TV problem}\label{sec:estimator}
We first introduce the local TV problem, which is used to analyze restrictions of the global estimator. The local TV problem is defined on an $m\times m$
square grid. Its interior grid, one-layer extension, and exterior frame are
\[
V_m=\{1,\ldots,m\}^2,\qquad
\overline V_m=\{0,\ldots,m+1\}^2,\qquad
F_m=\overline V_m\setminus V_m.
\]
In the grid graph, each vertex \((i,j)\in V_m\) has four neighbors $(i-1,j),\quad (i+1,j),\quad (i,j-1),\quad (i,j+1)$ in the extended grid.
Let $E_{\rm int}$ be the set of edges with both endpoints in $V_m$,
and let $E_{\rm fr}$ be the set of edges with one endpoint in $V_m$
and the other in $F_m$. Set $E=E_{\rm int}\cup E_{\rm fr}$.

Given observations \(\by=(y_v)_{v\in V_m}\) and exterior values \(\bb\in\mathbb R^{F_m}\), we define the local estimator by
\begin{equation}\label{eq:model}
\widehat{\bx}^{\lambda,\bb}
=\arg\min_{\bx\in\mathbb R^{V_m}}
\left\{
\|\by-\bx\|_2^2
+\lambda\sum_{\{v,w\}\in E_{\rm int}}|x_v-x_w|
+\lambda\sum_{\substack{\{v,w\}\in E_{\rm fr},\, v\in V_m}}
|x_v-b_w|
\right\}.
\end{equation}

The local TV problem is translation equivariant: subtracting the
same constant from the observations and every exterior value
subtracts that constant from the solution. Hence, if the signal is
constant with value $c_Q$ on $Q$, then
$\widehat{\btheta}^{\lambda}|_Q-c_Q\one_Q$ is the local estimator
with observations $(\by-\btheta)|_Q$ and exterior values
$b_w^Q-c_Q$. This reduction is deterministic and requires no
distributional assumption.

The following lemma identifies its connection to the global estimator~\eqref{eq:whole-image-model}.
\begin{lemma}[Restriction to a square]\label{lem:patch-restriction}
Let $\lambda>0$, and let $Q\subseteq[n]^2$ be an axis-aligned
square of side length $m$. After translating $Q$ to $V_m$,
the restriction $\widehat{\btheta}^{\lambda}|_Q$ is the solution to the local
TV problem~\eqref{eq:model} with observations $\by|_Q$ and suitable
exterior values. In the original grid coordinates, for each exterior
edge $\{v,w\}$ with $v\in Q$ and $w\notin Q$, assign
\[
b_w^Q=
\begin{cases}
\widehat\theta_w^\lambda,
    & w\in[n]^2\setminus Q,\\[2mm]
\widehat\theta_v^\lambda,
    & w\notin[n]^2.
\end{cases}
\]
Unused exterior corner vertices may be assigned arbitrarily.
\end{lemma}

The proof is given in Appendix~\ref{subsec:rectangular-auxiliary-proofs}.

The exterior values are inherited from the global estimator and may
therefore depend on both the observations and $\lambda$. Bounds
that hold simultaneously for all exterior values remain valid for
these data-dependent choices.

Throughout the paper, $i$ is the row index and $j$ the column index.
For a vertex $p=(i,j)\in V_m$, set
\begin{equation}\label{eq:distances}
 r_p=\min\{i,m+1-i\},\qquad
 s_p=\min\{j,m+1-j\},\qquad L_u=\log(eu)\quad(u\ge1).
\end{equation}
The quantities $r_p$ and $s_p$ are the distances from $p$ to
the nearer horizontal and vertical sides of the square, counting
boundary vertices at distance one; $\min(r_p,s_p)$ is the distance to the
nearest side. The geometric mean $\sqrt{r_p s_p}$ records proximity
to both sides meeting at the corresponding corner.
Sometimes we write $r_Q(p),s_Q(p)$ for the two directional distances within the rectangle $Q$.

\subsection{Directional perimeters and fluctuation control}\label{subsec:directional-parameter}
Let $R\subseteq\Z^2$ be a nonempty finite axis-aligned rectangle, and let
$\be_1=(1,0)$ and $\be_2=(0,1)$ denote the two coordinate directions.
For every $S\subseteq R$, including the empty set, define
\[
P_k(S)
=
\#\Bigl\{
\{v,w\}:
v\in S,\;
w\in\Z^2\setminus S,\;
w-v\in\{\be_k,-\be_k\}
\Bigr\},
\qquad k=1,2.
\]
Thus, $P_k(S)$ counts all edges leaving $S$ in the $k$-th
coordinate direction, including edges leaving $R$. For $R=V_m$,
these include the edges to the frame $F_m$. The ordinary perimeter
$\operatorname{Per}(S)=P_1(S)+P_2(S)$ therefore agrees with the
edge boundary of $S$ in $\Z^2$. For example, if $S$ is a rectangle
with $r$ rows and $s$ columns, then $P_1(S)=2s$ and $P_2(S)=2r$.

For an array $\by$ on $R$, define its directional fluctuation
parameter by
\begin{equation}\label{eq:diamond-norm}
\cutnorm{\by}
=
\max_{\varnothing\ne S\subseteq R}
\frac{\left|\sum_{v\in S}y_v\right|}{\sqrt{P_1(S)P_2(S)}}.
\end{equation}
We use $R=V_m$ for the local problem and $R=[n]^2$ for the global
problem.
Geometrically, $\cutnorm{\by}$ measures the largest cumulative
fluctuation over a subset $S$, normalized by its boundary exposure
in the two coordinate directions. The numerator records the total
imbalance of the observations on $S$, while $P_1(S)$ and $P_2(S)$
measure its interaction with its complement in each direction.
Their geometric mean treats the directions symmetrically and makes
the parameter sensitive to the shape of $S$. For a rectangle $S$
with $r$ rows and $s$ columns, the denominator is $2\sqrt{rs}$,
so the ratio associated with this set is
\[
 \frac{\left|\sum_{v\in S}y_v\right|
}{2\sqrt{rs}}.
\]

For independent observations $y_v\sim N(0,\sigma^2)$, the bounds of
Ding and Wirth~\cite{DW} give the ordinary-perimeter estimate
\[
\E \max_{\varnothing\ne S\subseteq V_m}
\frac{\left|\sum_{v\in S}y_v\right|}{P_1(S)+P_2(S)}
=O(\sigma(\log m)^{3/4}).
\]
Lemma~\ref{lem:dw-perimeter} explains how this estimate is adapted from their
lattice-animal formulation. The following theorem gives the same
order for $\cutnorm{\by}$, whose denominator is
$\sqrt{P_1(S)P_2(S)}$.

Compared with the Ding--Wirth statistic, \(\cutnorm{\by}\) retains the horizontal and vertical boundary contributions separately through the geometric mean \(\sqrt{P_1(S)P_2(S)}\). When \(P_1(S)\) and \(P_2(S)\) are highly unbalanced, this quantity can be much smaller than the ordinary perimeter \(P_1(S)+P_2(S)\), so a bound for the ordinary-perimeter statistic does not by itself control \(\cutnorm{\by}\). This directional normalization is also what allows the elementwise bound in Theorem~\ref{thm:elementwise} to depend separately on the horizontal and vertical distances \(r_p\) and \(s_p\).

\begin{theorem}[The directional fluctuation parameter under Gaussian noise]
\label{thm:gaussian-D}
There are universal constants $c,C_D>0$ such that the following holds.
Let $m\ge4$ and $\sigma>0$, and let $y_v$, $v\in V_m$, be independent
$N(0,\sigma^2)$ variables. Then
\begin{equation}\label{eq:gaussian-D-mean}
 c\sigma L_m^{3/4}
 \le\E\cutnorm{\by}
 \le C_D\sigma L_m^{3/4},\,\,\text{where  $L_m=\log(em)$.}
\end{equation}
Moreover, for every $0<\delta<1$,
\begin{equation}\label{eq:gaussian-D}
 \PP\{\cutnorm{\by}\le\sigma\kappa_{m,\delta}\}\ge1-\delta,
 \qquad
 \kappa_{u,\delta}
 :=C_D L_u^{3/4}+\sqrt{\frac{1}{2}\log\frac{1}{\delta}},\qquad u\ge1.
\end{equation}
\end{theorem}

\noindent\textbf{Proof idea.}
Stretching the lattice balances the two directional perimeters and
reduces the upper expectation bound to the ordinary-perimeter
estimate of Ding and Wirth. The lower bound follows from the
arithmetic--geometric mean inequality, and Lipschitz continuity
gives the high-probability bound by Gaussian concentration.
See Appendix~\ref{sec:gaussian-flows}.

\begin{remark}[Correlated and heteroscedastic Gaussian noise]
\label{rem:gaussian-covariance}
The upper bounds in Theorem~\ref{thm:gaussian-D} also hold for
$m\ge4$ when the observation vector in the local TV problem satisfies $\by\sim N(0,\bSigma)$, where
$\bSigma\preceq\sigma^2\bI$ and $\sigma>0$. Here $\bSigma$ is the
joint covariance matrix; the coordinates may be dependent and
have unequal variances. To see this, let
$\bu$ be an independent centered Gaussian vector with covariance
$\sigma^2\bI-\bSigma$. Convexity and conditional Jensen's inequality
give
\[
 \E\cutnorm{\by}
 \le \E\cutnorm{\by+\bu}
 \le C_D\sigma L_m^{3/4}.
\]
Writing $\by=\bSigma^{1/2}\mathbf Z$, with independent $N(0,1)$
entries in $\mathbf Z$, the map
$\mathbf Z\mapsto\cutnorm{\bSigma^{1/2}\mathbf Z}$ is
$\sigma/2$-Lipschitz by~\eqref{eq:norm-lipschitz}.
Gaussian concentration~\cite{Vershynin} therefore yields
\[
 \PP\{\cutnorm{\by}\le\sigma\kappa_{m,\delta}\}\ge1-\delta,
 \qquad 0<\delta<1.
\]
One may choose any $\sigma>0$ with
$\sigma^2\ge\|\bSigma\|_{\rm op}$.
\end{remark}

\section{Main deterministic results}\label{sec:main-results}
We first study the local TV problem~\eqref{eq:model}.
Section~\ref{subsec:main-cut-condition} bounds the magnitude of the
estimator at each vertex, and Section~\ref{subsec:main-geometry}
describes the sets where that magnitude exceeds a threshold.
Section~\ref{subsec:main-piecewise} applies these results to the
global estimator~\eqref{eq:whole-image-model} for signals constant
on rectangular partitions. Throughout this section, the observations $\by$ 
are deterministic, and no probabilistic noise model is assumed.
The local bounds depend on $\by$ through $\cutnorm{\by}$, defined
in~\eqref{eq:diamond-norm}; the global bound uses the same
parameter applied to the errors $\beps=\by-\btheta$.

\subsection{Elementwise and averaged bounds for the local estimator}\label{subsec:main-cut-condition}
\begin{theorem}[Elementwise bound for the local estimator]\label{thm:elementwise}
Consider the local TV problem~\eqref{eq:model}.
For every $m\ge4$ and $\by\in\mathbb R^{V_m}$, simultaneously
for all $p\in V_m$, all exterior values $\bb\in\mathbb R^{F_m}$, and all
$\lambda\ge\cutnorm{\by}$,
\begin{equation}\label{eq:det-main}
 \left|\widehat x_p^{\lambda,\bb}\right|
 \le \frac{4\lambda}{m}
 +\min\left\{
 2\cutnorm{\by},\,
 256L_m\left(
 \frac{\cutnorm{\by}}{\sqrt{r_p s_p}}
 +\frac{\cutnorm{\by}^2}{\lambda \min(r_p,s_p)}
 \right)\right\}.
\end{equation}
\end{theorem}
The proof is given in Appendix~\ref{subsec:det-main-proof}.

\noindent\textbf{Interpretation of the bound.}
The term $\lambda/m$ describes the effect of regularization and the
exterior values,
which can persist even when $\by=0$. The corner term
$\cutnorm{\by}/\sqrt{r_p s_p}$ depends on both horizontal and
vertical distances to the boundary, whereas the side term
$\cutnorm{\by}^2/(\lambda \min(r_p,s_p))$ depends on the smaller
distance. The minimum with $2\cutnorm{\by}$ bounds their combined
contribution uniformly over the square.
Since the bound is uniform over all exterior values, it also applies to data-dependent exterior values, such as those inherited from the global estimator.

Squaring and summing the elementwise bound gives the following
bound on the averaged squared norm of the local estimator.

\begin{corollary}[Local averaged squared norm]\label{cor:det-global}
Consider the local estimator~\eqref{eq:model}.
Let $m\ge4$, and fix $\by\in\mathbb R^{V_m}$.
For every $\lambda\ge\cutnorm{\by}$,
\begin{equation}\label{eq:det-global}
 \frac{1}{m^2}\sup_{\bb\in\R^{F_m}}
 \norm{\widehat{\bx}^{\lambda,\bb}}_2^2
 \le C\left\{
 \frac{\lambda^2}{m^2}
 +\frac{L_m^4\cutnorm{\by}^2}{m^2}
 +\frac{L_m^2\cutnorm{\by}^4}{m\lambda^2}\right\},
\end{equation}
where $C>0$ is universal.
\end{corollary}
Appendix~\ref{subsec:det-summation} gives the summation argument.
Section~\ref{subsec:main-piecewise} applies this bound to the global
estimator for signals constant on rectangular partitions, by covering
each rectangle with squares.

\subsection{Level-set geometry for the local estimator}\label{subsec:main-geometry}
We describe the level-set geometry of the local estimator~\eqref{eq:model}.
The following definition captures the geometric structure used below.

Partition each coordinate range into
\[
I_-=\{1,\ldots,\lfloor m/2\rfloor\},
\qquad
I_+=\{\lfloor m/2\rfloor+1,\ldots,m\}.
\]
Their Cartesian products form four quadrants of \(V_m\), each containing one of the corners
\[
(1,1),\qquad (1,m),\qquad (m,1),\qquad (m,m).
\]
For a corner \(c=(c_1,c_2)\in\{1,m\}^2\), a set \(S\subseteq V_m\) is
\emph{corner-monotone from \(c\)} if, whenever \(p=(i,j)\in S\),
every grid point in the rectangle between \(p\) and \(c\) also belongs to \(S\); that is,
\[
\bigl([\min\{i,c_1\},\max\{i,c_1\}]
\times
[\min\{j,c_2\},\max\{j,c_2\}]\bigr)\cap V_m
\subseteq S.
\]
Equivalently, membership in \(S\) is preserved under every coordinate step toward \(c\).

\begin{theorem}[Corner structure of level sets for the local estimator]\label{thm:geometry}
For the local estimator~\eqref{eq:model}, there is a universal
constant $C>0$ such that the following holds for every $m\ge4$. Suppose
\begin{equation}\label{eq:geo-corner-condition}
 \lambda\ge C\cutnorm{\by},\qquad
 t\ge C\left\{\frac{\lambda}{m}
 +\frac{L_m^2\cutnorm{\by}^2}{\lambda}\right\}.
\end{equation}
Then, for every choice of exterior values $\bb$, the set
\[
\{p\in V_m:|\widehat x_p^{\lambda,\bb}|>t\}
\]
is the union of at most four sets, with at most one in each quadrant.
Each nonempty set is corner-monotone from the corner of its quadrant,
and $\widehat{\bx}^{\lambda,\bb}$ has constant sign on that set.

The central band
\[
B_t=\{p\in V_m:|\widehat x_p^{\lambda,\bb}|\le t\}
\]
is nonempty and connected, contains the middle row(s) and column(s),
and intersects every row and column in an interval.
\end{theorem}

Estimator values whose magnitudes exceed $t$ are confined to
staircase regions attached to the corners, with a constant sign on
each nonempty region. The remaining vertices form the connected
central band $B_t$, which contains the middle rows and columns and
intersects every row and column in an interval.
Figure~\ref{fig:corner-structure} illustrates the corresponding
geometry for the global estimation error on a square where the
signal is constant. The proof is given in
Appendix~\ref{subsec:geo-theorem-proof}. 

\begin{remark}[Boundary contact for the local estimator]
Under the hypotheses of Theorem~\ref{thm:geometry},
for every exterior edge
$\{v,w\}\in E_{\rm fr}$ with $v\in V_m$ and $w\in F_m$,
\begin{equation}\label{eq:boundary-contact}
 \widehat x_v^{\lambda,\bb}>t\ \Longrightarrow\ b_w>t,
 \qquad
 \widehat x_v^{\lambda,\bb}<-t\ \Longrightarrow\ b_w<-t.
\end{equation}
\end{remark}
Step~4 of the proof in Appendix~\ref{subsec:geo-theorem-proof}
establishes the remark.

In Theorem~\ref{thm:geometry}, boundary contact has a concrete interpretation: a positive threshold set attached to a corner can extend along either adjacent side only through the consecutive exterior values above $t$, starting at that corner. In particular, both exterior values adjacent to its corner vertex must exceed $t$. The negative assertion uses values below $-t$.

The following result gives boundary attachment at every nonnegative
threshold. It requires only $\lambda\ge\cutnorm{\by}$ and applies
even when the threshold condition for corner structure is not
satisfied.

\begin{proposition}[Boundary attachment for the local estimator]
\label{prop:boundary-attachment}
For the local estimator~\eqref{eq:model}, suppose
$\lambda\ge\cutnorm{\by}$. Then, simultaneously
for all exterior values $\bb$ and all $t\ge0$, every nonempty
connected component $C$ of $\{\widehat{\bx}^{\lambda,\bb}>t\}$ contains a
vertex $v$ with an exterior neighbor $w\in F_m$ satisfying $b_w>t$.
Similarly, every nonempty connected component of
$\{\widehat{\bx}^{\lambda,\bb}<-t\}$ contains a vertex with an exterior
neighbor whose value is less than $-t$.
\end{proposition}

Neither the positive nor the negative threshold set can have an
isolated interior component. The proof in
Appendix~\ref{subsec:boundary-attachment} is based on summing the
optimality equations over a connected component.

\paragraph{Scope of the square formulation.}
Theorems~\ref{thm:elementwise} and~\ref{thm:geometry} are formulated
on squares. Square coverings give the global error bound for
rectangular partitions in Section~\ref{subsec:main-piecewise}, and
Corollary~\ref{cor:gaussian-signal-elementwise} gives the corresponding
elementwise bounds. The four-corner geometry result in Theorem~\ref{thm:geometry} is applied to square
subregions of the constant rectangles.

\FloatBarrier
\Needspace{12\baselineskip}
\subsection{Global squared-error bounds for piecewise-constant signals}\label{subsec:main-piecewise}
Let \(\by=\btheta+\beps\) on \([n]^2\), where \(\btheta\) is constant on each cell of a rectangular partition and \(\beps=\by-\btheta\) denotes the fixed error array. We consider the global estimator~\eqref{eq:whole-image-model} and apply the directional fluctuation parameter~\eqref{eq:diamond-norm} to \(\beps\) on the full grid. The following theorem bounds the global squared error
$\norm{\widehat{\btheta}^\lambda-\btheta}_2^2$.

\begin{theorem}[Global error bounds for piecewise-constant signals]
\label{thm:piecewise-constant}
Consider the global estimator~\eqref{eq:whole-image-model}.
Suppose $[{n}]^2$ is partitioned into $k$ axis-aligned rectangles
$R_1,\ldots,R_k$, with side lengths $1\le a_i\le b_i$, and
$\btheta$ is constant on each rectangle.
For every $\lambda\ge\cutnorm{\beps}$,
\begin{equation}\label{eq:rectangular-det}
 \frac{1}{N}\norm{\widehat{\btheta}^\lambda-\btheta}_2^2
 \le\frac{C}{N}\sum_{i=1}^k
 \left\{(\lambda^2+{L_n}^4\cutnorm{\beps}^2)\frac{b_i}{a_i}
       +\frac{{L_n}^2\cutnorm{\beps}^4b_i}{\lambda^2}\right\}.
\end{equation}
In particular, if $b_i/a_i\le\rho$ for every $i$, where $\rho\ge1$, then
\begin{equation}\label{eq:rectangular-det-aspect}
 \frac{1}{N}\norm{\widehat{\btheta}^\lambda-\btheta}_2^2
 \le C\rho\left\{
 \frac{k\lambda^2}{N}
 +\frac{k{L_n}^4\cutnorm{\beps}^2}{N}
 +\frac{{L_n}^2\cutnorm{\beps}^4\sqrt{k}}{\lambda^2\sqrt N}\right\}.
\end{equation}
Here $C>0$ is a universal constant.
\end{theorem}
The bounds in Theorem~\ref{thm:piecewise-constant} do not require knowledge of the partition or of the
signal values on its rectangles. The proof covers each rectangle
with squares and sums the local bounds from
Corollary~\ref{cor:det-global}; see
Appendix~\ref{subsec:piecewise-constant}.

\section{Probabilistic guarantees and sharpness under Gaussian noise}
\label{sec:gaussian-consequences}\label{subsec:gaussian-consequences}
This section gives probabilistic guarantees and worst-case sharpness
results under Gaussian noise. Theorem~\ref{thm:gaussian-D} turns
the deterministic bounds of Section~\ref{sec:main-results} into
high-probability guarantees and risk upper bounds. The sharpness/lower bounds
are based on Gaussian averages and deterministic comparisons. In particular, 
Section~\ref{subsec:gaussian-signal-deviations} gives elementwise
and level-set geometric guarantees for the global
estimator~\eqref{eq:whole-image-model}.
Section~\ref{subsec:gaussian-risk-results} establishes elementwise
and averaged risk upper bounds, together with matching lower bounds,
for the local estimator~\eqref{eq:model}.
Section~\ref{subsec:gaussian-rectangle-results} gives global risk
upper bounds on rectangular partitions, matching lower bounds on
regular square partitions, and a discussion of the rectangular
assumption.

\subsection{High-probability elementwise error and level-set geometry}
\label{subsec:gaussian-signal-deviations}
Consider the observation model~\eqref{eq:intro-observation-model}
with $n\ge4$, independent $\varepsilon_p\sim N(0,\sigma^2)$,
and $\sigma>0$. Assume that the deterministic signal $\btheta$
is piecewise constant on an axis-aligned rectangular partition
of $[n]^2$, and use the free-boundary estimator
\eqref{eq:whole-image-model}. Throughout this subsection,
$L_n=\log(en)$ and $\kappa_{n,\delta}$ is defined
in~\eqref{eq:gaussian-D}.

\begin{corollary}[Simultaneous elementwise error for the global estimator]
\label{cor:gaussian-signal-elementwise}
Consider the global estimator~\eqref{eq:whole-image-model}.
For each vertex $p$, let $m_p$ be the shorter side length of the
partition rectangle containing $p$, and let $r_p,s_p$ be its
distances to the nearer horizontal and vertical sides of that
rectangle, counting boundary vertices at distance one.
There is a universal constant $C>0$ such that, for every
$0<\delta<1$, with probability at least $1-\delta$,
simultaneously for all vertices $p$ and all
$\lambda\ge\sigma\kappa_{n,\delta}$,
\begin{equation}\label{eq:gaussian-signal-elementwise}
 \bigl|\widehat\theta_p^\lambda-\theta_p\bigr|
 \le C\left\{
 \frac{\lambda}{m_p}
 +\frac{\sigma L_n\kappa_{n,\delta}}{\sqrt{r_ps_p}}
 +\frac{\sigma^2L_n^2\kappa_{n,\delta}^2}
 {\lambda \min(r_p,s_p)}
 \right\}.
\end{equation}
\end{corollary}

\begin{corollary}[Simultaneous large-error geometry for the global estimator]
\label{cor:gaussian-signal-geometry}
For the global estimator~\eqref{eq:whole-image-model}, there are universal constants $C_0,C_1>0$ such that, for every
$0<\delta<1$, with probability at least $1-\delta$, the following
holds simultaneously for every axis-aligned square $Q$ of side
$m\ge4$ contained in a partition rectangle and all $\lambda,t$
satisfying
\begin{equation}\label{eq:gaussian-signal-geometry}
 \lambda\ge C_0\sigma\kappa_{n,\delta},\qquad
 t\ge C_1\left\{
 \frac{\lambda}{m}
 +\frac{\sigma^2L_n^2\kappa_{n,\delta}^2}{\lambda}
 \right\}.
\end{equation}
Write $e_p^\lambda=\widehat\theta_p^\lambda-\theta_p$. The set
\[
 E_t(Q)=\{p\in Q:|e_p^\lambda|>t\}
\]
is a union of at most four sets, at most one in each quadrant of
$Q$, each corner-monotone from the corner of its quadrant, and
$e_p^\lambda$ has constant sign on each nonempty set. The central band
$Q\setminus E_t(Q)$ is nonempty and connected, contains the middle
row(s) and column(s) of $Q$, and intersects every row and column of
$Q$ in an interval.
\end{corollary}
The proofs of Corollaries~\ref{cor:gaussian-signal-elementwise}
and~\ref{cor:gaussian-signal-geometry} are given in
Appendix~\ref{subsec:gaussian-signal-proofs}. The geometry in Corollary~\ref{cor:gaussian-signal-geometry} is
described on squares contained in regions where the signal is constant.
Both corollaries hold on the same event
$\{\cutnorm{\beps}\le\sigma\kappa_{n,\delta}\}$, so their
conclusions hold jointly with probability at least $1-\delta$.
Taking $\delta=n^{-10}$ gives
$\kappa_{n,n^{-10}}\le(C_D+\sqrt5)L_n^{3/4}$ and recovers
Theorems~\ref{thm:informal-piecewise}
and~\ref{thm:informal-piecewise-geometry} after omitting the
logarithmic factors.

\subsection{Risk bounds and sharpness for the local estimator}
\label{subsec:gaussian-risk-results}\label{subsec:gaussian-lower-bounds}
\label{subsec:gaussian-tuning-results}
\paragraph{Upper bounds.}
Consider the local TV problem~\eqref{eq:model} with independent
$y_v\sim N(0,\sigma^2)$ and $\sigma>0$. Theorems~\ref{thm:gaussian-D}
and~\ref{thm:elementwise} give~\eqref{eq:det-main}, with probability
at least $1-\delta$, uniformly over all vertices, exterior values,
and $\lambda\ge\sigma\kappa_{m,\delta}$. The expected-risk bounds
below use a fixed deterministic regularization parameter of at least
order $\sigma L_m^{3/4}$ and also bound the error on the complement of
this event.

\begin{corollary}[Local elementwise and averaged risk]\label{cor:elementwise-risk}
For the local estimator~\eqref{eq:model}, suppose that $y_v$, $v\in V_m$, are independent $N(0,\sigma^2)$
variables, where $\sigma>0$. There are universal constants
$C_0,C>0$ such that, for every $m\ge4$, every fixed deterministic
\[
 \lambda\ge C_0\sigma L_m^{3/4},
\]
and every $p\in V_m$,
\begin{equation}\label{eq:main-risk}
 \E\!\left[
 \sup_{\bb\in\R^{F_m}}|\widehat x_p^{\lambda,\bb}|^2
 \right]
 \le C\left\{
 \frac{\lambda^2}{m^2}
 +\frac{\sigma^2L_m^{7/2}}{r_p s_p}
 +\frac{\sigma^4L_m^5}{\lambda^2\min(r_p,s_p)^2}
 \right\}.
\end{equation}
Consequently,
\begin{equation}\label{eq:mse-risk}
 \frac1{m^2}\E\!\left[
 \sup_{\bb\in\R^{F_m}}\|\widehat{\bx}^{\lambda,\bb}\|_2^2
 \right]
 \le C\left\{
 \frac{\lambda^2}{m^2}
 +\frac{\sigma^2L_m^{11/2}}{m^2}
 +\frac{\sigma^4L_m^5}{m\lambda^2}
 \right\}.
\end{equation}
\end{corollary}
The proof is given in Appendix~\ref{subsec:gaussian-elementwise-risk}.
We note that the supremum inside the
expectation makes these upper bounds valid even when the exterior
values depend on the observations. 
\paragraph{Elementwise lower bounds and sharpness.}
For the stated range of $\lambda$ in Corollary~\ref{cor:elementwise-risk},
Theorem~\ref{thm:elementwise-risk-lower} matches the elementwise bound
at every vertex, up to logarithmic factors, for the positive extremal
solution defined below, and hence in the worst case over deterministic
exterior values.

We establish a three-term lower bound for the local elementwise
risk when $\lambda\ge C_0\sigma L_m^{3/4}$, matching the upper
bound up to logarithmic factors. We then determine the averaged
risk for $\lambda=\sigma m^\alpha$.
For each fixed $\lambda>0$ and observation array $\by$, define the positive extremal solution by 
\[
 \bx_+(\by):=\lim_{M\to\infty}
 \widehat{\bx}^{\lambda,M\one_{F_m}}(\by).
\]
 Lemma~\ref{lem:comparison} shows that the limit
is achieved at a finite $M$ that may depend on $\by$
and $\lambda$. We write
$(u)_+=\max\{u,0\}$ for the positive part.

\begin{theorem}[Gaussian elementwise lower bounds and sharpness]
\label{thm:elementwise-risk-lower}
Consider the local estimator~\eqref{eq:model} with independent
$y_v\sim N(0,\sigma^2)$, where $m\ge4$ and $\sigma>0$.
Write $L_m=\log(em)$.
There are universal constants $c,C_0>0$ such that, for every
fixed deterministic $\lambda\ge C_0\sigma L_m^{3/4}$ and every
vertex $p\in V_m$,
\begin{equation}\label{eq:local-elementwise-risk-lower}
 \E(x_{+,p})_+^2
 \ge c\left\{
 \frac{\lambda^2}{m^2}+\frac{\sigma^2}{r_p s_p}
 +\frac{\sigma^4}{\lambda^2\min(r_p,s_p)^2\log^2\!\bigl(e\min(r_p,s_p)\bigr)}
 \right\}.
\end{equation}
Consequently, the same lower bound holds for
$\sup_{\bb\in\R^{F_m}}\E|\widehat x_p^{\lambda,\bb}|^2$,
where the supremum is over deterministic exterior values.
\end{theorem}
The proof is given in Appendix~\ref{subsec:elementwise-lower-proof}. 
When $\lambda\ge C_0\sigma L_m^{3/4}$, the same positive extremal
solution gives all three lower bounds at every vertex. Their sum
matches the upper bound in~\eqref{eq:main-risk} up to logarithmic
factors, also in the worst case over deterministic exterior values.

To compare the bounds, set
\[
 A_p(\lambda)
 :=\frac{\lambda^2}{m^2}
   +\frac{\sigma^2}{r_p s_p}
   +\frac{\sigma^4}{\lambda^2\min(r_p,s_p)^2}.
\]
Taking $C_0$ large enough for both results,
Theorem~\ref{thm:elementwise-risk-lower} and
Corollary~\ref{cor:elementwise-risk} imply that, for
$\lambda\ge C_0\sigma L_m^{3/4}$, 
\begin{equation}\label{eq:elementwise-sharpness}
 \begin{aligned}
 \frac{cA_p(\lambda)}{\log^2\!\bigl(e\min(r_p,s_p)\bigr)}
 &\le\E(x_{+,p})_+^2
 \le\sup_{\bb\in\R^{F_m}}\E
      |\widehat x_p^{\lambda,\bb}|^2\\
 &\le\E\!\left[\sup_{\bb\in\R^{F_m}}
      |\widehat x_p^{\lambda,\bb}|^2\right]
 \le CL_m^5 A_p(\lambda).
 \end{aligned}
\end{equation}
Thus, for this range of $\lambda$, the local elementwise risk bound
is sharp up to logarithmic factors at every vertex for the positive
extremal solution and, consequently, in the worst case over 
exterior values. 

\paragraph{Local averaged-risk sharpness.}
For the local TV problem~\eqref{eq:model} with independent
$y_v\sim N(0,\sigma^2)$ and $\sigma>0$, define the averaged risk of
the positive extremal solution and the worst risk over deterministic
exterior values by
\[
 R_{m,\sigma}^+(\lambda):=\frac1{m^2}\E\|\bx_+\|_2^2,
 \qquad
 R_{m,\sigma}^*(\lambda)
 :=\sup_{\bb\in\R^{F_m}}\frac1{m^2}
 \E\|\widehat{\bx}^{\lambda,\bb}\|_2^2.
\]
Lemma~\ref{lem:local-risk-lower} gives lower bounds for the
averaged risk of the positive extremal solution for every
$\lambda>0$. Combining these bounds with
Corollary~\ref{cor:elementwise-risk} gives the following
dependence on $\lambda=\sigma m^\alpha$.

\begin{corollary}[Local averaged risk]
\label{cor:power-match}
For the local TV problem~\eqref{eq:model} with independent
$y_v\sim N(0,\sigma^2)$ and $\sigma>0$, fix $0<\alpha<1$ and set
\[
 \lambda=\sigma m^\alpha,
 \qquad \gamma(\alpha)=\min\{2-2\alpha,\,1+2\alpha\}.
\]
There are universal constants $c,C>0$ such that, for all
$m\ge m_0(\alpha)$ with $m_0(\alpha)$ sufficiently large,
\[
 \begin{aligned}
 c\sigma^2m^{-\gamma(\alpha)}
 &\le R_{m,\sigma}^+(\lambda)
 \le R_{m,\sigma}^*(\lambda)\\
 &\le \frac1{m^2}\E\!\left[
 \sup_{\bb\in\R^{F_m}}\|\widehat{\bx}^{\lambda,\bb}\|_2^2
 \right]
 \le C\sigma^2L_m^5 m^{-\gamma(\alpha)}.
 \end{aligned}
\]
\end{corollary}
The proof is given in Appendix~\ref{subsec:power-match-proof}.
The exponent $\gamma(\alpha)$ is uniquely maximized at
$\alpha=1/4$, where $\gamma(1/4)=3/2$. Among choices of the form
$\lambda=\sigma m^\alpha$ with $0<\alpha<1$,
$\lambda=\sigma m^{1/4}$ gives the fastest polynomial decay:
both averaged risks are of order $\sigma^2m^{-3/2}$,
up to logarithmic factors.

\Needspace{18\baselineskip}
\subsection{Global risk bounds and sharpness on rectangular partitions}
\label{subsec:gaussian-rectangle-results}\label{subsec:corner-risk-sharpness}
\paragraph{Upper bound.}
For the global estimator~\eqref{eq:whole-image-model}, assume
that the observation errors $\varepsilon_p$ are independent
$N(0,\sigma^2)$ variables, where $\sigma>0$, and that $\btheta$
is deterministic.
\begin{corollary}[Global Gaussian risk on rectangular partitions]
\label{cor:rectangle-gaussian}
\label{cor:piecewise-constant}
\label{cor:rectangular-risk-tuned}
Consider the global estimator~\eqref{eq:whole-image-model}.
Suppose $[{n}]^2$ is partitioned into $k$ axis-aligned rectangles
$R_1,\ldots,R_k$, and $\btheta$ is constant on each $R_i$.
Let their side lengths be $1\le a_i\le b_i$.
There are universal constants $C_0,C>0$ such that the following hold.

\begin{enumerate}[label=(\roman*),leftmargin=*]
\item \textit{General aspect ratios.} For every fixed deterministic
$\lambda\ge C_0\sigma{L_n}^{3/4}$,
\begin{equation}\label{eq:rectangular-risk}
 \frac1N\E\norm{\widehat{\btheta}^\lambda-\btheta}_2^2
 \le\frac{C}{N}\sum_{i=1}^k
 \left\{(\lambda^2+\sigma^2{L_n}^{11/2})\frac{b_i}{a_i}
       +\frac{\sigma^4{L_n}^5b_i}{\lambda^2}\right\}.
\end{equation}

\item \textit{Bounded aspect ratios.} If, in addition, $b_i/a_i\le\rho$
for every $i$, where $\rho\ge1$, then, for a fixed sufficiently large
universal constant $C_{\rm reg}$, the choice
\[
 \lambda=C_{\rm reg}\sigma{L_n}^{5/4}\left(\frac Nk\right)^{1/8}
\]
gives
\begin{equation}\label{eq:rectangular-risk-tuned}
 \frac1N\E\norm{\widehat{\btheta}^\lambda-\btheta}_2^2
 \le C\rho\sigma^2\left\{
       {L_n}^{5/2}\left(\frac{k}{N}\right)^{3/4}
       +{L_n}^{11/2}\frac{k}{N}\right\}.
\end{equation}
\end{enumerate}
\end{corollary}
The proof is given in Appendix~\ref{subsec:gaussian-rectangles}.
Since $k\le N$ and ${L_n}\ge1$, the bound is at most
$2C\rho\sigma^2{L_n}^{11/2}(k/N)^{3/4}$; and if $k=O(1)$, then the second
term in~\eqref{eq:rectangular-risk-tuned} is eventually dominated by the
first as $N\to\infty$.

\paragraph{Global MSE sharpness.}
Corollary~\ref{thm:global-risk-sharpness} shows that the rate
$(k/N)^{3/4}$ in~\eqref{eq:rectangular-risk-tuned} is sharp up to
logarithmic factors for the global estimator with deterministic
regularization, in the worst case over signals on regular
partitions into equal squares.
Section~\ref{subsec:two-dimensional-risk-comparison} compares
this global risk bound with existing guarantees for two-dimensional
TV denoising.

To formulate the lower bound, consider a regular square partition as follows. Let $q\ge3$ and $m\ge4$ be integers, and set
\[
n=qm,\qquad k=q^2,\qquad N=n^2=km^2=q^2m^2.
\]
Partition $[n]^2$ into a $q\times q$ array of squares of side length $m$, indexed by
$Q_{ab},a,b\in[q],$ and let $\Theta_{q,m}$ denote the class of signals that are constant on each $Q_{ab}$.

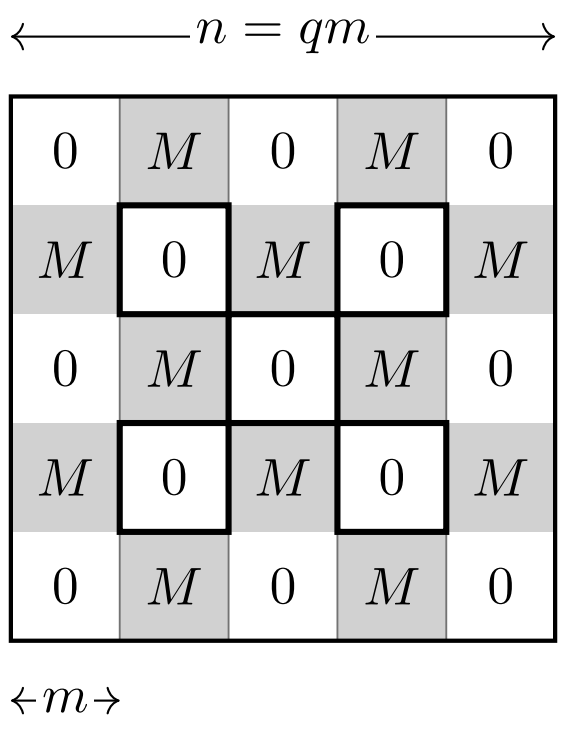
\begin{figure}[H]
\centering
\begin{tikzpicture}[x=0.72cm,y=-0.72cm,font=\small]
  \foreach \a in {1,...,5}{
    \foreach \b in {1,...,5}{
      \pgfmathtruncatemacro{\xx}{\b-1}
      \pgfmathtruncatemacro{\yy}{\a-1}
      \pgfmathtruncatemacro{\pp}{mod(\a+\b,2)}
      \ifnum\pp=0
        \fill[white] (\xx,\yy) rectangle ++(1,1);
        \node at (\xx+0.5,\yy+0.5) {$0$};
      \else
        \fill[black!18] (\xx,\yy) rectangle ++(1,1);
        \node at (\xx+0.5,\yy+0.5) {$M$};
      \fi
    }
  }
  \draw[step=1,black!55,line width=0.35pt] (0,0) grid (5,5);
  \draw[thick] (0,0) rectangle (5,5);
  \foreach \x/\y in {1/1,3/1,2/2,1/3,3/3}
      \draw[line width=1.15pt] (\x,\y) rectangle ++(1,1);
  \draw[<->] (0,-0.55)--(5,-0.55)
      node[midway,fill=white,inner sep=1pt] {$n=qm$};
  \draw[<->] (0,5.55)--(1,5.55)
      node[midway,fill=white,inner sep=1pt] {$m$};
\end{tikzpicture}
\caption{Checkerboard example used for the global lower bound, shown
for $q=5$. Each partition square has side length $m$; white squares
have signal value $0$ and shaded squares have signal value $M$.
The outlined interior zero squares are surrounded on all four sides
by squares with signal value $M$.}
\label{fig:checkerboard-family}
\end{figure}

We use the following checkerboard example to establish the lower bound. For $M>0$, define $\btheta^{(M)}\in\Theta_{q,m}$ by
\[
\theta^{(M)}|_{Q_{ab}}
=
\begin{cases}
0, & a+b\ \text{even},\\
M, & a+b\ \text{odd}.
\end{cases}
\]
Neighboring squares therefore have signal values $0$ and $M$,
with jump size $M$. As $M\to\infty$, each interior zero square is
surrounded by squares whose signal value becomes arbitrarily large.
Figure~\ref{fig:checkerboard-family} illustrates the checkerboard
example.

\begin{corollary}[Sharpness of the global MSE rate]
\label{thm:global-risk-sharpness}
Consider the global estimator~\eqref{eq:whole-image-model}
under the observation model~\eqref{eq:intro-observation-model}
with independent $N(0,\sigma^2)$ noise, where $\sigma>0$.
There is a universal constant $c>0$ such that, for every fixed
deterministic $\lambda>0$,
\begin{equation}\label{eq:global-risk-lower}
 \liminf_{M\to\infty}
 \frac1N\E_{\btheta^{(M)}}
 \norm{\widehat{\btheta}^{\lambda}-\btheta^{(M)}}_2^2
 \ge c\left\{
 \frac{\lambda^2}{m^2}
 +\frac{\sigma^2}{m\max\{1,(\lambda/\sigma)^2\}}
 \right\}.
\end{equation}
The same lower bound therefore holds for the supremum of the risk
over $\Theta_{q,m}$. Consequently,
\begin{equation}\label{eq:global-risk-optimal}
 \inf_{\lambda>0}\sup_{\btheta\in\Theta_{q,m}}
 \frac1N\E_{\btheta}
 \norm{\widehat{\btheta}^{\lambda}-\btheta}_2^2
 \ge c\sigma^2\left(\frac{k}{N}\right)^{3/4}.
\end{equation}
\end{corollary}
The proof is given in Appendix~\ref{subsec:global-risk-lower-proof}.
Chatterjee and Goswami~\cite[Theorem~II.3]{CGTV} established the
$N^{-3/4}$ exponent for constrained TV in the low-noise limit;
Section~\ref{subsec:two-dimensional-risk-comparison} explains the
connection to penalized TV. The bound here gives the dependence on
$k$ for regular square partitions and the two contributions as
functions of $\lambda$.
Together with Corollary~\ref{cor:rectangle-gaussian}, this proves
sharpness of $(k/N)^{3/4}$ up to logarithmic factors within the
penalized estimator family with deterministic regularization.

Let $Q$ be any interior square of the fixed checkerboard example.
For each fixed $\lambda>0$, as $M$ increases, the estimator on a
zero square becomes the positive extremal solution. On an $M$-square,
the estimator centered by $M$ becomes the negative extremal solution.
Gaussian sign symmetry gives the same squared-risk bounds in both
cases. This transfers the averaged bounds in
Lemma~\ref{lem:local-risk-lower} for every $\lambda>0$ and the
three-term elementwise bound in~\eqref{eq:local-elementwise-risk-lower}
for $\lambda\ge C_0\sigma L_m^{3/4}$.

\paragraph{Role of the rectangular assumption.}
\label{par:rectangular-assumption}
The global risk bound uses the axis-aligned rectangular structure
of the constant regions. A related obstruction for other region
shapes is known for constrained TV: Chatterjee and
Goswami~\cite[Theorem~II.6]{CGTV} establish an $N^{-1/2}$
lower-bound rate in the low-noise limit for a signal with a diagonal
boundary. This illustrates the role of region geometry, while our
penalized-estimator guarantee is stated for rectangular partitions.

\section{Comparison with existing results}\label{sec:examples-comparisons}
Sections~\ref{subsec:one-dimensional-comparison}
and~\ref{subsec:two-dimensional-risk-comparison} compare these results
with one-dimensional elementwise bounds and existing two-dimensional
risk guarantees.

{
\subsection{Comparison with one-dimensional elementwise bounds}
\label{subsec:one-dimensional-comparison}

We compare the squared-loss result of Zhang and Chatterjee~\cite{ZC},
which assumes independent, identically distributed centered
sub-Gaussian errors of scale $\sigma$, with our Gaussian grid bounds.

\paragraph{Elementwise bounds.}
Let $m$ be the length of an interval on which the true signal is
constant and $d\ge1$ the distance of coordinate $i$ to its nearer
endpoint, counting the endpoint as
distance one. Write $\lambda_1$ for the regularization parameter in~\cite{ZC}.
Their quadratic loss is half the squared loss used here, so the penalty
parameters are related by $\lambda=2\lambda_1$.
Up to logarithmic factors, their high-probability bound can be written as
\begin{equation}\label{eq:comparison-one-dimensional}
 |\widehat\theta_i-\theta_i|
 \lesssim_{\log}
 \frac{\sigma}{\sqrt d}
 +\frac{\sigma^2}{\lambda_1}
 +\frac{\lambda_1}{m}.
\end{equation}
Here $\lesssim_{\log}$ denotes an inequality up to universal constants
and logarithmic factors in the sample size and confidence parameter. The additional
block-average term $\sigma/\sqrt m$ is absorbed by the distance term.
For an $m\times m$ square on which the true signal is constant,
after centering at that signal value,
Theorems~\ref{thm:gaussian-D} and~\ref{thm:elementwise} give
\begin{equation}\label{eq:comparison-two-dimensional}
 \sup_{\bb}|\widehat x_p^{\lambda,\bb}|
 \lesssim_{\log}
 \frac{\lambda}{m}
 +\frac{\sigma}{\sqrt{r_p s_p}}
 +\frac{\sigma^2}{\lambda \min(r_p,s_p)},
\end{equation}
on the Gaussian event, for $\lambda$ satisfying the stated lower bound.
Both bounds separate
boundary-induced bias, a location term, and a term decreasing with
regularization. The bias scale $\lambda/m$ reflects the ratio of
order $m$ exterior edges to $m^2$ vertices.

The one-dimensional bound holds for every fixed $\lambda_1>0$;
our deterministic bound holds simultaneously for
$\lambda\ge\cutnorm{\by}$, while the expected-risk corollaries use
fixed deterministic $\lambda\ge C_0\sigma L_m^{3/4}$.
Neither bound imposes an upper bound on the regularization parameter,
although the increasing bias term limits its usefulness for large $\lambda$.

\paragraph{Corners and sides.}
An interval endpoint is analogous to a square corner: in the local
TV problem, exterior edges can offset the TV penalty on all inward
edges. The factor $1/\sqrt{r_p s_p}$ distinguishes corners from
straight sides. When $r_p$ is bounded and $s_p\asymp m$, the location
term decreases as $m^{-1/2}$, whereas it need not decrease when both
distances are bounded. The side term
$\sigma^2/(\lambda \min(r_p,s_p))$ also decreases as $\lambda$
grows, so a suitable choice of $\lambda$ permits consistency near a
straight side even when the distance to the boundary remains bounded.
Distinguishing a corner from a straight side therefore requires
both the horizontal and vertical distances to the boundary.
The local lower bounds in Theorem~\ref{thm:elementwise-risk-lower}
hold for the positive extremal solution and imply lower bounds for the
worst risk over deterministic exterior values; favorable boundary
values can yield smaller errors.

\paragraph{Additional geometry in the proof.}
Both analyses localize the global estimator to a region on which the true
signal is constant and obtain bounds that hold uniformly for all
boundary values; Lemma~\ref{lem:patch-restriction}
gives this reduction for squares. In one dimension, each connected
threshold component is an interval, and lowering its estimated values
changes the TV penalty only at its two ends. Contact with the endpoints of the containing interval on which the
true signal is constant yields the distance and bias terms.

On the grid, threshold sets can have holes or multiple occupied runs
in a row or column. Applying an interval argument separately to rows
leaves vertical penalties at every vertex. The boundary-crossing
inequality of Lemma~\ref{lem:leakage} and nested-rectangle propagation
in Lemma~\ref{lem:peeling} provide control in both directions without
assuming connectedness or monotonicity.

\paragraph{Fluctuation control.}
Zhang and Chatterjee~\cite{ZC} control interval sums using a
nonasymptotic law of the iterated logarithm. Our analysis requires
sums over arbitrary vertex sets, normalized by
$\sqrt{P_1(S)P_2(S)}$; Theorem~\ref{thm:gaussian-D} supplies this
control by combining the Ding--Wirth estimate with a lattice-stretching
argument that balances the two directional perimeters.
The expected-risk consequences are developed in
Section~\ref{sec:gaussian-consequences}.

}

\subsection{Comparison with two-dimensional risk bounds}
\label{subsec:two-dimensional-risk-comparison}

We compare the global risk bound in
Corollary~\ref{cor:rectangle-gaussian} with existing guarantees
for two-dimensional TV denoising under the Gaussian observation
model~\eqref{eq:intro-observation-model}, with $N=n^2$.

\paragraph{Oracle inequalities for penalized TV.}
H\"utter and Rigollet~\cite[Corollary~5]{HR} establish sharp oracle
inequalities for the penalized estimator.
With regularization parameter $\lambda\asymp\sigma\log N$ under
the loss normalization in~\eqref{eq:whole-image-model}, their bound
implies, with high probability,
\[
 \frac1N\norm{\widehat{\btheta}^{\lambda}-\btheta}_2^2
 \;\widetilde{\lesssim}\;
 \frac{\sigma^2\{J(\btheta)+1\}}{N},
\]
Here $J(\btheta)$ counts grid edges with different signal values at
their endpoints, and logarithmic factors are suppressed.
The bound allows arbitrary region shapes and jump sizes. It gives
order $\sigma^2N^{-1/2}$ when $J(\btheta)\asymp\sqrt N$ and
order $\sigma^2/N$ for a constant signal, up to logarithms.
For a fixed number of rectangles with bounded aspect ratios,
Corollary~\ref{cor:rectangle-gaussian} instead uses
$\lambda\asymp\sigma N^{1/8}$ and gives order
$\sigma^2N^{-3/4}$, up to logarithms in both expressions.
These guarantees use different regularization parameters; the
oracle inequality is a high-probability statement, whereas our
corollary bounds expected squared error.

\paragraph{Adaptation to rectangular partitions.}
Chatterjee and Goswami~\cite[Theorem~II.2]{CGTV} establish an
adaptive bound for the constrained TV estimator. Setting the
constraint radius to the true total variation gives an oracle
choice that depends on the unknown signal.
For signals constant on $k$ rectangles with aspect ratios at most
$\rho$, their expected-risk upper bound is
\[
 C_{\rm CG}(\rho)\sigma^2(\log N)^9 k^{5/4}N^{-3/4}.
\]
For the penalized estimator~\eqref{eq:whole-image-model},
Corollary~\ref{cor:rectangle-gaussian} gives, with the deterministic
choice of $\lambda$ stated there,
\[
 \frac1N\E\norm{\widehat{\btheta}^{\lambda}-\btheta}_2^2
 \le C\rho\sigma^2\left\{
 (\log N)^{5/2}(k/N)^{3/4}
 +(\log N)^{11/2}k/N\right\}.
\]
For fixed $k$ and $\rho$, both bounds have the power $N^{-3/4}$.
The regularization parameter here depends only on $N$, $k$, and
$\sigma$; it does not require the signal's total variation,
rectangle locations, or signal values. This answers the question
posed by Chatterjee and Goswami~\cite[Section~III-F]{CGTV} about
whether penalized TV attains the $N^{-3/4}$ rate, and supplies an
explicit deterministic choice of the regularization parameter.

Van de Geer's 2019 Oberwolfach Lecture~\cite[Section~5]{vdGTV}
announces the $N^{-3/4}$ rate, up to logarithmic factors, for
penalized grid TV in an equal two-rectangle example. The report
outlines a path construction and an effective-sparsity argument,
and her subsequent lecture slides~\cite[pp.~70--75]{vdGTVslides}
give further details for that example.
Corollary~\ref{cor:rectangle-gaussian} treats general rectangular
partitions with bounded aspect ratios and gives an expected-risk
upper bound of order $\sigma^2(k/N)^{3/4}$, up to logarithmic
factors.

\paragraph{Lower bounds and dependence on regularization.}
For even $n$, Chatterjee and Goswami~\cite[Theorem~II.3]{CGTV}
consider the fixed signal $\theta^*_{ij}=\one\{j>n/2\}$ on two
equal rectangles and prove
\[
 \lim_{\sigma\downarrow0}\sigma^{-2}
 \operatorname{MSE}(\widehat{\btheta}_{V^*},\btheta^*)
 \ge cN^{-3/4},
\]
where $V^*$ is the true normalized total variation and MSE is
expected squared error divided by $N$.
Their proof lower-bounds the Gaussian width of the tangent cone
at this signal. The general proximal-risk comparison~\cite{OHprox},
also discussed in~\cite[Section~III-F]{CGTV}, transfers the
low-noise obstruction to penalized TV when the penalty-to-noise
ratio is fixed.

For $a>0$ and a standard Gaussian array $\bz$, positive homogeneity
of TV gives
\[
 \widehat{\btheta}^{\lambda}(a\btheta^*+\sigma\bz)
 =a\widehat{\btheta}^{\lambda/a}
       (\btheta^*+(\sigma/a)\bz).
\]
Letting $a\to\infty$ therefore yields a
$c\sigma^2N^{-3/4}$ lower bound on the worst-case penalized risk
over all amplitudes of this two-rectangle signal, for every
$\sigma>0$ and deterministic $\lambda>0$.
The checkerboard example in
Corollary~\ref{thm:global-risk-sharpness}, with increasing jump
sizes, supplies the dependence on $k$ for regular square partitions
and the explicit bias and side contributions
in~\eqref{eq:global-risk-lower}.
Theorem~\ref{thm:elementwise-risk-lower} matches the combined
three-term elementwise risk bound, up to logarithmic factors, at
every vertex in its stated parameter range, using the same positive
extremal solution for all three lower bounds. These local conclusions therefore also hold for the
worst risk over deterministic exterior values.
The amplitude scaling uses unrestricted jump sizes. Minimax results
over fixed bounded-variation balls~\cite{SWT,SWHT} concern a
different signal class.

\section{Discussion}\label{sec:discussion}

For signals constant on rectangular partitions,
Theorem~\ref{thm:elementwise} gives a local elementwise bound that
separates regularization bias from side and corner effects, and
Corollary~\ref{cor:gaussian-signal-elementwise} transfers it to the
global estimator under Gaussian noise. Theorem~\ref{thm:geometry}
shows that the sets where the local estimator magnitude exceeds a
sufficiently large threshold have a corner-monotone structure, and
Corollary~\ref{cor:gaussian-signal-geometry} gives the corresponding
geometry of global estimation errors without requiring recovery of
the partition boundaries. Theorem~\ref{thm:elementwise-risk-lower}
establishes sharpness of the combined three-term risk bound at
every vertex in its stated range of $\lambda$, up to logarithmic
factors. The same positive extremal solution gives all three lower
bounds. Through the checkerboard example in
Corollary~\ref{thm:global-risk-sharpness}, they also apply to the
worst-case global risk at vertices of interior partition squares,
with unrestricted jump sizes.

Corollary~\ref{cor:rectangle-gaussian} establishes the global
expected-risk rate $\sigma^2(k/N)^{3/4}$, up to logarithmic
factors, for the standard penalized estimator on rectangular
partitions of bounded aspect ratio. For fixed $k$, this gives the
$N^{-3/4}$ rate with an explicit deterministic regularization
parameter depending on $N$, $k$, and $\sigma$, but not on the
signal's total variation, rectangle locations, or signal values.
Together with Corollary~\ref{thm:global-risk-sharpness}, this
establishes the optimal polynomial dependence on $k/N$ within
the TV-penalized estimator family under deterministic choice of regularization parameter $\lambda$.

Several questions remain. The additional logarithmic factors introduced by the propagation and summation arguments may be improvable. It would also be interesting to extend the spatial error bounds and threshold-set geometry beyond axis-aligned rectangular regions. More broadly, the local analysis developed here may be adaptable to other graph structures, noise distributions, and loss functions, including robust or quantile-regression losses. Understanding which aspects of the side--corner decomposition and large-error geometry persist in these settings is an interesting direction for future work. 
\clearpage
\appendix

\section{Shared deterministic tools}\label{sec:preliminaries}

This appendix gives the exterior comparison and threshold-set
optimality lemmas used in the proofs, together with a mean identity
for the positive extremal solution.
For fixed exterior values $\bb$, let $\mathcal J_{\bb}$ denote the
objective in~\eqref{eq:model}.

\subsection{Auxiliary lemmas}

\subsubsection{Extremal exterior comparison}

For fixed $\by$ and $\lambda\ge0$, define the extremal solutions by
\begin{equation}\label{eq:envelopes}
 \bx_+
 :=
 \lim_{M\to\infty}
 \widehat{\bx}^{\lambda,M\one_{F_m}},
 \qquad
 \bx_-
 :=
 \lim_{M\to\infty}
 \widehat{\bx}^{\lambda,-M\one_{F_m}}.
\end{equation}
The lemma below proves that these limits exist and stabilize for
finite $M$. Sign reversal in the objective gives
$\widehat{\bx}^{\lambda,-\bb}(-\by)
=-\widehat{\bx}^{\lambda,\bb}(\by)$, and hence
\begin{equation}\label{eq:extremal-sign}
 \bx_-(\by)=-\bx_+(-\by).
\end{equation}

\begin{lemma}[Extremal exterior comparison]\label{lem:comparison}
The limits in~\eqref{eq:envelopes} exist and are attained by
sufficiently large finite constant exterior values.
In addition, for every finite exterior vector $\bb\in \R^{F_m}$,
\begin{equation}\label{eq:sandwich}
 \bx_-
 \le
 \widehat{\bx}^{\lambda,\bb}
 \le
 \bx_+
 \qquad\text{coordinatewise}.
\end{equation}
For every vertex $p\in V_m$, the map
$\by\mapsto \sup_{\bb}|\widehat x_p^{\lambda,\bb}(\by)|$ 
is measurable and satisfies
\begin{equation}\label{eq:crude}
 \sup_{\bb}
 \bigl|\widehat x_p^{\lambda,\bb}\bigr|
 \le |y_p|+2\lambda.
\end{equation}
These conclusions also hold on a rectangle with any fixed set of
specified exterior edges; exterior penalties are omitted on free
sides. For a specified sign pattern
$\boldsymbol\eta=(\eta_w)$, with $\eta_w\in\{-1,1\}$, the solution
with exterior values $b_w=M\eta_w$ stabilizes whenever
$M>\norm{\by}_\infty+2\lambda$. Denote this stabilized map by
$T_{\boldsymbol\eta}(\by)$. It is measurable and satisfies
\[
 T_{\boldsymbol\eta}(\by+c\one)
 =T_{\boldsymbol\eta}(\by)+c\one,
 \qquad c\in\R.
\]
For the square $V_m$ with positive exterior values on all four
sides, the stabilized solution also satisfies
\begin{equation}\label{eq:test-set-balance}
 \sum_{v\in V_m}x_{+,v}
 =
 \sum_{v\in V_m}y_v+2\lambda m.
\end{equation}
Thus $\overline x_+=\overline y+2\lambda/m$, where bars denote
averages over $V_m$.
\end{lemma}

\subsubsection{A perturbation inequality}

The connection between threshold sets of discrete TV minimizers
and parametric minimum cuts is classical; see Chambolle and
Darbon~\cite[Proposition~3.3 and Section~4.1]{CD}.
The direct perturbation argument below gives the inequality needed
here, including fixed exterior values.

Write $\widehat{\bx}=\widehat{\bx}^{\lambda,\bb}$.
Extend the estimated values to the frame by setting
$\widetilde x_v=\widehat x_v$ for $v\in V_m$ and
$\widetilde x_w=b_w$ for $w\in F_m$.
The inequality below is obtained by decreasing all estimated values on a set
\(C\subseteq\{\widehat x>t\}\) by a common small amount
\(\varepsilon>0\). Since \(\widehat{\bx}\) minimizes the original TV
objective, this perturbation cannot decrease the objective. For an edge
leaving \(C\), its TV contribution decreases by \(\varepsilon\) if the
other endpoint has value at most \(t\), and can increase by at most
\(\varepsilon\) if the other endpoint has value above \(t\). Counting
these two types of edges gives the inequality in the following lemma.

\begin{lemma}[Optimality on a threshold subset]
\label{lem:geo-binary}
Fix $t\in\R$ and let
$C\subseteq\{v\in V_m:\widehat x_v>t\}$.
Define
\[
 \begin{aligned}
 N_-(C,t)
 &:=
 \#\bigl\{\{v,w\}\in E_{\rm int}\cup E_{\rm fr}:
       v\in C,\ w\notin C,\ \widetilde x_w\le t\bigr\},\\
 N_+(C,t)
 &:=
 \#\bigl\{\{v,w\}\in E_{\rm int}\cup E_{\rm fr}:
       v\in C,\ w\notin C,\ \widetilde x_w>t\bigr\}.
 \end{aligned}
\]
Each edge is counted once, with $v$ denoting its endpoint in $C$.
In particular,
$N_-(C,t)+N_+(C,t)=\operatorname{Per}(C)$.
Then
\begin{equation}\label{eq:threshold-cut-optimality}
 t|C|-\sum_{v\in C}y_v
 +\frac{\lambda}{2}\{N_-(C,t)-N_+(C,t)\}
 \le0.
\end{equation}
\end{lemma}

\subsection{Proofs of the auxiliary lemmas}
\label{subsec:shared-tools-proofs}
\begin{proof}[Proof of Lemma~\ref{lem:comparison}]
\emph{Uniform bound and stabilization.}
The coordinatewise optimality conditions for~\eqref{eq:model} give
\[
 0
 =
 2(\widehat x_p^{\lambda,\bb}-y_p)
 +\lambda\sum_{\{p,q\}\in E}s_{pq},
 \qquad |s_{pq}|\le1.
\]
Since any vertex has at most four neighbors in $E$, 
\[
 |\widehat x_p^{\lambda,\bb}-y_p|
 \le
 \frac{\lambda}{2}\sum_{\{p,q\}\in E}|s_{pq}|
 \le2\lambda.
\]
This proves~\eqref{eq:crude}.

Put $B=\norm{\by}_\infty+2\lambda$. Every estimated value lies
in $[-B,B]$. If $M>B$, each exterior-edge subgradient is $-1$
for the exterior vector $M\one_{F_m}$ and $+1$ for
$-M\one_{F_m}$. The optimality equations are therefore independent
of $M$ in each case. Uniqueness of the minimizer proves the stated
stabilization and the existence of the limits in~\eqref{eq:envelopes}.

\medskip
\noindent
\emph{Monotonicity and comparison.}
Suppose $\bb\le\bb'$, and let
\[
 \bu=\widehat{\bx}^{\lambda,\bb},
 \qquad
 \bv=\widehat{\bx}^{\lambda,\bb'}.
\]
Write $\bu\wedge\bv$ and $\bu\vee\bv$ for the coordinatewise
minimum and maximum. For real numbers $a,b,c,d$, the lattice
inequality
\[
 \begin{aligned}
 &|\min\{a,c\}-\min\{b,d\}|
 +|\max\{a,c\}-\max\{b,d\}|\\
 &\hspace{3em}\le |a-b|+|c-d|
 \end{aligned}
\]
applies to each internal edge. It also applies to each frame
edge, because $\bb\le\bb'$ implies
$\bb\wedge\bb'=\bb$ and $\bb\vee\bb'=\bb'$.
The sum of the two squared-loss terms is unchanged by taking
coordinatewise minima and maxima. Consequently,
\[
 \mathcal J_{\bb}(\bu\wedge\bv)
 +
 \mathcal J_{\bb'}(\bu\vee\bv)
 \le
 \mathcal J_{\bb}(\bu)
 +
 \mathcal J_{\bb'}(\bv).
\]
On the other hand, minimality gives
\[
 \mathcal J_{\bb}(\bu)
 \le
 \mathcal J_{\bb}(\bu\wedge\bv),
 \qquad
 \mathcal J_{\bb'}(\bv)
 \le
 \mathcal J_{\bb'}(\bu\vee\bv).
\]
Both inequalities must therefore be equalities. Since the
squared loss makes each objective strictly convex, its minimizer
is unique. Hence $\bu=\bu\wedge\bv$, proving monotonicity in
the exterior values. For any finite $\bb$, choose
$M>\max\{B,\norm{\bb}_\infty\}$. Comparing
$-M\one_{F_m}\le\bb\le M\one_{F_m}$ and using stabilization
gives~\eqref{eq:sandwich}.

\medskip
\noindent
\emph{Measurability.}
For fixed $\bb$, let $\bu$ and $\bv$ be the solutions associated with
observations $\by$ and $\by'$, respectively. Subtracting their
optimality conditions and using monotonicity of the subgradient
of the convex penalty gives
\[
 \norm{\bu-\bv}_2^2
 \le
 \ip{\by-\by'}{\bu-\bv}.
\]
By Cauchy--Schwarz,
\[
 \norm{\bu-\bv}_2\le\norm{\by-\by'}_2.
\]
Thus each estimator with finite exterior values is continuous in the observations.
Taking $M$ through the positive integers in~\eqref{eq:envelopes}
shows that $\bx_+$ and $\bx_-$ are elementwise limits of continuous
maps and are therefore measurable.
The sandwich bound and attainment of both extremal solutions give
\[
 \sup_{\bb}\bigl|\widehat x_p^{\lambda,\bb}\bigr|
 =\max\{x_{+,p},-x_{-,p}\},
\]
which proves the required measurability.
For fixed observations, the solution is also continuous in $\bb$:
the uniform coordinate bound gives bounded minimizers, and the
objectives converge uniformly on bounded sets when $\bb$ converges.
Every subsequential limit therefore minimizes the limiting objective,
and uniqueness identifies that limit. Consequently,
\[
 \sup_{\bb\in\R^{F_m}}\norm{\widehat{\bx}^{\lambda,\bb}(\by)}_2^2
 =\sup_{\bb\in\mathbb Q^{F_m}}\norm{\widehat{\bx}^{\lambda,\bb}(\by)}_2^2.
\]
The right-hand side is a countable supremum of continuous functions
of $\by$, so the norm supremum used in
Corollary~\ref{cor:elementwise-risk} is measurable as well.

\medskip
\noindent
\emph{Rectangles and specified exterior signs.}
The coordinate bound, comparison, and continuity arguments use
only the absolute-difference penalties and the fact that each
vertex has at most four incident edges. They therefore apply to
rectangles with any fixed set of specified exterior edges.
For a specified sign pattern, every exterior value with
$\eta_w b_w>B$ gives the same exterior-edge subgradient. The
unique solution thus depends only on the signs, once the magnitudes
exceed $B$. In particular, $T_{\boldsymbol\eta}$ is the elementwise
limit through positive integers $M$ of continuous finite-boundary
maps, so it is measurable.

For any $c\in\R$, shifting the observations and all specified
exterior values by $c$ shifts the finite-boundary minimizer by
$c\one$, since the objective is unchanged when the estimated values
are shifted as well. Choose $M$ sufficiently large for both
observation arrays $\by$ and $\by+c\one$ and for the exterior
values $M\boldsymbol\eta$ and $M\boldsymbol\eta+c\one$.
Independence from the sufficiently large exterior magnitudes then
gives the stated translation equivariance of
$T_{\boldsymbol\eta}$.

For the mean identity, prescribe the constant exterior value
$M>\norm{\by}_\infty+2\lambda$ on the square $V_m$.
For $t$ sufficiently close to zero, all entries of
$\bx_++t\one$ remain below $M$. Internal TV is unchanged by
this translation, and each of the $4m$ exterior-edge terms changes
by $-\lambda t$. Differentiating the objective at its minimum
$t=0$ gives
\[
 0=2\sum_{v\in V_m}(x_{+,v}-y_v)-4\lambda m,
\]
which proves~\eqref{eq:test-set-balance}.
\end{proof}

\begin{proof}[Proof of Lemma~\ref{lem:geo-binary}]
If $C=\varnothing$, the result is immediate. Assume $C\ne\varnothing$.
Choose
$0<\varepsilon<\min_{v\in C}(\widehat x_v-t)$
and compare $\widehat{\bx}$ with
$\widehat{\bx}-\varepsilon\one_C$.
Only edges leaving $C$ can change their total variation
contribution.
Each edge counted by $N_-(C,t)$ decreases its contribution by
$\varepsilon$, since its endpoint in $C$ remains above $t$
and its other endpoint has value at most $t$.
Each edge counted by $N_+(C,t)$ increases its contribution by
at most $\varepsilon$, by the triangle inequality.
Optimality and expansion of the squared loss therefore give
\[
 \begin{aligned}
 0
 &\le
 \mathcal J_{\bb}(\widehat{\bx}-\varepsilon\one_C)
 -\mathcal J_{\bb}(\widehat{\bx})\\
 &\le
 2\varepsilon\sum_{v\in C}(y_v-\widehat x_v)
 +\varepsilon^2|C|
 +\lambda\varepsilon\{N_+(C,t)-N_-(C,t)\}.
 \end{aligned}
\]
Dividing by $2\varepsilon$ and letting $\varepsilon\downarrow0$
yields
\[
 \sum_{v\in C}(\widehat x_v-y_v)
 +\frac{\lambda}{2}\{N_-(C,t)-N_+(C,t)\}
 \le0.
\]
Since $\widehat x_v>t$ for $v\in C$, 
\eqref{eq:threshold-cut-optimality} is proved.
\end{proof}

\section{Proof of Theorem~\ref{thm:gaussian-D}}\label{sec:gaussian-flows}
We prove Theorem~\ref{thm:gaussian-D} by reducing directional
fluctuations to the ordinary-perimeter Gaussian statistic of Ding and
Wirth. The reduction stretches the lattice in one coordinate direction
and represents each original Gaussian variable as a normalized block sum.

\subsection{Auxiliary lemmas}
\subsubsection{Geometric bounds and Lipschitz continuity}
\begin{lemma}[Area and directional perimeters]\label{lem:perimeter-preliminaries}
For every nonempty $S\subseteq V_m$,
\begin{equation}\label{eq:area-perimeter}
 |S|\le\frac{P_1(S)P_2(S)}4,
\end{equation}
and
\begin{equation}\label{eq:perimeter-ratio}
 \frac1m\le\frac{P_1(S)}{P_2(S)}\le m.
\end{equation}
\end{lemma}

The area bound makes $\cutnorm{\cdot}$ a $1/2$-Lipschitz function. We will use it and concentration of measure to pass from the bound on the expectation of \(\cutnorm{\by}\) to a high-probability bound.
\begin{lemma}[Lipschitz continuity of the fluctuation parameter]
\label{lem:cutnorm-lipschitz}
For any arrays \(\by,\by'\) on \(V_m\),
\begin{equation}\label{eq:norm-lipschitz}
\big|\cutnorm{\by}-\cutnorm{\by'}\big|
\le \frac12\|\by-\by'\|_2.
\end{equation}
\end{lemma}

\subsubsection{The ordinary-perimeter Gaussian bound}
For a finite lattice rectangle $R\subset\Z^2$ and an array $\bz$ on
$R$, define
\begin{equation}\label{eq:dw-statistic}
 \mathcal G_R(\bz)
 :=\max_{\varnothing\ne S\subseteq R}
       \frac{\left|\sum_{v\in S}z_v\right|}{\operatorname{Per}(S)},
 \qquad
 \mathcal G(\by):=\mathcal G_{V_m}(\by).
\end{equation}
Here $\operatorname{Per}(S)$ counts all lattice edges leaving $S$,
including edges leaving $R$.
The following lemma adapts the bounds in \cite{DW} to our setting.
\begin{lemma}[Ding--Wirth perimeter Gaussian bound]\label{lem:dw-perimeter}
There are universal constants $c,C>0$ such that, for independent
$y_v\sim N(0,\sigma^2)$ on $V_m$, with $m\ge4$ and $\sigma>0$,
\begin{equation}\label{eq:dw-square-mean}
 c\sigma L_m^{3/4}
 \le \E\mathcal G(\by)
 \le C\sigma L_m^{3/4}.
\end{equation}
\end{lemma}

\subsection{Proof of Theorem~\ref{thm:gaussian-D} by lattice stretching}
\begin{proof}[Proof of Theorem~\ref{thm:gaussian-D}]
\emph{Step 1: a finite family of auxiliary  statistics $K_a(\by)$.}
For $a>0$, introduce the auxiliary statistic
\[
 K_a(\by)
 :=\max_{\varnothing\ne S\subseteq V_m}
 \frac{\left|\sum_{v\in S}y_v\right|}
      {P_1(S)/\sqrt a+\sqrt a\,P_2(S)},
 \qquad
 \mathcal A_m:=\{2^k:k\in\Z,\ |k|\le\lceil\log_2m\rceil\}.
\]
For every nonempty $S$, set $r=P_1(S)/P_2(S)$ and choose
$a=2^{\lfloor\log_2r\rfloor}$. By~\eqref{eq:perimeter-ratio},
$a\in\mathcal A_m$, and $x=r/a\in[1,2)$. Hence
\[
 \frac{P_1(S)/\sqrt a+\sqrt a\,P_2(S)}
      {\sqrt{P_1(S)P_2(S)}}
 =\sqrt x+\frac1{\sqrt x}\le\frac3{\sqrt2}.
\]
Therefore
\begin{equation}\label{eq:dyadic-set-comparison}
\cutnorm{\by}
 \le\frac3{\sqrt2}\max_{a\in\mathcal A_m}K_a(\by).
\end{equation}

\emph{Step 2: upper bound of $\E K_a(\by)$ through
comparison with $\mathcal G$ on a stretched grid.}
Let $q\ge1$ be an integer. Take independent $N(0,1)$ variables
$Z_{r,j}$ on $R_q:=\{1,\ldots,qm\}\times\{1,\ldots,m\}$ and put
\[
 y^{(q)}_{i,j}
 :=\frac{\sigma}{\sqrt q}\sum_{r=q(i-1)+1}^{qi}Z_{r,j},
 \qquad 1\le i,j\le m.
\]
The entries of $\by^{(q)}$ are independent $N(0,\sigma^2)$ variables,
so $\by^{(q)}$ has the same distribution as $\by$. For a nonempty
$S\subseteq V_m$, replace each vertex by a block of $q$ vertices
in the first coordinate direction:
\[
 T_q(S)
 :=\bigcup_{(i,j)\in S}
    \bigl(\{q(i-1)+1,\ldots,qi\}\times\{j\}\bigr).
\]
This stretches $S$ in the first coordinate direction, and 
\[
 P_1(T_q(S))=P_1(S),\qquad
 P_2(T_q(S))=qP_2(S),\qquad
 \sum_{v\in S}y_v^{(q)}
 =\frac{\sigma}{\sqrt q}\sum_{w\in T_q(S)}Z_w.
\]
Consequently,
\[
 \frac{\left|\sum_{v\in S}y_v^{(q)}\right|}
      {P_1(S)/\sqrt q+\sqrt q\,P_2(S)}
 =\sigma\frac{\left|\sum_{w\in T_q(S)}Z_w\right|}
       {\operatorname{Per}(T_q(S))},
\]
so $K_q(\by^{(q)})\le\sigma\mathcal G_{R_q}(\mathbf Z)$.
Extend $\mathbf Z$ independently from $R_q$ to the square
$\{1,\ldots,qm\}^2$. The maximum over that square is at least
$\mathcal G_{R_q}(\mathbf Z)$, since the perimeter of each set
counts all lattice edges leaving it. Lemma~\ref{lem:dw-perimeter},
with variance one and side length $qm$, therefore gives
\[
 \E K_q(\by)\le C\sigma\{\log(eqm)\}^{3/4}.
\]
Transposing coordinates exchanges $K_q$ and $K_{1/q}$, so they
have the same expectation. Each $a\in\mathcal A_m$ is an integer
$q\le2m$ or its reciprocal. Since $\log(eqm)\le C L_m$,
\begin{equation}\label{eq:uniform-balance-mean}
 \max_{a\in\mathcal A_m}\E K_a(\by)\le C\sigma L_m^{3/4}.
\end{equation}

\emph{Step 3: the Gaussian expectation and probabilistic bound.}
By~\eqref{eq:area-perimeter} and AM--GM,
\[
 \frac{\sqrt{|S|}}{P_1(S)/\sqrt a+\sqrt a\,P_2(S)}
 \le\frac{\sqrt{|S|}}{2\sqrt{P_1(S)P_2(S)}}\le\frac14.
\]
Thus every \(K_a\) is \(1/4\)-Lipschitz. By the Gaussian concentration inequality in moment-generating-function form~\cite[Theorem~5.5]{BLM}, for every \(s>0\),
\[
\E\exp\!\left\{s\bigl(K_a(\by)-\E K_a(\by)\bigr)\right\}
\le e^{\sigma^2s^2/32}.
\]
Writing \(X_a=K_a(\by)-\E K_a(\by)\), Jensen's inequality gives
\[
\begin{aligned}
\E\max_{a\in\mathcal A_m}X_a
&\le \frac1s\log\E e^{s\max_a X_a}\\
&\le \frac1s\log\!\left(\sum_{a\in\mathcal A_m}\E e^{sX_a}\right)\\
&\le \frac{\log|\mathcal A_m|}{s}+\frac{\sigma^2s}{32}.
\end{aligned}
\]
Optimizing over \(s>0\) yields
\[
\E\max_{a\in\mathcal A_m}
\bigl\{K_a(\by)-\E K_a(\by)\bigr\}
\le
\frac{\sigma}{4}\sqrt{2\log|\mathcal A_m|}
\le C\sigma L_m^{3/4},
\]
where the last inequality follows from $|\mathcal A_m|=2\lceil\log_2m\rceil+1$.
Combining it with~\eqref{eq:dyadic-set-comparison} and
\eqref{eq:uniform-balance-mean}, 
\[
 \E\cutnorm{\by}\le C_D\sigma L_m^{3/4}
\]
for a universal constant $C_D$. On the other hand, AM--GM gives
$\cutnorm{\by}\ge2\mathcal G(\by)$, so
\eqref{eq:dw-square-mean} gives the stated lower bound.

Finally, \eqref{eq:norm-lipschitz} and Gaussian concentration for Lipschitz functions~\cite[Theorem~5.6]{BLM} imply
\[
 \PP\{\cutnorm{\by}>\E\cutnorm{\by}+t\}\le e^{-2t^2/\sigma^2},
 \qquad t>0.
\]
Taking $t=\sigma\sqrt{\frac12\log(1/\delta)}$ proves
\eqref{eq:gaussian-D} with
$\kappa_{m,\delta}=C_D L_m^{3/4}+\sqrt{\frac12\log(1/\delta)}$.
\end{proof}

\subsection{Proofs of the auxiliary lemmas}

\begin{proof}[Proof of Lemma~\ref{lem:perimeter-preliminaries}]
If $S$ occupies $r$ rows and $c$ columns, then $|S|\le rc$,
$P_1(S)\ge2c$, and $P_2(S)\ge2r$, proving
\eqref{eq:area-perimeter}. Also,
$P_1(S)\le2|S|\le2mr\le mP_2(S)$, since every vertex in $S$
has at most two incident edges counted by $P_1(S)$.
Exchanging the coordinates proves the other half of
\eqref{eq:perimeter-ratio}.
\end{proof}

\begin{proof}[Proof of Lemma~\ref{lem:cutnorm-lipschitz}]
For each nonempty \(S\subseteq V_m\), define
$F_S(\by)
:=
\frac{\left|\sum_{v\in S}y_v\right|}
{\sqrt{P_1(S)P_2(S)}}.$
Since \(\cutnorm{\by}=\max_S F_S(\by)\), we have
\[
\big|\cutnorm{\by}-\cutnorm{\by'}\big|
\le
\max_{\varnothing\ne S\subseteq V_m}
\big|F_S(\by)-F_S(\by')\big|.
\]
In addition, for each \(S\), the Cauchy--Schwarz inequality and \eqref{eq:area-perimeter} give
\[
\big|F_S(\by)-F_S(\by')\big|
\le
\frac{\left|\sum_{v\in S}(y_v-y_v')\right|}
{\sqrt{P_1(S)P_2(S)}}
\le
\frac{\sqrt{|S|}}
{\sqrt{P_1(S)P_2(S)}}\|\by-\by'\|_2\leq \frac12\|\by-\by'\|_2.
\]
\end{proof}

\begin{proof}[Proof of Lemma~\ref{lem:dw-perimeter}]
If \(C_1,\ldots,C_r\) are the connected components of a nonempty
\(A\subseteq V_m\), then perimeter additivity and the triangle inequality give
\[
\frac{\left|\sum_{v\in A}y_v\right|}{\operatorname{Per}(A)}
\le
\max_{1\le j\le r}
\frac{\left|\sum_{v\in C_j}y_v\right|}
     {\operatorname{Per}(C_j)}.
\]
Hence it suffices to consider connected sets.

Write $\Lambda_R=[-R,R]^2\cap\Z^2$.
Proposition~2.2 of Ding--Wirth~\cite{DW} controls all simply
connected lattice animals contained in $\Lambda_R$, without requiring
them to contain a prescribed vertex. Equation~(11) of~\cite{DW}
shows that the absolute maximum is unchanged when all connected
lattice animals are allowed. Our perimeter counts all lattice edges
leaving a set, including edges leaving the containing box, and thus
agrees with the edge boundary $|\partial A|$ used there.
Embed $V_m$ in $\Lambda_m$ and extend the Gaussian field independently
outside $V_m$. Apply Proposition~2.2 to both signs of the normalized
field $\by/\sigma$; no additional union bound over starting vertices
is needed. Rescaling gives
\[
\PP\!\left\{
\mathcal G(\by)>C\sigma L_m^{3/4}+u
\right\}
\le 2e^{-u^2/(2\sigma^2)},
\qquad u>0.
\] Integrating this tail bound yields
\[
\E\mathcal G(\by)\le C\sigma L_m^{3/4}.
\]

For the reverse inequality, suppose first that $m\ge7$ and put
$\ell=\lfloor(m-1)/2\rfloor$. A translate of $\Lambda_\ell$ lies
inside $V_m$, and $\ell\ge3$, as required in Theorem~1.6 of~\cite{DW}.
Restrict the maximum defining $\mathcal G(\by)$ to connected lattice
animals in this subbox that contain its center. The rooted lower
bound in that theorem, applied to $\by/\sigma$, together with
$\log\ell\asymp L_m$, gives after rescaling
\[
\E\mathcal G(\by)\ge c\sigma L_m^{3/4}.
\]
For $4\le m\le6$, a singleton gives
$\E\mathcal G(\by)\ge\E|y_v|/4=\sigma\sqrt{2/\pi}/4$,
so the same bound holds after decreasing the universal constant $c$.
\end{proof}

\section{Proofs of Theorem~\ref{thm:elementwise} and Corollaries~\ref{cor:det-global} and~\ref{cor:elementwise-risk}}\label{sec:deterministic}

\subsection{Auxiliary definitions and lemmas}
For coordinate calculations in this appendix, use
$r_i=r_p$, $s_j=s_p$, and $d_{ij}=\min\{r_i,s_j\}$
for $p=(i,j)\in V_m$. Here $r_i$ depends only on the row index
and $s_j$ only on the column index. We also write
$\widehat x_{ij}^{\lambda,\bb}=\widehat x_p^{\lambda,\bb}$.

The deterministic argument is based on the following directional
cut condition. Let $\beta_1,\beta_2>0$ and suppose that
\begin{equation}\label{eq:noise-cut}
 \left|\sum_{v\in S}y_v\right|
 \le
 \beta_1P_1(S)+\beta_2P_2(S)
 \qquad
 \text{for every }S\subseteq V_m.
\end{equation}
\begingroup

We prove Theorem~\ref{thm:elementwise} using two auxiliary lemmas.
Lemma~\ref{lem:leakage} derives a bound for the part of a threshold set
inside a rectangle from~\eqref{eq:noise-cut}.
Lemma~\ref{lem:peeling} is a combinatorial result that uses such bounds
on nested rectangles to exclude vertices far from the boundary.
Applying the boundary-crossing inequality to the full grid gives the uniform estimate;
combining it with the propagation lemma gives the spatial estimate.
Only the second argument introduces a logarithmic factor.
In the proof of the theorem, the directional coefficients
$\beta_1,\beta_2$ are chosen using $\cutnorm{\by}$.
\endgroup

\subsubsection{Rectangle-crossing edges and the boundary-crossing inequality}\label{subsec:leakage}
Let $A\subseteq V_m$ be any set, let $R\subseteq V_m$ be a lattice rectangle, and put $S=A\cap R$. For $k=1,2$, define directly
\begin{equation}\label{eq:leakage-count}
 \ell_k(A,R)=\#\bigl\{\{v,w\}:v\in A\cap R,
     \ w\in(A\setminus R)\cup F_m,\ w-v=\pm \be_k\bigr\}.
\end{equation}
These counts include edges to the exterior frame and edges joining
vertices of $A$ on opposite sides of the rectangle boundary.

\begin{lemma}[Boundary-crossing inequality for a threshold set]\label{lem:leakage}
Assume the cut inequality~\eqref{eq:noise-cut} with $\beta_1,\beta_2>0$ and $\lambda\ge2\max(\beta_1,\beta_2)$. For any real $t$, put $A=\{v\in V_m:x_{+,v}>t\}$. If $R\subseteq V_m$ is any axis-aligned lattice rectangle with side lengths $a_1,a_2\ge1$, then
\begin{equation}\label{eq:leakage-inequality}
 \left(t-\frac{(\lambda-2\beta_1)}{a_1}
          -\frac{(\lambda-2\beta_2)}{a_2}\right)|A\cap R|
 \le 2\beta_1\ell_1(A,R)+2\beta_2\ell_2(A,R).
\end{equation}
\end{lemma}

\subsubsection{A nested-rectangle propagation lemma}\label{subsec:peeling}
The next lemma is purely combinatorial. It shows that if, on every nested
rectangle, the number of vertices of \(A\) remaining inside is controlled
by the corresponding boundary-crossing counts, then
\(A\) cannot extend too far into the interior. Applied to the threshold
set from Lemma~\ref{lem:leakage}, this converts the boundary-crossing bound into a
location-dependent exclusion estimate.

\begin{lemma}[Interior exclusion by nested-rectangle propagation]
\label{lem:peeling}
Let $m\ge4$, $J=\lfloor m/4\rfloor$, and $A\subseteq V_m$.
For
\[
 \mathbf h=(h_1,h_2)\in\{0,\ldots,J\}^2,
\]
define
\[
 R_{\mathbf h}
 =
 [h_1+1,m-h_1]\times[h_2+1,m-h_2],
 \qquad
 S_{\mathbf h}=A\cap R_{\mathbf h},
\]
and
\[
 \ell_k(\mathbf h):=\ell_k(A,R_{\mathbf h}),
 \qquad k=1,2.
\]
Suppose that, for some $c_1,c_2>0$,
\begin{equation}\label{eq:basic-recursion}
 |S_{\mathbf h}|
 \le
 c_1\ell_1(\mathbf h)+c_2\ell_2(\mathbf h)
 \qquad
 \text{for every }\mathbf h\in\{0,\ldots,J\}^2.
\end{equation}
Then $A$ contains no vertex $(i,j)$ satisfying
\begin{equation}\label{eq:peeling-exclusion}
 \frac{c_1}{r_i}+\frac{c_2}{s_j}
 \le \frac{1}{128L_m}.
\end{equation}
\end{lemma}

\begingroup

\subsection{Proof of Theorem~\ref{thm:elementwise}}\label{subsec:det-main-proof}
\begin{proof}[Proof of Theorem~\ref{thm:elementwise}]
First suppose $\cutnorm{\by}>0$ and fix $\lambda\ge\cutnorm{\by}$.
Write $L_m=\log(em)$. By~\eqref{eq:sandwich} and the identity
$\bx_-(\by)=-\bx_+(-\by)$, it suffices to establish the upper bounds
for $\bx_+$: the same arguments apply to $-\by$ because the directional
cut condition is unchanged by sign reversal.

For every $a>0$ and every nonempty $S\subseteq V_m$,
\begin{equation}\label{eq:amgm-cuts}
 \left|\sum_{v\in S}y_v\right|
 \le\cutnorm{\by}\sqrt{P_1(S)P_2(S)}
 \le\frac{\cutnorm{\by}}{2\sqrt a}P_1(S)
       +\frac{\cutnorm{\by}\sqrt a}{2}P_2(S).
\end{equation}
Thus the directional cut condition~\eqref{eq:noise-cut} holds with
\begin{equation}\label{eq:directional_parameters}
 \beta_1=\frac{\cutnorm{\by}}{2\sqrt a},\qquad
 \beta_2=\frac{\cutnorm{\by}\sqrt a}{2}.
\end{equation}

\emph{Uniform estimate.}
Take $a=1$, so $\beta_1=\beta_2=\cutnorm{\by}/2$ and
$\lambda\ge2\max(\beta_1,\beta_2)$.
For any real $t$, put $A=\{v\in V_m:x_{+,v}>t\}$.
Each vertex of $A$ is incident to at most one frame edge in each
coordinate direction, so $\ell_k(A,V_m)\le|A|$ for $k=1,2$.
Lemma~\ref{lem:leakage}, applied with $R=V_m$ and $a_1=a_2=m$, gives
\[
 \left(t-\frac{2(\lambda-\cutnorm{\by})}{m}\right)|A|
 \le2\cutnorm{\by}|A|.
\]
If $A\ne\varnothing$, division by $|A|$ yields
\[
 t\le\frac{2\lambda}{m}
       +2\left(1-\frac1m\right)\cutnorm{\by}.
\]
Consequently, $A=\varnothing$ whenever $t$ exceeds this quantity,
which bounds every coordinate of $\bx_+$.
Sign reversal and comparison with the extremal solutions therefore give
\begin{equation}\label{eq:uniform-cap}
\sup_{\bb}\norm{\widehat{\bx}^{\lambda,\bb}}_\infty
 \le\frac{2\lambda}{m}+2\left(1-\frac1m\right)\cutnorm{\by} \leq 4\lambda/m+2\cutnorm{\by}.
\end{equation}
The factor $1-1/m$ comes from this full-square comparison:
the exterior contribution $2\cutnorm{\by}$ is reduced by $2\cutnorm{\by}/m$
in the rectangle-bias term. In particular, to combine this bound with the spatial
estimate below, we further relax it to $4\lambda/m+2\cutnorm{\by}$,
so that both estimates have the same leading term.

\emph{Spatial estimate.}
Fix a vertex $(i,j)$, put $\Lambda=\lambda/\cutnorm{\by}\ge1$, and choose
\begin{equation}\label{eq:a-choice}
 a=\min\{\Lambda^2,\max\{\Lambda^{-2},s_j/r_i\}\}.
\end{equation}
We clip the value of $a$ so that $a\in[\Lambda^{-2},\Lambda^2]$, to make sure that both directional coefficients in \eqref{eq:directional_parameters}
are at most $\lambda/2$, satisfying the condition in Lemma~\ref{lem:leakage}. Set
\[
 \eta=256L_m\left(\frac{\beta_1}{r_i}+\frac{\beta_2}{s_j}\right),
 \qquad t=\frac{4\lambda}{m}+\eta,
 \qquad A=\{v\in V_m:x_{+,v}>t\}.
\]
Every rectangle $R_{\bh}$ in Lemma~\ref{lem:peeling} has two side lengths
$a_k=m-2h_k\ge m/2$ for $k=1,2$. Since $\lambda\ge2\max(\beta_1,\beta_2)$,
\[
 \frac{\lambda-2\beta_1}{a_1}
 +\frac{\lambda-2\beta_2}{a_2}
 \le\lambda\left(\frac1{a_1}+\frac1{a_2}\right)
 \le\frac{4\lambda}{m}.
\]
Lemma~\ref{lem:leakage} therefore gives, for every such rectangle,
\[
 \eta|A\cap R_{\bh}|
 \le2\beta_1\ell_1(A,R_{\bh})+2\beta_2\ell_2(A,R_{\bh}).
\]
With $c_k=2\beta_k/\eta$, this is the hypothesis~\eqref{eq:basic-recursion}
of Lemma~\ref{lem:peeling}. Moreover,
\[
 \frac{c_1}{r_i}+\frac{c_2}{s_j}
 =\frac2\eta\left(\frac{\beta_1}{r_i}+\frac{\beta_2}{s_j}\right)
 =\frac1{128L_m}.
\]
That lemma gives $(i,j)\notin A$, so $x_{+,(i,j)}\le t$.
The same sign-reversal and comparison argument as above yields
\[
 |\widehat x_{ij}^{\lambda,\bb}|
 \le\frac{4\lambda}{m}
 +128L_m\cutnorm{\by}\left(\frac1{r_i\sqrt a}+\frac{\sqrt a}{s_j}\right).
\]
Without the constraint, the expression in parentheses is minimized at
$\sqrt a=\sqrt{s_j/r_i}$, where it equals $2/\sqrt{r_i s_j}$.
Clipping the value of $\sqrt{a}$ to $[\Lambda^{-1},\Lambda]$ can increase its first
summand only at the upper endpoint, where that summand is
$1/(\Lambda r_i)$, and its second only at the lower endpoint, where it
is $1/(\Lambda s_j)$. Hence
\[
 \frac1{r_i\sqrt a}+\frac{\sqrt a}{s_j}
 \le\frac2{\sqrt{r_i s_j}}+\frac1{\Lambda r_i}+\frac1{\Lambda s_j}
 \le\frac2{\sqrt{r_i s_j}}+\frac2{\Lambda d_{ij}}.
\]
Substitution gives
\begin{equation}\label{eq:before-optimize}
 |\widehat x_{ij}^{\lambda,\bb}|
 \le\frac{4\lambda}{m}
 +256L_m\left\{\frac{\cutnorm{\by}}{\sqrt{r_i s_j}}
                  +\frac{\cutnorm{\by}^2}{\lambda d_{ij}}\right\}.
\end{equation}
Taking the smaller of this estimate and the relaxed uniform estimate in \eqref{eq:uniform-cap} proves
\eqref{eq:det-main}. Each argument is valid for every vertex, every
exterior vector, and every $\lambda$ satisfying the theorem's hypotheses,
so the conclusion is simultaneous in all of them.

If $\cutnorm{\by}=0$, the singleton bounds imply $\by=0$.
For $A=\{v:x_{+,v}>t\}$ and $R=V_m$, the summed-optimality
inequality~\eqref{eq:sum-kkt} and the geometric estimate~\eqref{eq:pk_lk}
give
\[
 \left(t-\frac{2\lambda}{m}\right)|A|\le0.
\]
The same comparison and sign-reversal argument therefore yields
$\sup_{\bb}\norm{\widehat{\bx}^{\lambda,\bb}}_\infty
\le2\lambda/m$, which implies~\eqref{eq:det-main}.
\end{proof}
\endgroup

\subsection{Proof of Corollary~\ref{cor:det-global}}\label{subsec:det-summation}
\begin{proof}[Proof of Corollary~\ref{cor:det-global}]
Write $X_{ij}=\sup_{\bb}|\widehat x_{ij}^{\lambda,\bb}|$. Squaring the location-dependent bound in
Theorem~\ref{thm:elementwise} gives
\[
 X_{ij}^2\le C\left\{\frac{\lambda^2}{m^2}
       +\frac{L_m^2\cutnorm{\by}^2}{r_i s_j}
       +\frac{L_m^2\cutnorm{\by}^4}{\lambda^2d_{ij}^2}\right\}.
\]
The distances satisfy
\begin{equation}\label{eq:distance-sums}
 \sum_{i,j}\frac1{r_i s_j}
 =\left(\sum_i\frac1{r_i}\right)^2\le4L_m^2,
 \qquad
 \sum_{i,j}\frac1{d_{ij}^2}
 \le4m\sum_{d=1}^{\infty}d^{-2}\le8m.
\end{equation}
The second estimate uses the fact that at most $4m$ vertices have any
given value of $d_{ij}$. Finally,
$\sup_{\bb}\sum_{i,j}|\widehat x_{ij}^{\lambda,\bb}|^2
 \le\sum_{i,j}X_{ij}^2$.
Summing the elementwise bound and dividing by $m^2$ proves
\eqref{eq:det-global}.
\end{proof}

\subsection{Proof of Corollary~\ref{cor:elementwise-risk}}\label{subsec:gaussian-elementwise-risk}

\begin{proof}[Proof of Corollary~\ref{cor:elementwise-risk}]
\textbf{Proof of \eqref{eq:main-risk}}
Take $\delta=m^{-12}$ and put $C_*=C_D+\sqrt6$, where $C_D$ is
the universal constant in~\eqref{eq:gaussian-D}. Since $L_m\ge1$ and
$\log m\le L_m$,
\[
 \kappa_{m,m^{-12}}
 =C_D L_m^{3/4}+\sqrt{6\log m}
 \le C_*L_m^{3/4}.
\]
Choose a universal $C_0\ge\max\{1,C_*\}$. Every deterministic
$\lambda\ge C_0\sigma L_m^{3/4}$ then satisfies
$\lambda\ge\sigma\kappa_{m,m^{-12}}$. On the event
$\mathcal E=\{\cutnorm{\by}\le\sigma\kappa_{m,m^{-12}}\}$, the deterministic
elementwise theorem applies; Theorem~\ref{thm:gaussian-D} gives
$\PP(\mathcal E^c)\le m^{-12}$.

Fix $p=(i,j)$ and write $X_p=\sup_{\bb}|\widehat x_p^{\lambda,\bb}|$. On $\mathcal E$, apply~\eqref{eq:det-main} and replace $\cutnorm{\by}$ on its right-hand side by $\sigma\kappa_{m,m^{-12}}$. Squaring and using $(a+b+c)^2\le3(a^2+b^2+c^2)$ gives
\begin{equation}\label{eq:good-risk}
 \E[X_p^2\one_{\mathcal E}]
 \le C\left\{
 \frac{\lambda^2}{m^2}
 +\frac{\sigma^2L_m^{7/2}}{r_i s_j}
 +\frac{\sigma^4L_m^5}{\lambda^2d_{ij}^2}\right\}.
\end{equation}
By~\eqref{eq:crude}, $X_p^2\le2y_p^2+8\lambda^2$. Since $\E y_p^4=3\sigma^4$, Cauchy--Schwarz gives
\begin{align}
 \E[X_p^2\one_{\mathcal E^c}]
 &\le2(\E y_p^4)^{1/2}\PP(\mathcal E^c)^{1/2}
       +8\lambda^2\PP(\mathcal E^c)\notag\\
 &\le2\sqrt3\,\sigma^2m^{-6}+8\lambda^2m^{-12}.
 \label{eq:bad-risk}
\end{align}
Because $\lambda\ge C_0\sigma L_m^{3/4}\ge\sigma$, both terms are bounded by a universal multiple of $\lambda^2/m^2$. Adding~\eqref{eq:good-risk} and~\eqref{eq:bad-risk} proves~\eqref{eq:main-risk}.

\textbf{Proof of \eqref{eq:mse-risk}}

For every realization,
\[
 \sup_{\bb}\sum_{i,j}|\widehat x_{ij}^{\lambda,\bb}|^2
 \le\sum_{i,j}\sup_{\bb}|\widehat x_{ij}^{\lambda,\bb}|^2.
\]
Take expectations, apply~\eqref{eq:main-risk} at each vertex, and use
the distance sums in~\eqref{eq:distance-sums}. After division by
$m^2$, the result is~\eqref{eq:mse-risk}.
\end{proof}

\subsection{Proofs of the auxiliary lemmas}

\begin{proof}[Proof of Lemma~\ref{lem:leakage}]
Choose $M>\max\{\norm{\by}_\infty+2\lambda,t\}$.
By the stabilization statement in Lemma~\ref{lem:comparison},
$\bx_+=\widehat{\bx}^{\lambda,M\one_{F_m}}$.
Every exterior value therefore exceeds $t$.
Applying Lemma~\ref{lem:geo-binary} with \(C=S=A\cap R\), we have
\[
 N_+(S,t)=\ell_1(A,R)+\ell_2(A,R),
 \qquad
 N_-(S,t)=\operatorname{Per}(S)-N_+(S,t).
\]
Hence
\begin{equation}\label{eq:sum-kkt}
 0\le
 -t|S|+\sum_{v\in S}y_v
 -\frac{\lambda}{2}\operatorname{Per}(S)
 +\lambda\{\ell_1(A,R)+\ell_2(A,R)\}.
\end{equation}
Apply~\eqref{eq:noise-cut} and rearrange:
\begin{equation}\label{eq:before-perimeter}
 t|S|\le
 -\sum_{k=1}^2\left(\frac{\lambda}{2}-\beta_k\right)
        \{P_k(S)-2\ell_k(A,R)\}
 +2\sum_{k=1}^2\beta_k\ell_k(A,R).
\end{equation}
We claim that
\begin{equation}\label{eq:pk_lk}
 P_k(S)-2\ell_k(A,R)
 \ge -\frac{2|S|}{a_k},
 \qquad k=1,2.
\end{equation}
If \eqref{eq:pk_lk} holds, then multiplying it by $-(\lambda/2-\beta_k)\le0$ and substituting into
\eqref{eq:before-perimeter} proves~\eqref{eq:leakage-inequality}.

Now we prove \eqref{eq:pk_lk}. By symmetry, it suffices to consider \(k=1\), so the relevant
segments are the columns of \(R\). For each column \(I\) of $R$, let
\(p_I\) be its contribution to \(P_1(S)\), and let \(\ell_I\) be
the number of edges counted by \(\ell_1(A,R)\).

If \(I\not\subseteq S\), then each edge counted by \(\ell_I\)
is accompanied by another boundary edge of \(S\) farther inside
the column. When both endpoints of \(I\) belong to \(S\), these
two inner edges are distinct because \(I\) contains a vertex outside
\(S\). Hence
\[
 p_I\ge 2\ell_I.
\]
If \(I\subseteq S\), then \(p_I=2\) and \(\ell_I\le2\), so
\[
 p_I-2\ell_I\ge -2.
\]

Let \(N_1\) be the number of columns of \(R\) that are entirely
contained in \(S\). These columns are disjoint and each contains
\(a_1\) vertices, so
\[
 a_1N_1\le |S|.
\]
Summing over all columns gives
\[
 P_1(S)-2\ell_1(A,R)
 \ge -2N_1
 \ge -\frac{2|S|}{a_1}.
\]
The same argument with rows in place of columns proves \eqref{eq:pk_lk}.\end{proof}

\begin{proof}[Proof of Lemma~\ref{lem:peeling}]
Intuitively, if a set survives deep inside the grid, the hypothesis
forces it to have many occupied vertices in the surrounding layers.
Repeating this argument outward eventually requires more occupied
vertices than the grid contains.

By reflecting the grid in either coordinate if necessary, we may assume that
the vertex lies in the upper-left quadrant. Thus we may write
\[
 \bd=(d_1,d_2),
 \qquad
 1\le d_1,d_2\le \lceil m/2\rceil,
\]
where
\[
 d_1=r_i,\qquad d_2=s_j.
\]
It therefore suffices to show that
$\bd\notin A$ whenever
\begin{equation}\label{eq:peeling-exclusion-d}
 \frac{c_1}{d_1}+\frac{c_2}{d_2}
 \le \frac{1}{128L_m}.
\end{equation}

\medskip
\noindent\emph{Step 1: Crossing edges are controlled by the adjacent outer layer.}
Suppose $h_k\ge1$. The rectangle
$R_{\bh-\be_k}$ is obtained from $R_{\bh}$ by adding one layer in
direction $k$. Every edge counted by $\ell_k(\bh)$ has its endpoint
outside $R_{\bh}$ in this added layer. Moreover, different crossing
edges have different such endpoints. Therefore
\[
 \ell_k(\bh)
 \le
 |S_{\bh-\be_k}\setminus S_{\bh}|
 =
 |S_{\bh-\be_k}|-|S_{\bh}|,
\]
and hence
\begin{equation}\label{eq:c1-shell-count}
 \ell_k(\bh)
 \le
 |S_{\bh-\be_k}|-|S_{\bh}|,
 \qquad h_k\ge1.
\end{equation}

If $h_k=0$, there is no larger rectangle in direction $k$ within
$V_m$. Each edge counted by $\ell_k(\bh)$ is then incident to a
distinct vertex of $S_{\bh}$, so
\begin{equation}\label{eq:c1-boundary-count}
 \ell_k(\bh)\le |S_{\bh}|.
\end{equation}
For any $\mathcal I\subseteq\{1,2\}$, suppose that $h_k\ge1$ for
$k\in\mathcal I$ and $h_k=0$ for $k\notin\mathcal I$.
Combining~\eqref{eq:c1-shell-count} and~\eqref{eq:c1-boundary-count}
with~\eqref{eq:basic-recursion} gives, for
$\zeta:=1-\sum_{k\notin\mathcal I}c_k$,
\begin{equation}\label{eq:active-recursion}
 \left(
 \zeta+\sum_{k\in\mathcal I}c_k
 \right)|S_{\bh}|
 \le
 \sum_{k\in\mathcal I}c_k|S_{\bh-\be_k}|.
\end{equation}

\medskip
\noindent\emph{Step 2: Select the directions in which layers will be removed.}
Define the active directions by
\[
 \mathcal I
 :=
 \{k\in\{1,2\}:d_k\ge32L_m\}.
\]
These are the directions in which the vertex $\bd=(d_1,d_2)$ is
sufficiently far from the boundary to allow repeated removal of
layers while keeping $\bd$ inside the rectangle.

For an inactive direction $k\notin\mathcal I$,
\eqref{eq:peeling-exclusion} gives
\[
 c_k
 \le \frac{d_k}{128L_m}
 <\frac14.
\]
Hence
\[
 \zeta
 =
 1-\sum_{k\notin\mathcal I}c_k
 >\frac12.
\]

If $\mathcal I=\varnothing$, then \eqref{eq:active-recursion} with $\bh=(0,0)$ and
$S_{\bh}=A$ implies
\[
 |A|
 \le (c_1+c_2)|A|
 <\frac12|A|,
\]
which forces $A=\varnothing$.
Thus it remains only to consider $\mathcal I\ne\varnothing$.

\medskip
\noindent\emph{Step 3: Geometric decay in the active directions.}

We distinguish the two possible nonempty choices of $\mathcal I$.

\smallskip
\noindent\emph{Case 1: $\mathcal I$ contains one direction.}
By symmetry, WLOG assume $\mathcal I=\{1\}$. Then $h_2=0$, and
\eqref{eq:active-recursion} implies
\[
 (\zeta+c_1)|S_{(h_1,0)}|
 \le
 c_1|S_{(h_1-1,0)}|.
\]
Iterating gives
\[
 |S_{(h_1,0)}|
 \le
 m^2
 \left(\frac{c_1}{\zeta+c_1}\right)^{h_1}.
\]

Choose the largest inward depth allowed by both the construction and
the requirement that $\bd$ remain inside:
\[
 h_1^*:=\min\{d_1-1,J\},
 \qquad
 \bh^*:=(h_1^*,0).
\]
Then $\bd\in R_{\bh^*}$ and
\[
 h_1^*\ge\frac{d_1}{4}.
\]
Moreover, the exclusion condition \eqref{eq:peeling-exclusion} and $\zeta>1/2$ imply
\[
 \frac{\zeta}{c_1}\ge\frac{64L_m}{d_1}.
\]
Applying $\log(1+x)\geq x/(1+x)$,
\[
 \begin{aligned}
 h_1^*\log\left(1+\frac{\zeta}{c_1}\right)
 &\ge
 \frac{d_1}{4}
 \log\left(1+\frac{64L_m}{d_1}\right)\ge 2L_m.
 \end{aligned}
\]
Therefore
\[
 |S_{\bh^*}|
 \le m^2e^{-2L_m}
 <1.
\]
Thus $A\cap R_{\bh^*}=S_{\bh^*}=\varnothing$, and
$\bd\in R_{\bh^*}$ implies $\bd\notin A$.

\smallskip
\noindent\emph{Case 2: $\mathcal I=\{1,2\}$.}
Here $\zeta=1$, and \eqref{eq:active-recursion} becomes
\[
 (1+c_1+c_2)|S_{(h_1,h_2)}|
 \le
 c_1|S_{(h_1-1,h_2)}|
 +
 c_2|S_{(h_1,h_2-1)}|.
\]
When $h_1=0$ or $h_2=0$, one of the two powers below equals one,
so the bound follows from $|S_{(h_1,h_2)}|\le m^2$.
For $h_1,h_2\ge1$, the displayed recursion and induction on
$h_1+h_2$ yield
\[
 |S_{(h_1,h_2)}|
 \le
 m^2
 \left[
 \left(\frac{c_1}{1+c_1}\right)^{h_1}
 +
 \left(\frac{c_2}{1+c_2}\right)^{h_2}
 \right].
\]

Now set
\[
 h_k^*:=\min\{d_k-1,J\},
 \qquad k=1,2,
\]
and let $\bh^*:=(h_1^*,h_2^*)$. Then $\bd\in R_{\bh^*}$ and
\[
 h_k^*\ge\frac{d_k}{4},
 \qquad k=1,2.
\]
The exclusion condition \eqref{eq:peeling-exclusion} gives
\[
 \frac1{c_k}\ge\frac{128L_m}{d_k},
\]
and therefore
\[
 h_k^*\log\left(1+\frac1{c_k}\right)\ge2L_m,
 \qquad k=1,2.
\]
It follows that
\[
 |S_{\bh^*}|
 \le
 2m^2e^{-2L_m}
 <1.
\]
Again, $A\cap R_{\bh^*}=S_{\bh^*}=\varnothing$, and
$\bd\in R_{\bh^*}$ implies $\bd\notin A$.
\end{proof}

\section{Proofs of Theorem~\ref{thm:geometry} and Proposition~\ref{prop:boundary-attachment}}\label{sec:threshold-geometry}

Above the geometric threshold, the positive and negative threshold sets
form at most four corner-monotone regions, each with a constant sign,
surrounding a connected central band whose intersection with each
row or column is an interval. We first collect {three auxiliary lemmas}
and use them to prove Theorem~\ref{thm:geometry}.
The auxiliary lemma proofs follow the proofs of the theorem and
Proposition~\ref{prop:boundary-attachment}. The key corner argument
uses a single count of binary transitions to establish both monotonicity
and boundary contact; Proposition~\ref{prop:boundary-attachment} follows
separately from summed optimality equations.

\subsection{Auxiliary definitions and lemmas for corner geometry}
\label{subsec:geo-statistics}

We collect the definitions and {three auxiliary lemmas} used in the
proof of Theorem~\ref{thm:geometry}.
The auxiliary lemma proofs are given at the end of this appendix.
Throughout, $m\ge4$, $L_m=\log(em)$, and
$\widehat{\bx}=\widehat{\bx}^{\lambda,\bb}$ denotes the local estimator on $V_m$.
For coordinate calculations, use $r_i=r_p$, $s_j=s_p$, and
$\widehat x_{ij}=\widehat x_p$ when $p=(i,j)\in V_m$;
$r_i$ depends only on the row index and $s_j$ only on the column index.

\paragraph{Corner neighborhoods and corner monotonicity.}
Put $h=\lfloor m/2\rfloor$.
For a corner $c=(c_1,c_2)\in\{1,m\}^2$ and an integer
$1\le a\le h$, define
\[
 Q_{c,a}
 :=
 \bigl\{
 p\in V_m:
 |p_k-c_k|\le a-1,\ k=1,2
 \bigr\},
 \qquad
 H_c:=Q_{c,h}.
\]
Thus $Q_{c,a}$ is the $a\times a$ square attached to $c$.
Let $\mathcal V_c$ denote the quadrant associated with $c$ in
Section~\ref{subsec:main-geometry}; then
$Q_{c,a}\subseteq H_c\subseteq\mathcal V_c$.

For $p,q\in V_m$, write
\[
 p\preceq_c q
 \quad\Longleftrightarrow\quad
 |p_k-c_k|\le |q_k-c_k|,
 \qquad k=1,2.
\]
Thus $p\preceq_c q$ means that $p$ is coordinatewise closer to $c$
than $q$ is.
A set $S\subseteq V_m$ is corner-monotone from $c$ if
\[
 q\in S,\quad p\preceq_c q
 \quad\Longrightarrow\quad p\in S.
\]
In particular, every nonempty corner-monotone set contains $c$.

\begingroup

The following lemma converts the elementwise bound into the
localization needed for the geometric theorem. It confines large
errors to four corner regions before the row and column rearrangements
used below.

\begin{lemma}[Localization of threshold sets]\label{lem:excursion-localization}
Suppose $m\ge4$ and $0<\cutnorm{\by}\le\lambda$. If
\begin{equation}\label{eq:excursion-localization-level}
 t\ge2048\left\{
 \frac{\lambda}{m}
 +\frac{L_m^2\cutnorm{\by}^2}{\lambda}
 \right\},
 \qquad
 a_t:=
 \left\lceil
 \left(\frac{768L_m\cutnorm{\by}}{t}\right)^2
 \right\rceil,
\end{equation}
then $1\le a_t\le\lfloor m/2\rfloor$ and, for every $(i,j)\in V_m$,
\begin{equation}\label{eq:excursion-localization}
 r_i s_j\ge a_t
 \quad\Longrightarrow\quad
 |\widehat x_{ij}^{\lambda,\bb}|<t.
\end{equation}
\end{lemma}

For the corner-monotone sets defined above, put
\begin{equation}\label{eq:geo-width-statistic}
 \Phi(\by)
 :=\max_{c\in\{1,m\}^2}
   \max_{\substack{\varnothing\ne S\subseteq H_c\\
                    S\text{ corner-monotone from }c}}
   \frac{\left|\sum_{p\in S}y_p\right|}{\sqrt{|S|}}.
\end{equation}

\begin{lemma}[Corner-set fluctuations from directional fluctuations]
\label{lem:corner-fluctuation}
For every $\by\in\mathbb R^{V_m}$,
\begin{equation}\label{eq:corner-fluctuation-comparison}
 \Phi(\by)^2\le18L_m\cutnorm{\by}^2.
\end{equation}
\end{lemma}

\begin{lemma}[Corner monotonicity and boundary contact]
\label{lem:geo-binary-repair}
Let $Q=Q_{c,a}$, where $2\le a\le\lfloor m/2\rfloor$, and let
$\widehat{\bx}=\widehat{\bx}^{\lambda,\bb}$ with $\lambda>0$.
Suppose that $t>0$ and
\[
 \widehat x_w\le t
 \quad\text{for every }w\in V_m\setminus Q
 \text{ adjacent to a vertex of }Q.
\]
If
\begin{equation}\label{eq:geo-monotonicity-condition}
 \lambda>
 4\sqrt a\,\cutnorm{\by}+\frac{\Phi(\by)^2}{4t},
\end{equation}
then $A=\{p\in Q:\widehat x_p>t\}$ is corner-monotone from $c$.
Moreover,
\begin{equation}\label{eq:local-boundary-contact}
 p\in A,\quad w\in F_m,\quad p\sim w
 \quad\Longrightarrow\quad b_w>t.
\end{equation}
\end{lemma}

\subsection{Proof of Theorem~\ref{thm:geometry}}
\label{subsec:geo-theorem-proof}

\begin{proof}[Proof of Theorem~\ref{thm:geometry}]

Fix $\by$, $\lambda$, and $t$ satisfying
\eqref{eq:geo-corner-condition}, with the universal constant
chosen so that $C\ge8192$.
Fix any exterior array $\bb$, and abbreviate
$\widehat{\bx}=\widehat{\bx}^{\lambda,\bb}$.

If $\cutnorm{\by}=0$, Theorem~\ref{thm:elementwise} gives
\[
 |\widehat x_p|
 \le\frac{4\lambda}{m}
 <t
 \qquad(p\in V_m).
\]
Thus the large-error set is empty and $B_t=V_m$, proving the
conclusion in this case. Henceforth assume $\cutnorm{\by}>0$.

\medskip
\noindent\emph{Step 1: Localize the large-error set to four corner squares.}
The assumptions imply
\[
 \lambda\ge\cutnorm{\by},
 \qquad
 t\ge2048\left\{
 \frac{\lambda}{m}+\frac{L_m^2\cutnorm{\by}^2}{\lambda}
 \right\}.
\]
{Lemma~\ref{lem:excursion-localization}} therefore applies with
\[
 a:=
 \left\lceil
 \left(\frac{768L_m\cutnorm{\by}}{t}\right)^2
 \right\rceil.
\]
It gives
\begin{equation}\label{eq:geo-high-support}
 1\le a\le \lfloor m/2\rfloor,
 \qquad
 r_i s_j\ge a
 \quad\Longrightarrow\quad
 |\widehat x_{ij}|<t.
\end{equation}
Consequently,
\[
 \{p:|\widehat x_p|>t\}
 \subseteq
 \{(i,j):r_i s_j<a\}.
\]

If $a=1$, then $r_i s_j\ge1$ at every vertex, so again the
large-error set is empty and $B_t=V_m$.
Assume henceforth that $a\ge2$.

Since $r_i,s_j\ge1$, the inequality $r_i s_j<a$ implies
$r_i<a$ and $s_j<a$.
Every vertex in the large-error set therefore belongs to one of the four
corner squares $Q_{c,a}$.
Because $a\le \lfloor m/2\rfloor$, these squares are disjoint and each is
contained in its corresponding quadrant $\mathcal V_c$.

\medskip
\noindent\emph{Step 2: Verify the numerical condition in
Lemma~\ref{lem:geo-binary-repair}.}
The definition of $a$ gives
\[
 \begin{aligned}
 \cutnorm{\by}\sqrt a
 &\le \frac{768L_m\cutnorm{\by}^2}{t}+\cutnorm{\by}\\
 &\le \frac{768}{8192L_m}\lambda
       +\frac{1}{8192}\lambda\\
 &\le \frac{769}{8192}\lambda
 <\frac{\lambda}{8}.
 \end{aligned}
\]
Here we used
\[
 t\ge\frac{8192L_m^2\cutnorm{\by}^2}{\lambda},
 \qquad
 \lambda\ge8192\cutnorm{\by},
 \qquad L_m\ge1.
\]
Hence $4\sqrt a\,\cutnorm{\by}<\lambda/2$.

Lemma~\ref{lem:corner-fluctuation} also gives
\[
 \frac{\Phi(\by)^2}{4t}
 \le\frac{9L_m\cutnorm{\by}^2}{2t}
 \le\frac{9}{16384L_m}\lambda
 <\frac{\lambda}{2}.
\]
Adding these inequalities yields
\[
 4\sqrt a\,\cutnorm{\by}+\frac{\Phi(\by)^2}{4t}<\lambda,
\]
which verifies~\eqref{eq:geo-monotonicity-condition}.

\medskip

\noindent\emph{Step 3: Verify the inward-neighbor condition in
Lemma~\ref{lem:geo-binary-repair}.}
This condition ensures that the augmented rows and columns in
Lemma~\ref{lem:geo-binary-repair} have zero threshold indicators at
their inward ends.

Fix a corner $c$. We claim that
\[
 |\widehat x_w|<t
 \qquad
 \text{for every }w\in V_m\setminus Q_{c,a}\text{ adjacent to }Q_{c,a}.
\]
WLOG it suffices to consider $c=(1,1)$, so
$Q_{c,a}=\{1,\ldots,a\}^2$. Every neighbor of $Q_{c,a}$ on its inward sides then
has the form $(a+1,j)$ or $(i,a+1)$, with $1\le i,j\le a$.

For $w=(a+1,j)$,
\[
 r_{a+1}=\min\{a+1,m-a\}\ge a,
\]
since $a\le\lfloor m/2\rfloor$, while $s_j\ge1$. Hence
$r_{a+1}s_j\ge a$. The same argument gives
$r_i s_{a+1}\ge a$ for $w=(i,a+1)$.
Therefore~\eqref{eq:geo-high-support} yields
\[
 |\widehat x_w|<t
\]
for every inward neighbor of $Q_{c,a}$. By reflection, the same holds at
all four corners. Thus the inward-neighbor condition holds for both
$\widehat{\bx}$ and the sign-reversed solution $-\widehat{\bx}$.

\medskip
\noindent\emph{Step 4: Establish corner monotonicity of the positive and negative threshold sets.}
For each corner $c$, let
\[
 A_c^\pm:=\{p\in Q_{c,a}:\ \pm\widehat x_p>t\}.
\]
By Steps 2 and 3, Lemma~\ref{lem:geo-binary-repair} applies
to $\widehat{\bx}$. It also applies to the sign-reversed solution, since
\[
 \widehat{\bx}^{\lambda,-\bb}(-\by)
 =-\widehat{\bx}^{\lambda,\bb}(\by),
 \qquad
 \cutnorm{-\by}=\cutnorm{\by},
 \qquad
 \Phi(-\by)=\Phi(\by).
\]
Thus both $A_c^+$ and $A_c^-$ are corner-monotone from $c$.
Every nonempty corner-monotone set contains $c$, so these two
sets cannot both be nonempty. Consequently,
$A_c:=A_c^+\cup A_c^-$ is either empty or corner-monotone
with a single sign. Step 1 gives
\[
 \{p\in V_m:|\widehat x_p|>t\}
 =\bigcup_{c\in\{1,m\}^2}A_c,
\]
with each $A_c$ contained in its corresponding quadrant.
The boundary conclusion of the lemma, applied to both signs,
also gives $b_w>t$ whenever $\widehat x_v>t$ and $b_w<-t$
whenever $\widehat x_v<-t$, for every frame edge $\{v,w\}$
with $v\in V_m$ and $w\in F_m$.

\medskip
\noindent\emph{Step 5: Geometry of the central band.}
Let $B_t:=\{p\in V_m:|\widehat x_p|\le t\}$, and let
$\mathcal X_m$ be the union of the rows and columns indexed by
$\lfloor m/2\rfloor$ and $\lfloor m/2\rfloor+1$.
Each of these indices has boundary distance at least
$\lfloor m/2\rfloor\ge a$, so
\eqref{eq:geo-high-support} gives $\mathcal X_m\subseteq B_t$.

Within each quadrant, $B_t$ is the complement of a
corner-monotone set. Membership in $B_t$ is therefore preserved
under coordinate steps away from that quadrant's corner
that remain in the quadrant.
Every vertex of $B_t$ can consequently be joined to the
connected cross $\mathcal X_m$ while remaining in $B_t$.
Thus $B_t$ is nonempty and connected. The same monotonicity
shows that each half-row and half-column section of $B_t$
is an interval adjoining the cross. Since the cross belongs
to $B_t$, the two halves join to form an interval.
Hence each row or column intersects $B_t$ in an interval, and $B_t$
contains the middle row(s) and column(s).
\end{proof}

\subsection{Proof of Proposition~\ref{prop:boundary-attachment}}
\label{subsec:boundary-attachment}

\begin{proof}[Proof of Proposition~\ref{prop:boundary-attachment}]
Fix $\bb$, $t\ge0$, and a nonempty connected component $C$ of
$\{\widehat{\bx}^{\lambda,\bb}>t\}$. Suppose that $b_w\le t$ for every
exterior neighbor $w$ of $C$. Every interior neighbor outside $C$
also has estimated value at most $t$, by the definition of a component.
Thus each edge leaving $C$ has outward subgradient $+1$.
Summing the optimality equations and cancelling internal-edge terms gives
\[
 \sum_{v\in C}y_v
 =\sum_{v\in C}\widehat x_v^{\lambda,\bb}
  +\frac{\lambda}{2}\bigl(P_1(C)+P_2(C)\bigr)
 >t|C|+\frac{\lambda}{2}\bigl(P_1(C)+P_2(C)\bigr).
\]
On the other hand, the definition of $\cutnorm{\by}$ and AM--GM give
\[
 \sum_{v\in C}y_v
 \le\cutnorm{\by}\sqrt{P_1(C)P_2(C)}
 \le\frac{\lambda}{2}\bigl(P_1(C)+P_2(C)\bigr),
\]
contradicting $t\ge0$. Apply the same argument to $-\by,-\bb$, using
$\cutnorm{-\by}=\cutnorm{\by}$, to prove the negative assertion.
The argument is pathwise and holds for all $\bb$ and $t\ge0$ simultaneously.
\end{proof}

\subsection{Proofs of the auxiliary lemmas}

\subsubsection{Proof of threshold-set localization}

\begin{proof}[Proof of Lemma~\ref{lem:excursion-localization}]
The location-dependent bound in Theorem~\ref{thm:elementwise}, using $\min\{r_i,s_j\}\ge1$,
gives
\[
 |\widehat x_{ij}^{\lambda,\bb}|
 \le
 \frac{4\lambda}{m}
 +\frac{256L_m\cutnorm{\by}}{\sqrt{r_i s_j}}
 +\frac{256L_m\cutnorm{\by}^2}{\lambda}.
\]
If $r_i s_j\ge a_t$, the second term is at most $t/3$.
The first is at most $t/512<t/3$, and the third at most
$t/(8L_m)<t/3$, proving \eqref{eq:excursion-localization}.

By AM--GM,
$t\ge4096L_m\cutnorm{\by}/\sqrt m$, and hence
\[
 a_t
 \le\left\lceil\frac{9m}{256}\right\rceil
 \le\left\lfloor\frac m2\right\rfloor
 \qquad(m\ge4).
\]
\end{proof}

\subsubsection{Proof of the corner-set fluctuation bound}
\label{subsec:geo-fluctuation-proof}

\begin{proof}[{Proof of Lemma~\ref{lem:corner-fluctuation}}]
Fix a nonempty corner-monotone set $S\subseteq H_c$ and, by reflection,
assume $c=(1,1)$. Write its nonempty row lengths as
$\ell_1\ge\ell_2\ge\cdots\ge\ell_r\ge1$, so
$|S|=\sum_{i=1}^r\ell_i$.

Let $K=\lfloor\log_2 r\rfloor$ and, for $0\le k\le K$, set
$I_k=\{2^k,\ldots,\min(2^{k+1}-1,r)\}$ and
$S_k=\{(i,j)\in S:i\in I_k\}$. Each occupied row and column of $S_k$
is an interval, hence $P_1(S_k)=2\ell_{2^k}$ and
$P_2(S_k)=2|I_k|$. By the definition of $\cutnorm{\by}$,
\[
 \left|\sum_{p\in S_k}y_p\right|
 \le 2\cutnorm{\by}\sqrt{|I_k|\ell_{2^k}}.
\]

For $k\ge1$, $|I_k|\le2^k$, while the $2^{k-1}$ rows in $I_{k-1}$
all have length at least $\ell_{2^k}$; hence
$|I_k|\ell_{2^k}\le2|S_{k-1}|$. For $k=0$,
$|I_0|\ell_1=|S_0|$. Therefore
\[
 \sum_{k=0}^K |I_k|\ell_{2^k}\le3|S|.
\]
Using the triangle inequality and Cauchy--Schwarz,
\[
 \left|\sum_{p\in S}y_p\right|
 \le
 2\cutnorm{\by}\sum_{k=0}^K\sqrt{|I_k|\ell_{2^k}}
 \le
2\cutnorm{\by} \sqrt{(K+1)}\sqrt{\sum_{k=0}^K |I_k|\ell_{2^k}}\le
2\cutnorm{\by} \sqrt{(K+1)}\sqrt{3|S|}.
\]

Since $S\subseteq H_c$, we have $r\le\lfloor m/2\rfloor$, so
$K+1\le\log_2 m$. Hence
$12(K+1)\le(12/\log2)\log m<18L_m$. 
Since this bound holds for every corner $c$ and every nonempty
corner-monotone $S\subseteq H_c$, the definition of $\Phi(\by)$ gives
\[
 \Phi(\by)^2\le18L_m\cutnorm{\by}^2.
\]
\end{proof}

\subsubsection{Proof of corner monotonicity and boundary contact}
\label{subsec:geo-monotonicity-proof}

\begin{proof}[Proof of Lemma~\ref{lem:geo-binary-repair}]

Reflect the coordinates so that $c=(1,1)$ and $Q=[a]^2$. Let
\[
 A:=\{p\in Q:\widehat x_p>t\},
 \qquad
 u_{ij}:=\one_{\{(i,j)\in A\}}
       =\one_{\{\widehat x_{ij}>t\}},
 \quad 1\le i,j\le a.
\]
Extend this notation to the neighboring vertices of $Q$ by setting
\[
 u_{i0}:=\one_{\{b_{(i,0)}>t\}},\qquad
 u_{0j}:=\one_{\{b_{(0,j)}>t\}},
\]
and
\[
 u_{i,a+1}:=0,\qquad
 u_{a+1,j}:=0,
 \qquad 1\le i,j\le a.
\]
The last two assignments of zeros are valid because by the hypothesis of Lemma~\ref{lem:geo-binary-repair}, every neighbor of \(Q\) lying in \(V_m\setminus Q\) has estimated value at most \(t\). 

\emph{A single count of violations.}
Define
\[
 K:=
 \sum_{i=1}^a\sum_{j=0}^a
   (u_{i,j+1}-u_{ij})_+
 +
 \sum_{j=1}^a\sum_{i=0}^a
   (u_{i+1,j}-u_{ij})_+ .
\]
Thus $K$ is the total number of $0\to1$ transitions in the rows and
columns when read away from the corner $(1,1)$, including transitions
across the outward boundary of $Q$.

If $K=0$, then
$u_{i,j+1}\le u_{ij}$ and $u_{i+1,j}\le u_{ij},$ throughout the augmented rows and columns, so $A$ is corner-monotone from $(1,1)$. Moreover, if an occupied
vertex lies on the left or top boundary of $Q$, its adjacent outward
exterior indicator must equal one. Thus $K=0$ would also give the
boundary-contact conclusion \eqref{eq:local-boundary-contact}.

We next relate $K$ to the edge counts in
Lemma~\ref{lem:geo-binary}. Note that to define $ N_\pm(A,t)$, we continue to use the edge sets $E_{\rm int}$ and $E_{\rm fr}$ of the full square $V_m$
from Section~\ref{sec:fluctuation}. Since $Q=[a]^2$ is the upper-left corner square of
$V_m$, every edge in $E_{\rm int}\cup E_{\rm fr}$ leaving
$A\subseteq Q$ appears in one of the augmented rows or columns:
edges on the left and top sides meet the frame $F_m$, while edges on
the right and bottom sides meet $V_m\setminus Q$.

Along any augmented row or column,
with labels $u_0,\ldots,u_{a+1}$ and $u_{a+1}=0$,
\[
 \sum_{j=0}^{a}|u_{j+1}-u_j|
 =
 u_0+2\sum_{j=0}^{a}(u_{j+1}-u_j)_+.
\]
For each line, the contribution to $N_--N_+$ is
\[
N_--N_+= \sum_{j=0}^{a}|u_{j+1}-u_j|-u_0.
\]
Summing this identity over all rows and columns 
gives
\[
 N_-(A,t)-N_+(A,t)=2K.
\]

Indeed, an edge leaving $A$ contributes $+1$ to the left-hand side
when its other endpoint has threshold indicator zero, and $-1$ when
that indicator is one.

Applying Lemma~\ref{lem:geo-binary} with $C=A$ therefore yields
\begin{equation}\label{eq:geo-binary-identity}
 t|A|-\sum_{p\in A}y_p+\lambda K\le0.
\end{equation}
If $A=\varnothing$, the conclusion is trivial, so assume henceforth
that $A\ne\varnothing$.

\emph{Coordinatewise rearrangement toward the corner.}
Every horizontal or vertical interval $I\subseteq Q$ has directional
perimeters $2$ and $2|I|$. Hence
\begin{equation}\label{eq:geo-cut-interval}
 \left|\sum_{p\in I}y_p\right|
 \le 2\cutnorm{\by}\sqrt{|I|}.
\end{equation}

Let
\[
 H(\bu):=\sum_{i=1}^a\sum_{j=1}^{a-1}(u_{i,j+1}-u_{ij})_+,
 \qquad
 V(\bu):=\sum_{i=1}^{a-1}\sum_{j=1}^a(u_{i+1,j}-u_{ij})_+.
\]
Since $K$ also counts possible $0\to1$ transitions across the outward
boundary,
\[
 H(\bu)+V(\bu)\le K.
\]

First compress the ones in each row to the left, preserving the number
of ones in that row. Let the resulting array be $\bu^{\rm H}$.
A run beginning in column $1$ does not move. Every occupied run
that moves therefore begins at a distinct internal horizontal
$0\to1$ transition, so at most $H(\bu)$ runs move. If a run $I$ is
moved to an interval $J$, then $|I|=|J|\le a$, and
\eqref{eq:geo-cut-interval} gives
\[
 \left|\sum_{p\in I}y_p-\sum_{p\in J}y_p\right|
 \le 4\cutnorm{\by}\sqrt a.
\]
Therefore
\begin{equation}\label{eq:compress1}
 \bigl|\ip{\by_Q}{\bu-\bu^{\rm H}}\bigr|
 \le 4\cutnorm{\by}\sqrt a\,H(\bu).
\end{equation}

Write $N_i=\sum_j u_{ij}$. Between two consecutive compressed rows,
the number of vertical $0\to1$ transitions is
\[
 (N_{i+1}-N_i)_+
 \le \sum_{j=1}^a (u_{i+1,j}-u_{ij})_+.
\]
Thus
\[
 V(\bu^{\rm H})\le V(\bu).
\]

Now compress each column of $\bu^{\rm H}$ upward, and let
$\bu^{\rm C}$ be the resulting array. By the same argument,
\begin{equation}\label{eq:compress2}
 \bigl|\ip{\by_Q}{\bu^{\rm H}-\bu^{\rm C}}\bigr|
 \le 4\cutnorm{\by}\sqrt a\,V(\bu^{\rm H})
 \le 4\cutnorm{\by}\sqrt a\,V(\bu).
\end{equation}
If
\[
 M_j:=\sum_{i=1}^a u^{\rm H}_{ij},
\]
then
\[
 M_j=\#\{i:N_i\ge j\},
\]
so $M_1\ge\cdots\ge M_a$. Hence the support
\[
 S:=\{(i,j):u^{\rm C}_{ij}=1\}
\]
is corner-monotone from $(1,1)$. Both coordinate compressions preserve
cardinality, so $|S|=|A|$. By the triangle inequality, \eqref{eq:compress1}, and \eqref{eq:compress2}, 
\begin{equation}\label{eq:geo-sorting-bounds}
 \left|\sum_{p\in A}y_p-\sum_{p\in S}y_p\right|
 \le  \bigl|\ip{\by_Q}{\bu-\bu^{\rm H}}\bigr|+ \bigl|\ip{\by_Q}{\bu^{\rm H}-\bu^{\rm C}}\bigr|\le 
 4\cutnorm{\by}\sqrt a\,[H(\bu)+V(\bu)]
 \le
 4\cutnorm{\by}\sqrt a\,K.
\end{equation}

Since $S\subseteq Q\subseteq H_c$ is nonempty and corner-monotone,
\eqref{eq:geo-width-statistic} gives
\[
 \sum_{p\in A}y_p
 \le
 \Phi(\by)\sqrt{|A|}
 +4\cutnorm{\by}\sqrt a\,K.
\]
Combining this with \eqref{eq:geo-binary-identity},
\[
 \begin{aligned}
 \bigl(\lambda-4\cutnorm{\by}\sqrt a\bigr)K
 &\le \Phi(\by)\sqrt{|A|}-t|A|\\
 &\le \frac{\Phi(\by)^2}{4t}\\
 &<\lambda-4\cutnorm{\by}\sqrt a,
 \end{aligned}
\]
where the last inequality follows from
\eqref{eq:geo-monotonicity-condition}. Thus $K<1$, and hence $K=0$, proving the lemma.
\end{proof}

\endgroup

\section{Proofs of Theorem~\ref{thm:piecewise-constant} and Corollaries~\ref{cor:gaussian-signal-elementwise}, \ref{cor:gaussian-signal-geometry}, and~\ref{cor:rectangle-gaussian}}\label{sec:rectangular-proofs}
We prove the deterministic error bound and the Gaussian deviation
and risk bounds for signals on rectangular partitions.
Lemma~\ref{lem:patch-restriction} applies local bounds that hold
uniformly over exterior values to restrictions of a global estimator,
including squares that meet its free boundary.
Theorem~\ref{thm:piecewise-constant} and the rectangle lemma below
require no distributional assumption.
We group the deterministic error and Gaussian risk proofs, which
share the square-covering argument, before the Gaussian deviation
proofs. The proofs of the restriction and rectangle lemmas come last.

\subsection{Auxiliary definitions and a rectangle bound}\label{subsec:local-patch}
Use the free-boundary estimator~\eqref{eq:whole-image-model}.
For a square $Q=\{i_0+1,\ldots,i_0+m\}\times\{j_0+1,\ldots,j_0+m\}\subseteq[{n}]^2$ and $p=(i_0+i,j_0+j)\in Q$, write
\[
 \begin{gathered}
 r_Q(p)=\min\{i,m+1-i\},\qquad s_Q(p)=\min\{j,m+1-j\},\\
 d_Q(p)=\min\{r_Q(p),s_Q(p)\},\qquad L_m=\log(em).
 \end{gathered}
\]

\paragraph{Bounds on squares with constant signal.}
For a fixed square $Q$ of side $m\ge4$ on which the true signal is
constant with value $c_Q$, Lemma~\ref{lem:patch-restriction} and
translation equivariance transfer the deterministic elementwise
bound of Theorem~\ref{thm:elementwise} and the Gaussian risk bound
of Corollary~\ref{cor:elementwise-risk} to
$\widehat{\btheta}^{\lambda}|_Q-c_Q\one_Q$.
Use the local size and distances defined above, under the respective
conditions on $\lambda$.

\begin{lemma}[elementwise error on rectangles with constant signal]
\label{lem:rectangle-elementwise}
Let $\by=\btheta+\beps$ on $[n]^2$, with $n\ge4$ and arbitrary
errors $\beps\in\mathbb R^{[n]^2}$, and suppose that $\btheta$ is piecewise constant
on a rectangular partition. For $p$, let $R(p)$ be
its rectangle and let $m_p$ be the shorter side length of $R(p)$.
Let $r_p,s_p$ be its distances to the nearer horizontal and vertical
sides of $R(p)$, counting boundary vertices at distance one, and put
$d_p=\min\{r_p,s_p\}$ and $L_{m_p}=\log(em_p)$.
There is a universal constant $C>0$ such that the free-boundary
estimator~\eqref{eq:whole-image-model} satisfies, simultaneously over
all vertices $p$ and all $\lambda\ge\cutnorm{\beps}$,
\begin{equation}\label{eq:rectangle-elementwise-det}
 \bigl|\widehat\theta_p^\lambda-\theta_p\bigr|
 \le C\left[
 \frac{\lambda}{m_p}
 +\min\left\{
 \cutnorm{\beps},\,
 \frac{L_{m_p}\cutnorm{\beps}}{\sqrt{r_ps_p}}
 +\frac{L_{m_p}^2\cutnorm{\beps}^2}{\lambda d_p}
 \right\}\right].
\end{equation}
\end{lemma}

\subsection{Proof of Theorem~\ref{thm:piecewise-constant}}\label{subsec:piecewise-constant}
The proof covers each rectangle by squares, then uses the
deterministic mean squared error bound on each square.
\begin{proof}[Proof of Theorem~\ref{thm:piecewise-constant}]
For an axis-aligned square $Q\subseteq[{n}]^2$, let
$\mathcal D_Q(\beps|_Q)$ denote the directional fluctuation
parameter after translating $Q$ to a square grid of the same size.
Every nonempty $S\subseteq Q$ has the same directional perimeters
whether computed in $Q$ or in $[{n}]^2$: both count all 
edges leaving $S$, including edges to the corresponding exterior layer.
Consequently,
\begin{equation}\label{eq:patch-D-restriction}
 \mathcal D_Q(\beps|_Q)
 \le\cutnorm{\beps}.
\end{equation}

Suppose $Q$ has side $a\ge4$ and $\btheta=c_Q$ on $Q$.
Lemma~\ref{lem:patch-restriction} and translation equivariance
identify $\widehat{\btheta}^\lambda|_Q-c_Q\one_Q$ with the local
estimator with observations $\beps|_Q$.
Corollary~\ref{cor:det-global}, together with
\eqref{eq:patch-D-restriction} and $L_a\le L_n$, therefore gives
\begin{equation}\label{eq:patch-det-sum}
 \sum_{p\in Q}|\widehat\theta_p^\lambda-\theta_p|^2
 \le C\left\{
  \lambda^2+{L_n}^4\cutnorm{\beps}^2
  +\frac{{L_n}^2\cutnorm{\beps}^4a}{\lambda^2}
 \right\}.
\end{equation}

If $a_i\ge4$, cover $R_i$ by
$q_i=\lceil b_i/a_i\rceil\le2b_i/a_i$ squares of side $a_i$ contained
in $R_i$: place successive squares along its longer side and place the
last square against the far end.  Only the last square may overlap its
predecessor.  Summing~\eqref{eq:patch-det-sum} over these squares is an
upper bound for the squared error on $R_i$, since all squared errors are
nonnegative.  The result is the $i$-th summand of
\eqref{eq:rectangular-det}, up to a universal constant.

If $a_i\le3$, coordinatewise optimality of
\eqref{eq:whole-image-model}, and the fact that each vertex has at most
four neighbors, give
$|\widehat\theta_p^\lambda-y_p|\le2\lambda$.
Singleton sets in the definition of $\cutnorm{\beps}$ give
$|\varepsilon_p|\le2\cutnorm{\beps}$, and hence
\[
 |\widehat\theta_p^\lambda-\theta_p|
 \le2\lambda+2\cutnorm{\beps}\le4\lambda.
\]
Since $a_ib_i=a_i^2(b_i/a_i)\le9b_i/a_i$, we obtain
\[
 \sum_{p\in R_i}|\widehat\theta_p^\lambda-\theta_p|^2
 \le144\lambda^2\frac{b_i}{a_i}.
\]
The first term in~\eqref{eq:rectangular-det} therefore absorbs these
rectangles as well.  Summing over the partition and dividing by $N$
proves~\eqref{eq:rectangular-det}.

Finally, when $b_i/a_i\le\rho$,
\[
 \sum_{i=1}^k\frac{b_i}{a_i}\le\rho k,
 \qquad
 \sum_{i=1}^k b_i
 \le\sqrt{k\sum_{i=1}^k b_i^2}
 \le\sqrt{\rho k\sum_{i=1}^k a_i b_i}
 =\sqrt{\rho kN}.
\]
Substitution and $\sqrt\rho\le\rho$ give~\eqref{eq:rectangular-det-aspect}
with the same universal constant $C$.
\end{proof}

\subsection{Proof of Corollary~\ref{cor:rectangle-gaussian}}\label{subsec:gaussian-rectangles}

\begin{proof}[Proof of Corollary~\ref{cor:rectangle-gaussian}]
Write $\beps=\by-\btheta$, so its entries are independent
$N(0,\sigma^2)$ variables. By
Theorem~\ref{thm:gaussian-D}, the event
\[
 \mathcal E=\{\cutnorm{\beps}\le K\},
 \qquad
 K=\sigma\kappa_{{n},{n}^{-12}}
 \le C_*\sigma {L_n}^{3/4},
\]
has $\PP(\mathcal E^c)\le {n}^{-12}$, where $C_*=C_D+\sqrt6$.
Choose $C_0\ge\max\{1,C_*\}$. For every fixed deterministic
$\lambda\ge C_0\sigma {L_n}^{3/4}$, the deterministic rectangle theorem
applies on $\mathcal E$. Using $\cutnorm{\beps}\le K$ in
Theorem~\ref{thm:piecewise-constant} yields
\[
 \frac1N\E\left[
  \|\widehat{\btheta}^\lambda-\btheta\|_2^2\one_{\mathcal E}\right]
 \le\frac{C}{N}\sum_{i=1}^k
 \left\{
  (\lambda^2+\sigma^2{L_n}^{11/2})\frac{b_i}{a_i}
  +\frac{\sigma^4{L_n}^5b_i}{\lambda^2}
 \right\}.
\]
On the exceptional event, coordinatewise optimality gives
$|\widehat\theta_p^\lambda-\theta_p|^2
 \le2\varepsilon_p^2+8\lambda^2$.
Since $\E\varepsilon_p^4=3\sigma^4$, Cauchy--Schwarz implies
\[
 \frac1N\E\left[
  \|\widehat{\btheta}^\lambda-\btheta\|_2^2\one_{\mathcal E^c}\right]
 \le2\sqrt3\,\sigma^2{n}^{-6}+8\lambda^2{n}^{-12}
 \le C\frac{\lambda^2}{N}.
\]
Here the last inequality uses $\lambda\ge C_0\sigma{L_n}^{3/4}$.
Because $\sum_i b_i/a_i\ge1$, this term is absorbed by the first
summand of the claimed bound.  Adding the contributions proves
\eqref{eq:rectangular-risk}.

Under $b_i/a_i\le\rho$,
\[
 \sum_{i=1}^k\frac{b_i}{a_i}\le\rho k,
 \qquad
 \sum_{i=1}^k b_i
 \le\sqrt{k\sum_{i=1}^k b_i^2}
 \le\sqrt{\rho k\sum_{i=1}^k a_i b_i}
 =\sqrt{\rho kN}.
\]
Applying~\eqref{eq:rectangular-risk} and using $\sqrt\rho\le\rho$
therefore gives, for every fixed deterministic
$\lambda\ge C_0\sigma{L_n}^{3/4}$,
\begin{equation}\label{eq:rectangular-risk-aspect}
 \frac1N\E\norm{\widehat{\btheta}^\lambda-\btheta}_2^2
 \le C\rho\left\{\frac{k\lambda^2}{N}
       +\frac{k\sigma^2{L_n}^{11/2}}{N}
       +\frac{\sigma^4{L_n}^5\sqrt{k}}{\lambda^2\sqrt N}\right\}.
\end{equation}
Balancing the first and third terms gives the stated choice of $\lambda$. This choice
satisfies the required lower bound because ${L_n}\ge1$ and $k\le N$, so
$\lambda\ge C_0\sigma{L_n}^{3/4}$ for a sufficiently large fixed $C_{\rm reg}$.
Substitution proves~\eqref{eq:rectangular-risk-tuned}.
\end{proof}

\subsection{Proofs of Corollaries~\ref{cor:gaussian-signal-elementwise}
and~\ref{cor:gaussian-signal-geometry}}
\label{subsec:gaussian-signal-proofs}

\begin{proof}[Proof of Corollaries~\ref{cor:gaussian-signal-elementwise}
and~\ref{cor:gaussian-signal-geometry}]
Fix $0<\delta<1$. By Theorem~\ref{thm:gaussian-D}, the event
\[
 \mathcal E_\delta
 =\{\cutnorm{\beps}\le\sigma\kappa_{n,\delta}\}
\]
has probability at least $1-\delta$. We prove both conclusions on
this event.

\emph{Elementwise error.}
For every $\lambda\ge\sigma\kappa_{n,\delta}$, we have
$\lambda\ge\cutnorm{\beps}$ on $\mathcal E_\delta$.
Lemma~\ref{lem:rectangle-elementwise} therefore applies
simultaneously at every vertex. Using $L_{m_p}\le L_n$ and bounding
the minimum in~\eqref{eq:rectangle-elementwise-det} by its second
argument gives
\[
 \bigl|\widehat\theta_p^\lambda-\theta_p\bigr|
 \le C\left\{
 \frac{\lambda}{m_p}
 +\frac{\sigma L_n\kappa_{n,\delta}}{\sqrt{r_ps_p}}
 +\frac{\sigma^2L_n^2\kappa_{n,\delta}^2}
 {\lambda d_p}
 \right\},
\]
which proves Corollary~\ref{cor:gaussian-signal-elementwise}.
The lemma also covers rectangles whose shorter side has length at most three.

\emph{Large-error geometry.}
Let $Q$ be any square in Corollary~\ref{cor:gaussian-signal-geometry},
with side $m\ge4$ and constant signal value $c_Q$.
By Lemma~\ref{lem:patch-restriction} and translation equivariance,
the centered restriction $\widehat{\btheta}^{\lambda}|_Q-c_Q\one_Q$
is a local estimator with observations $\beps|_Q$.
Every nonempty subset of $Q$ has the same two directional
perimeters when viewed in $Q$ or in $[n]^2$, since both
definitions count edges to the exterior frame. Hence, on
$\mathcal E_\delta$,
\[
 \mathcal D_Q(\beps|_Q)
 \le\cutnorm{\beps}
 \le\sigma\kappa_{n,\delta}.
\]
Choose $C_0$ and $C_1$ at least as large as the constant in
Theorem~\ref{thm:geometry}. Since $L_m\le L_n$,
condition~\eqref{eq:gaussian-signal-geometry} implies that theorem's
conditions on $\lambda$ and $t$ for the square $Q$.
Its corner-monotonicity, constant-sign, and central-band conclusions
give Corollary~\ref{cor:gaussian-signal-geometry}.
These exterior values may depend on the observations and
on $\lambda$; the theorem is uniform over those values.
The single event $\mathcal E_\delta$ thus controls all vertices,
squares, regularization parameters, and levels in the two corollaries,
without a union bound over the partition or the squares.

Finally, for $\delta=n^{-10}$,
\[
 \kappa_{n,n^{-10}}
 =C_D L_n^{3/4}+\sqrt{5\log n}
 \le(C_D+\sqrt5)L_n^{3/4}.
\]
The stochastic terms in~\eqref{eq:gaussian-signal-elementwise} are
therefore bounded by universal constants times
$\sigma L_n^{7/4}/\sqrt{r_ps_p}$ and
$\sigma^2L_n^{7/2}/(\lambda d_p)$, respectively,
and the second term in the level bound is at most a universal
constant times $\sigma^2L_n^{7/2}/\lambda$.
Since $N=n^2$ and $n\ge4$, we have $L_n\le\log N$.
These bounds imply Theorems~\ref{thm:informal-piecewise}
and~\ref{thm:informal-piecewise-geometry} with their stated
probability and conditions on $\lambda$.
\end{proof}

\subsection{Proofs of the restriction and rectangle lemmas}
\label{subsec:rectangular-auxiliary-proofs}

\begin{proof}[Proof of Lemma~\ref{lem:patch-restriction}]
Hold the estimated values $\widehat\theta_w^\lambda$ fixed for
$w\in[n]^2\setminus Q$. As a function of the coordinates on $Q$,
the global objective~\eqref{eq:whole-image-model}, up to terms
independent of these coordinates, is
\[
 \sum_{v\in Q}(y_v-x_v)^2
 +\lambda\!\!\sum_{\substack{\{v,w\}\subseteq Q\\v\sim w}}
       |x_v-x_w|
 +\lambda\!\!\sum_{\substack{v\in Q,\;w\in[n]^2\setminus Q\\v\sim w}}
       |x_v-\widehat\theta_w^\lambda|.
\]
Global optimality implies that $\widehat{\btheta}^{\lambda}|_Q$
minimizes this expression.

For each frame edge $\{v,w\}$ with $v\in Q$ and $w\notin[n]^2$,
add the term $\lambda|x_v-\widehat\theta_v^\lambda|$.
These added terms are nonnegative and vanish at
$\widehat{\btheta}^{\lambda}|_Q$, so the restriction remains a
minimizer. The resulting objective is the local TV problem
\eqref{eq:model} with the exterior values stated in the lemma,
after translating $Q$ to $V_m$. The squared-loss term is strictly
convex, so the minimizer is unique. Unused exterior corner vertices
do not enter the objective and may be assigned arbitrarily.
\end{proof}

\begin{proof}[Proof of Lemma~\ref{lem:rectangle-elementwise}]
Fix a vertex $p$, and abbreviate $m=m_p$, $r=r_p$, $s=s_p$, and
$d=d_p$. If $m\le3$, coordinatewise optimality and the singleton
sets in the definition of $\cutnorm{\beps}$ give
\[
 |\widehat\theta_p^\lambda-\theta_p|
 \le2\lambda+|\varepsilon_p|
 \le2\lambda+2\cutnorm{\beps}
 \le4\lambda
 \le12\frac{\lambda}{m},
\]
which proves the claim in this case.

Suppose $m\ge4$. Interchange the coordinates if needed so that $r$
is the distance in the short direction of $R(p)$; the claimed bound
is unchanged by this interchange. Put $h=\lceil m/2\rceil$ and
$s_* =\min\{s,h\}$. There is an $m\times m$ square $Q\subseteq R(p)$
containing $p$ whose local distances are $r_Q(p)=r$ and
$s_Q(p)=s_*$. To choose its interval in the long direction, place
that interval against the nearer side when $s\le h$, and otherwise
place $p$ in a middle position. Since $r\le h$, we also have
$d_Q(p)=\min\{r,s_*\}=d$.

The directional perimeters of a set $S\subseteq Q$ are the same
whether computed in $Q$ or in the full grid, so
$\mathcal D_Q(\beps|_Q)\le\cutnorm{\beps}$.
Lemma~\ref{lem:patch-restriction} and translation equivariance
identify the centered restriction on $Q$ with a local estimator
having observations $\beps|_Q$.
Theorem~\ref{thm:elementwise}, with $L_m=\log(em)$, therefore gives
\[
 |\widehat\theta_p^\lambda-\theta_p|
 \le\frac{4\lambda}{m}
 +\min\left\{
 2\cutnorm{\beps},\,
 256L_m\left(
 \frac{\cutnorm{\beps}}{\sqrt{rs_*}}
 +\frac{\cutnorm{\beps}^2}{\lambda d}
 \right)\right\}.
\]
If $s\le h$, then $s_*=s$, and $L_m\ge1$ gives the required bound.
If $s>h$, then $d=r$ and $h\ge m/2$. By AM--GM,
\[
 \frac{L_m\cutnorm{\beps}}{\sqrt{rh}}
 \le\frac{\sqrt2\,L_m\cutnorm{\beps}}{\sqrt{mr}}
 \le\frac1{\sqrt2}\left\{
 \frac{\lambda}{m}
 +\frac{L_m^2\cutnorm{\beps}^2}{\lambda r}
 \right\}.
\]
In both cases, the spatial bound is thus at most
\[
 C\left\{
 \frac{\lambda}{m}
 +\frac{L_m\cutnorm{\beps}}{\sqrt{rs}}
 +\frac{L_m^2\cutnorm{\beps}^2}{\lambda d}
 \right\}.
\]
Combining this estimate with the uniform bound
$4\lambda/m+2\cutnorm{\beps}$ proves
\eqref{eq:rectangle-elementwise-det}, after enlarging the universal
constant $C$. Every step is deterministic and uniform over all
exterior values, so the conclusion holds simultaneously over all
vertices and all $\lambda$ satisfying the lemma's hypotheses.
\end{proof}

\section{Proofs of the Gaussian sharpness results}\label{sec:completion}

\subsection{Auxiliary lemmas}\label{subsec:local-risk-lower}

\begin{lemma}[Lower bounds for local risk]\label{lem:local-risk-lower}
For the local estimator~\eqref{eq:model} with independent
$y_v\sim N(0,\sigma^2)$, let $R_{m,\sigma}^+(\lambda)$ and
$R_{m,\sigma}^*(\lambda)$ be the risks defined before
Corollary~\ref{cor:power-match}.
There is a universal $c>0$ such that, for every $m\ge4$,
$\sigma>0$, and $\lambda>0$,
\[
 \begin{aligned}
 R_{m,\sigma}^+(\lambda)
 &\ge\frac{4\lambda^2+\sigma^2}{m^2},\\
 R_{m,\sigma}^+(\lambda)
 &\ge\frac{c\sigma^2}{m\max\{1,(\lambda/\sigma)^2\}}.
 \end{aligned}
\]
Both lower bounds also hold for $R_{m,\sigma}^*(\lambda)$, since
$R_{m,\sigma}^*(\lambda)\ge R_{m,\sigma}^+(\lambda)$.
\end{lemma}

\begin{lemma}[Effect of the two end columns]\label{lem:strip-end-comparison}
Let $H,\ell\ge1$ be integers, let $\lambda>0$, and put
$w=2\ell+1$ and $Q=[H]\times[w]$. Write
\[
 \TV_Q(\bz)=\sum_{\substack{\{a,b\}\subseteq Q\\a\sim b}}
                 |z_a-z_b|,
\]
where each interior edge is counted once.
For observations $\by\in\R^Q$, choose
$M>\norm{\by}_\infty+2\lambda$ and define
\[
 J_M(\bz;\by)=\norm{\bz-\by}_2^2+\lambda\TV_Q(\bz)
 +\lambda\sum_{j=1}^w
       \bigl(|z_{1j}-M|+|z_{Hj}+M|\bigr).
\]
Let
\begin{align*}
 \bv(\by)&=\arg\min_{\bz\in\R^Q}J_M(\bz;\by),\\
 \bu(\by)&=\arg\min_{\bz\in\R^Q}
 \left\{J_M(\bz;\by)+\lambda\sum_{i=1}^H
       \bigl(|z_{i1}+M|+|z_{iw}+M|\bigr)\right\}.
\end{align*}
These are the stabilized solutions of Lemma~\ref{lem:comparison},
so their values are independent of $M$ in the displayed range.
Both solutions use $+M$ above $Q$ and $-M$ below it.
The left and right sides are free for $\bv$, whereas $\bu$
uses $-M$ on those sides. Then, for every row $i$,
\begin{equation}\label{eq:strip-end-comparison}
 0\le v_{i,\ell+1}-u_{i,\ell+1}
 \le\frac{18\lambda\log(eH)}{w}.
\end{equation}
\end{lemma}

\subsection{Proof of Theorem~\ref{thm:elementwise-risk-lower}}
\label{subsec:elementwise-lower-proof}

\begin{proof}[Proof of Theorem~\ref{thm:elementwise-risk-lower}]
Write $d_p=\min(r_p,s_p)$ for this proof.
For any vertex $p$, Lemma~\ref{lem:comparison} and dominated
convergence give
\begin{equation}\label{eq:elementwise-extremal-lower}
 \begin{aligned}
 \E(x_{+,p})_+^2
 &\le\E x_{+,p}^2
 =\lim_{M\to\infty}
 \E|\widehat x_p^{\lambda,M\one_{F_m}}|^2\\
 &\le\sup_{\bb\in\R^{F_m}}\E|\widehat x_p^{\lambda,\bb}|^2.
 \end{aligned}
\end{equation}
The exterior values in this limit are deterministic, and
\eqref{eq:crude} supplies the integrable bound
$(|y_p|+2\lambda)^2$. We first bound the positive part of $x_{+,p}$
for the corner and side terms for every $\lambda>0$, then add the
bias bound under the stronger assumption on $\lambda$.

\emph{A Gaussian mean calculation.}
We will use the following calculation for two stabilized
auxiliary TV solutions on rectangles. Let $\mathcal T(\by_Q)$
be such a solution, and suppose its value at the center $c$ of $Q$
satisfies $\E\mathcal T_c(\by_Q)=0$.
Translation equivariance of the stabilized map in
Lemma~\ref{lem:comparison} gives
\[
 \mathcal T_c(\by_Q)=\overline y_Q+W,\qquad
 W=\mathcal T_c(\by_Q-\overline y_Q\one_Q).
\]
The Gaussian mean $\overline y_Q\sim N(0,\sigma^2/|Q|)$ is
independent of the centered observations, hence of $W$, and
$\E W=0$. The coordinate bound in Lemma~\ref{lem:comparison}
ensures the required moments are finite.
For every deterministic $b_0\ge0$, conditional Jensen's inequality
gives
\[
 \E\bigl(\mathcal T_c(\by_Q)-b_0\bigr)_+^2
 \ge\E(\overline y_Q-b_0)_+^2
 =\frac{\sigma^2}{|Q|}
   \E\left(Z-\frac{b_0\sqrt{|Q|}}{\sigma}\right)_+^2,
 \qquad Z\sim N(0,1).
\]
The center-mean condition will follow from the boundary
symmetries specified below.

\emph{Corner term at a fixed vertex.}
The next argument holds for every $\lambda>0$.
Reflect the square so that $p=(r,s)$, where $r=r_p$ and $s=s_p$,
and let
\[
 Q=[1,2r-1]\times[1,2s-1],\qquad A=|Q|\le4rs.
\]
Then $p$ is the center of $Q$. Let $\bu(\by_Q)$ be the
stabilized solution of the ordinary TV problem on $Q$ with exterior values $+M$
along its top and left sides and $-M$ along its bottom and right
sides, as in Lemma~\ref{lem:comparison}.
Choose $M>\norm{\by}_\infty+2\lambda$ so that the solution on the full square with exterior values $M\one_{F_m}$ is $\bx_+$.
Its restriction to $Q$ has the same top and left exterior values
as the auxiliary solution, and exterior values at least $-M$ on the
other two sides. Boundary comparison therefore gives
$x_{+,p}\ge u_p$.
The argument includes rectangles with one row or one column;
the exterior edges on opposite sides are counted separately.

Rotation of $Q$ through $180$ degrees followed by sign reversal
preserves ordinary TV and the specified exterior configuration.
It fixes $p$ and changes the sign of the estimated value there.
Since the Gaussian observations have the same distribution after
this transformation, $\E u_p=0$. Applying the Gaussian mean
calculation with $\mathcal T=\bu$ and $b_0=0$ yields
\begin{equation}\label{eq:elementwise-corner-lower}
 \E(x_{+,p})_+^2\ge\E(u_p)_+^2
 \ge\frac{\sigma^2}{2A}\ge\frac{\sigma^2}{8rs},
 \qquad \lambda>0.
\end{equation}

\emph{Side term at a fixed vertex.}
Let $\lambda>0$ be arbitrary. Reflect and interchange the
coordinates so that $p=(r,s)$, where
$r=d_p$ and $s=\max\{r_p,s_p\}$.
Put $H=2r-1$ and
\[
 a=324\max\{1,(\lambda/\sigma)^2\}H\log^2(eH),
\]
and let $w$ be the smallest odd integer at least $a$.
Then $H\le m$, $w\ge3$, and $a\le w\le a+2\le2a$.

Suppose first that $w\le2s-1$. Let $Q$ be the rectangle of
height $H$ attached to the top side of the square and of width
$w$ centered horizontally at $p$. After translating its column
indices, $p$ is the center $c=((H+1)/2,(w+1)/2)$ of $Q$.
Give $Q$ exterior values $+M$ above and $-M$ below and on its
left and right sides. Its stabilized solution of the ordinary TV problem is $\bu$
from Lemma~\ref{lem:strip-end-comparison}.
The same boundary comparison as in the corner argument gives
$x_{+,p}\ge u_c$.

The reference solution $\bv$ from Lemma~\ref{lem:strip-end-comparison}
has the same specified top and bottom values and free left and
right sides. Reversing the row order and all signs,
$z_{ij}\mapsto-z_{H+1-i,j}$, preserves ordinary TV and these
boundary conditions and changes the sign of the estimated value at
$c$. The Gaussian observations have the same distribution after
this transformation, so $\E v_c=0$.
Lemma~\ref{lem:strip-end-comparison} and the Gaussian mean
calculation with $\mathcal T=\bv$ and
$b_0=18\lambda\log(eH)/w$ give
\begin{align}\label{eq:elementwise-gaussian-strip}
 \E(x_{+,p})_+^2
 &\ge\E\left(v_c-\frac{18\lambda\log(eH)}w\right)_+^2\notag\\
 &\ge\frac{\sigma^2}{Hw}
   \E\left(Z-\frac{18\lambda\sqrt H\log(eH)}
                       {\sigma\sqrt w}\right)_+^2,
 \qquad Z\sim N(0,1).
\end{align}
By the definition of $a$, the standardized shift is at most
$(\lambda/\sigma)/\sqrt{\max\{1,(\lambda/\sigma)^2\}}\le1$.
Using $w\le2a$, $H\le2r$, and $\log(eH)\le2\log(er)$ gives
\[
 \E(x_{+,p})_+^2
 \ge\frac{\sigma^2}{Hw}\E(Z-1)_+^2
 \ge\frac{c\sigma^2}
 {r^2\log^2(er)\max\{1,(\lambda/\sigma)^2\}}.
\]

If $w>2s-1$, minimality of the odd integer $w$ gives
$2s-1<a$, hence $s<a$. Consequently,
\[
 s\le C\max\{1,(\lambda/\sigma)^2\}\,r\log^2(er),
\]
and the corner bound~\eqref{eq:elementwise-corner-lower} gives
the same conclusion. We have therefore proved, at every vertex,
\begin{equation}\label{eq:elementwise-side-lower}
 \E(x_{+,p})_+^2
 \ge\frac{c\sigma^2}
 {r^2\log^2(er)\max\{1,(\lambda/\sigma)^2\}},
 \qquad \lambda>0.
\end{equation}
In particular, at every boundary vertex $r=1$, this is
$c\sigma^2/\max\{1,(\lambda/\sigma)^2\}$.

The corner and side estimates above hold for the same positive
extremal solution at every vertex and for every $\lambda>0$.
We now establish the bias estimate in the theorem's stated range.

\emph{Bias at every vertex.}
Let $t=\min_v x_{+,v}$ and $S=\{v:x_{+,v}=t\}$.
For a sufficiently large positive exterior value, every edge
leaving $S$ leads to a strictly larger value. Summing the
optimality conditions over $S$ gives
\[
 t|S|=\sum_{v\in S}y_v
       +\frac{\lambda}{2}\{P_1(S)+P_2(S)\}
 \ge\frac{\lambda-\cutnorm{\by}}{2}
       \{P_1(S)+P_2(S)\}.
\]
The inequality uses the definition of $\cutnorm{\by}$ and
$\sqrt{P_1(S)P_2(S)}\le\{P_1(S)+P_2(S)\}/2$.
By~\eqref{eq:area-perimeter},
$P_1(S)+P_2(S)\ge4\sqrt{|S|}\ge4|S|/m$.
Thus, when $\lambda\ge\cutnorm{\by}$,
\begin{equation}\label{eq:elementwise-bias-lower}
 \min_v x_{+,v}\ge\frac{2\{\lambda-\cutnorm{\by}\}}{m}.
\end{equation}
For $\lambda\ge C_0\sigma L_m^{3/4}$ with $C_0$ sufficiently
large, Theorem~\ref{thm:gaussian-D} gives
$\PP\{\cutnorm{\by}\le\lambda/2\}\ge1/2$.
Consequently, $\E(x_{+,p})_+^2\ge\lambda^2/(2m^2)$ at every
vertex.

Enlarge $C_0$ if necessary so that
$\lambda\ge C_0\sigma L_m^{3/4}$ implies $\lambda\ge\sigma$ and
Corollary~\ref{cor:elementwise-risk} also applies.
Then~\eqref{eq:elementwise-side-lower} has the form stated in
\eqref{eq:local-elementwise-risk-lower}, with no logarithmic loss when
$d_p=1$. The bias, corner, and side bounds concern the same positive
extremal solution, so their maximum controls their sum up to a
universal factor. This proves~\eqref{eq:local-elementwise-risk-lower};
the worst-case consequence follows from
\eqref{eq:elementwise-extremal-lower}.
The comparison with the upper bound follows from
Corollary~\ref{cor:elementwise-risk}, using
$1\le\log(ed_p)\le L_m$ and $1\le L_m^{7/2}\le L_m^5$.
\end{proof}

\subsection{Proof of Corollary~\ref{cor:power-match}}\label{subsec:power-match-proof}

\begin{proof}[Proof of Corollary~\ref{cor:power-match}]
Fix $\alpha\in(0,1)$. For all sufficiently large $m$,
$m^\alpha\ge C_0L_m^{3/4}$, so Corollary~\ref{cor:elementwise-risk}
applies with $\lambda=\sigma m^\alpha$. Substituting into~\eqref{eq:mse-risk}
gives
\[
 \frac1{m^2}\E\!\left[
 \sup_{\bb\in\R^{F_m}}
 \|\widehat{\bx}^{\,\sigma m^\alpha,\bb}\|_2^2
 \right]
 \le C\sigma^2\left\{
 m^{-(2-2\alpha)}
 +L_m^5m^{-(1+2\alpha)}
 +L_m^{11/2}m^{-2}
 \right\}.
\]
By the definition of $\gamma(\alpha)$ and $L_m\ge1$, each of the first
two terms in braces is at most $L_m^5m^{-\gamma(\alpha)}$. Moreover,
$\gamma(\alpha)\le3/2$ and $L_m\le m$, so
\[
 L_m^{11/2}m^{-2}
 =L_m^5m^{-3/2}\left(\frac{L_m}{m}\right)^{1/2}
 \le L_m^5m^{-3/2}
 \le L_m^5m^{-\gamma(\alpha)}.
\]
This proves the upper bound with some universal constant $C$.
Lemma~\ref{lem:local-risk-lower} gives
$R_{m,\sigma}^+(\lambda)\le R_{m,\sigma}^*(\lambda)$.
The next inequality follows by taking the supremum inside the
expectation.

For the lower bound, the first inequality in
Lemma~\ref{lem:local-risk-lower} gives
\[
 R_{m,\sigma}^+(\sigma m^\alpha)\ge4\sigma^2m^{-(2-2\alpha)}.
\]
The term $4\lambda^2/m^2$ is the \emph{squared boundary-induced bias term}:
the positive extremal solution satisfies
$\overline x_+=\overline y+2\lambda/m$, so sufficiently large positive
exterior values induce a mean shift of $2\lambda/m$. Its square
contributes to the risk, and the shift persists when $\by=0$.
If $1/4\le\alpha<1$, then $\gamma(\alpha)=2-2\alpha$, so this
lower bound proves the claim.

If $0<\alpha<1/4$, the second inequality in
Lemma~\ref{lem:local-risk-lower} and $m^\alpha\ge1$ yield
\[
 R_{m,\sigma}^+(\sigma m^\alpha)
 \ge\frac{c\sigma^2}{m\max\{1,m^{2\alpha}\}}
 =c\sigma^2m^{-(1+2\alpha)}
 =c\sigma^2m^{-\gamma(\alpha)}.
\]
Choosing $m_0(\alpha)$ sufficiently large proves the corollary.
\end{proof}

\subsection{Proof of Corollary~\ref{thm:global-risk-sharpness}}
\label{subsec:global-risk-lower-proof}

\begin{proof}[Proof of Corollary~\ref{thm:global-risk-sharpness}]
We realize the positive extremal solution on a fixed fraction of
the squares in a global checkerboard signal, and then apply the
local lower-bound construction.

For $a,b\in[q]$, write
$Q_{ab}=[(a-1)m+1,am]\times[(b-1)m+1,bm]$.
Use $\btheta^{(M)}$ from the checkerboard example in
Section~\ref{subsec:gaussian-rectangle-results}.
Each zero square with $2\le a,b\le q-1$ is surrounded on all
four sides by squares with signal value $M$. Their number is
\[
 J=\left\lceil\frac{(q-2)^2}{2}\right\rceil\ge\frac{k}{9}.
\]
For $q=3$, $J=1=k/9$; for $q\ge4$,
$J\ge(q-2)^2/2\ge q^2/8$.

Fix $\lambda>0$ and couple the observations as
$\by^{(M)}=\btheta^{(M)}+\beps$, using the same Gaussian noise
array for all $M$. Denote the resulting global estimator by
$\widehat{\btheta}^{\lambda,M}$.
Since every vertex has at most four adjacent vertices,
coordinatewise optimality gives
\[
 \bigl|\widehat\theta_p^{\lambda,M}-y_p^{(M)}\bigr|
 \le2\lambda.
\]
If $M>2\norm{\beps}_\infty+4\lambda$, every estimated exterior
value adjacent to an interior zero square $Q$ exceeds
$\norm{\beps_Q}_\infty+2\lambda$.
Lemma~\ref{lem:patch-restriction} and the stabilization argument
in Lemma~\ref{lem:comparison} therefore show that
$\widehat{\btheta}^{\lambda,M}_Q=\bx_+(\beps_Q)$, after
translating $Q$ to $[m]^2$.
Thus this equality holds for every noise realization once $M$ is
sufficiently large. On each such zero square,
$|\widehat\theta_p^{\lambda,M}|\le|\varepsilon_p|+2\lambda$,
so dominated convergence applies to its squared error.

Summing over the $J$ interior zero squares and taking the limit
through deterministic values of $M$ gives
\begin{align*}
 \liminf_{M\to\infty}\frac1N
 \E\norm{\widehat{\btheta}^{\lambda,M}-\btheta^{(M)}}_2^2
 &\ge
 \lim_{M\to\infty}\frac1N
 \sum_{\substack{2\le a,b\le q-1\\a+b\text{ even}}}
 \E\norm{\widehat{\btheta}^{\lambda,M}_{Q_{ab}}}_2^2\\
 &=\frac{J}{k}R_{m,\sigma}^+(\lambda)
 \ge\frac19R_{m,\sigma}^+(\lambda).
\end{align*}

Lemma~\ref{lem:local-risk-lower} gives, for every $\lambda>0$,
\[
 R_{m,\sigma}^+(\lambda)\ge c\left\{
 \frac{\lambda^2}{m^2}
 +\frac{\sigma^2}{m\max\{1,(\lambda/\sigma)^2\}}
 \right\},
\]
which proves~\eqref{eq:global-risk-lower}. Since
$\btheta^{(M)}\in\Theta_{q,m}$ for every $M$, the same bound holds
for the supremum over $\Theta_{q,m}$.

If $\lambda\le\sigma$, the second term in
\eqref{eq:global-risk-lower} is $\sigma^2/m\ge\sigma^2m^{-3/2}$.
If $\lambda\ge\sigma$, the arithmetic--geometric mean inequality
gives
\[
 \frac{\lambda^2}{m^2}+\frac{\sigma^4}{m\lambda^2}
 \ge\frac{2\sigma^2}{m^{3/2}}.
\]
Since $m^{-3/2}=(k/N)^{3/4}$, taking the infimum over
deterministic $\lambda>0$ proves~\eqref{eq:global-risk-optimal}.

Finally, fix a vertex $p$ of any interior partition square $Q$
in the fixed checkerboard example. If $Q$ has signal value zero,
the centered restriction is $\bx_+(\beps_Q)$ for sufficiently large $M$.
If $Q$ has signal value $M$, centering by $M$ makes every adjacent
exterior value less than $-\norm{\beps_Q}_\infty-2\lambda$ under
the same condition $M>2\norm{\beps}_\infty+4\lambda$.
The centered restriction is then $\bx_-(\beps_Q)$.
By~\eqref{eq:extremal-sign} and Gaussian sign symmetry,
$\E x_{-,p}^2=\E x_{+,p}^2$.
In both cases, the bound
$|\widehat\theta_p^{\lambda,M}-\theta_p^{(M)}|
\le |\varepsilon_p|+2\lambda$
allows dominated convergence, giving
\[
 \lim_{M\to\infty}\E_{\btheta^{(M)}}
 |\widehat\theta_p^\lambda-\theta_p^{(M)}|^2
 =\E x_{+,p}^2.
\]
The supremum of the elementwise risk over $\Theta_{q,m}$ is at
least this limit.
Thus the corner and side bounds in
\eqref{eq:elementwise-corner-lower} and~\eqref{eq:elementwise-side-lower},
proved for the positive extremal solution, transfer to every vertex
of $Q$ for every $\lambda>0$, with distances measured within $Q$. When
$\lambda\ge C_0\sigma L_m^{3/4}$, the additional elementwise
bias bound in~\eqref{eq:local-elementwise-risk-lower} also transfers.
\end{proof}

\subsection{Proof of Lemma~\ref{lem:local-risk-lower}}\label{subsec:local-risk-lower-proof}

\begin{proof}[Proof of Lemma~\ref{lem:local-risk-lower}]
The mean identity~\eqref{eq:test-set-balance} gives
$\overline x_+=\overline y+2\lambda/m$.
Since $\Var(\overline y)=\sigma^2/m^2$ and
$m^{-2}\norm{\bx_+}_2^2\ge\overline x_+^2$, this proves the
first lower bound for $R_{m,\sigma}^+(\lambda)$ for every $\lambda>0$.
For the second, sum~\eqref{eq:elementwise-side-lower}, which was
proved for every $\lambda>0$, over the top row:
\[
 R_{m,\sigma}^+(\lambda)
 \ge\frac1{m^2}\sum_{j=1}^m\E(x_{+,(1,j)})_+^2
 \ge\frac{c\sigma^2}{m\max\{1,(\lambda/\sigma)^2\}}.
\]
Finally, taking deterministic exterior values $M\one_{F_m}$ and
letting $M\to\infty$, Lemma~\ref{lem:comparison}, \eqref{eq:crude},
and dominated convergence give
$R_{m,\sigma}^*(\lambda)\ge R_{m,\sigma}^+(\lambda)$.
\end{proof}

\subsection{Proof of Lemma~\ref{lem:strip-end-comparison}}
\label{subsec:strip-end-comparison-proof}

\begin{proof}[Proof of Lemma~\ref{lem:strip-end-comparison}]
Let $J^u,J^v$ denote the objectives defining $\bu,\bv$, put
$\boldsymbol{\Delta}=\bv-\bu$ and $Q_{\rm end}=[H]\times\{1,w\}$,
and let $e(C,D)$ count the edges joining disjoint sets $C,D\subseteq Q$.
Fix $t\in\R$, set $S=\{a\in Q:\Delta_a>t\}$, and take a nonempty
$C\subseteq S$. As in the proof of Lemma~\ref{lem:geo-binary}, perturb the two solutions
in opposite directions. For sufficiently small $\varepsilon>0$,
\[
 0\le J^u(\bu+\varepsilon\one_C)-J^u(\bu)
       +J^v(\bv-\varepsilon\one_C)-J^v(\bv).
\]
The two squared losses contribute
$-2\varepsilon\sum_{a\in C}\Delta_a+2\varepsilon^2|C|$.
Lemma~\ref{lem:comparison} bounds both solutions in
$[-\norm{\by}_\infty-2\lambda,\norm{\by}_\infty+2\lambda]$,
so the perturbed values remain in $(-M,M)$ for sufficiently small
$\varepsilon$. The common top and bottom penalties therefore cancel,
and the end penalties in $J^u$ contribute
$\lambda\varepsilon|C\cap Q_{\rm end}|$.
Each edge from $C$ to $S\setminus C$ contributes at most
$2\lambda\varepsilon$. For an edge from $a\in C$ to $b\in S^c$,
the combined contribution is
\[
\begin{aligned}
 &\lambda\bigl(|u_a-u_b+\varepsilon|-|u_a-u_b|
              +|v_a-v_b-\varepsilon|-|v_a-v_b|\bigr)\\
 &\qquad=\lambda\bigl(g_\varepsilon(u_a-u_b)
                    -g_\varepsilon(v_a-v_b-\varepsilon)\bigr)\le0,
\end{aligned}
\]
where $g_\varepsilon(x)=|x+\varepsilon|-|x|$ is nondecreasing and
$u_a-u_b=v_a-v_b-(\Delta_a-\Delta_b)<v_a-v_b-\varepsilon$
when $\varepsilon<\Delta_a-\Delta_b$.
Dividing by $2\varepsilon$ and letting $\varepsilon\downarrow0$ gives
\begin{equation}\label{eq:strip-cut-comparison}
 \sum_{a\in C}\Delta_a
 \le\lambda e(C,S\setminus C)
       +\frac\lambda2|C\cap Q_{\rm end}|.
\end{equation}

For the lower bound, let
\[
C=\{a\in Q:\Delta_a<0\}.
\]
Use the opposite perturbations
\[
u\mapsto u-\varepsilon\mathbf 1_C,
\qquad
v\mapsto v+\varepsilon\mathbf 1_C.
\]
The same comparison as above gives
\[
\sum_{a\in C}\Delta_a
\ge \frac{\lambda}{2}|C\cap Q_{\rm end}|
\ge0.
\]
Since every summand on the left is negative, \(C\) must be empty. Thus \(\Delta_a\ge0\) for all \(a\in Q\).

For the upper bound, set
\[
S=\{a\in Q:\Delta_a>\lambda/2\}.
\]
Applying~\eqref{eq:strip-cut-comparison} with \(t=\lambda/2\) and \(C=S\) gives
\[
\frac{\lambda}{2}|S|
<
\sum_{a\in S}\Delta_a
\le
\frac{\lambda}{2}|S\cap Q_{\rm end}|
\le
\frac{\lambda}{2}|S|,
\]
unless \(S=\varnothing\). Hence
\begin{equation}\label{eq:strip-delta-range}
0\le\Delta_a\le\lambda/2
\qquad(a\in Q).
\end{equation}

Fix a row \(i\) and suppose \(0<t<\Delta_{i,\ell+1}\). Define
\[
S=\{a\in Q:\Delta_a>t\},\qquad
B_k=[H]\times[\ell+1-k,\ell+1+k],\qquad
N_k=|S\cap B_k|,
\]
for \(0\le k\le \ell\). Since \((i,\ell+1)\in S\), we have \(N_0\ge1\).

For \(k<\ell\), apply~\eqref{eq:strip-cut-comparison} to
\[
C=S\cap B_k.
\]
The strip \(B_k\) does not contain either end column, so \(C\cap Q_{\rm end}=\varnothing\). Moreover, since \(\Delta_a>t\) for every \(a\in C\),
\[
tN_k<\sum_{a\in C}\Delta_a
\le \lambda e(C,S\setminus C).
\]
Any edge from \(C\) to \(S\setminus C\) must leave \(B_k\) through one of the two adjacent columns \(B_{k+1}\setminus B_k\). Each vertex in these new columns has at most one neighbor in \(B_k\), and therefore
\[
e(C,S\setminus C)
\le |S\cap(B_{k+1}\setminus B_k)|
=N_{k+1}-N_k.
\]
Consequently,
\[
tN_k<\lambda(N_{k+1}-N_k),
\]
or
\begin{equation}\label{eq:strip-count-growth}
N_{k+1}>(1+t/\lambda)N_k.
\end{equation}
Since \(t>0\), this implies \(N_{k+1}>N_k\). The \(N_k\) are integers and \(N_0\ge1\), so induction gives
\[
N_k\ge k+1.
\]
With $k_0=\lfloor \ell/2\rfloor$, we have
$N_{k_0}\ge \ell/2$, $\ell-k_0\ge \ell/2$, and $N_\ell\le3H\ell$.
Consequently,
\[
 (1+t/\lambda)^{\ell/2}\le\frac{N_\ell}{N_{k_0}}\le6H, \,\,\text{i.e.,}\,\, \frac{\ell}{2}\log\left(1+\frac{t}{\lambda}\right)\le \log(6H).
\]
Since \(0<t<\Delta_{i,\ell+1}\le\lambda/2\) by~\eqref{eq:strip-delta-range}, we have
\(\log(1+t/\lambda)\ge 2t/(3\lambda)\), and hence
\[
t\le \frac{3\lambda\log(6H)}{\ell}.
\]
Because \(w=2\ell+1\le3\ell\) and \(\log(6H)\le2\log(eH)\),
\[
t\le \frac{18\lambda\log(eH)}{w}.
\]
This holds for every \(0<t<\Delta_{i,\ell+1}\); letting \(t\uparrow\Delta_{i,\ell+1}\) proves the upper bound. The lower bound follows from~\eqref{eq:strip-delta-range}.
\end{proof}

\end{document}